\documentclass[11pt]{article}

\usepackage{tikz}
\usetikzlibrary{calc}
\usetikzlibrary{matrix}
\usepackage{amsmath,amssymb,amsfonts} 
\usepackage{verbatim}
\usepackage{color}                    
\usepackage{hyperref}                 

\usepackage{graphicx}

\def \p{\partial}
\newcommand{\V}[1]{\mbox{\boldmath $ #1 $}}
\newcommand{\bey}{\begin{eqnarray}}
\newcommand{\eey}{\end{eqnarray}}
\newcommand{\nn}{\nonumber}
\newcommand{\beq}{\begin{equation}}
\newcommand{\eeq}{\end{equation}}
\newtheorem{thm}{\hspace{6mm}Theorem}[section]

\newtheorem{lem}{\hspace{6mm}Lemma}[section]
\newtheorem{co}{\hspace{6mm}Corollary}[section]

\newtheorem{exam}{\hspace{6mm}Example}[section]
\newtheorem{rem}{\hspace{6mm}Remark}[section]
\newcommand{\proofend}{\mbox{ }\hfill \raisebox{.4ex}{\framebox[1ex]{}}}
\vspace{5mm}

\begin{document}

\date{}
\title{An anisotropic mesh adaptation method for the finite element solution of heterogeneous anisotropic diffusion problems}
\author{Xianping Li \thanks{
Department of Mathematics, the University of Kansas, Lawrence, KS 66045,
U.S.A. ({\tt lxp@math.ku.edu})}
\and Weizhang Huang \thanks{
Department of Mathematics, the University of Kansas, Lawrence, KS 66045,
U.S.A. ({\tt huang@math.ku.edu})}
}
\maketitle

\vspace{10pt}

\begin{abstract}
Heterogeneous anisotropic diffusion problems arise in the various areas of science and engineering including
plasma physics, petroleum engineering, and image processing. Standard numerical methods can produce
spurious oscillations when they are used to solve those problems. A common approach to avoid this difficulty
is to design a proper numerical scheme and/or a proper mesh so that the numerical solution validates
the discrete counterpart (DMP) of the maximum principle satisfied by the continuous solution.
A well known mesh condition for the DMP satisfaction
by the linear finite element solution of isotropic diffusion problems
is the non-obtuse angle condition that requires the dihedral angles of
mesh elements to be non-obtuse. In this paper, a generalization of the condition, the so-called
anisotropic non-obtuse angle condition, is developed for the finite element solution of heterogeneous anisotropic
diffusion problems.  The new condition is essentially the same as
the existing one except that the dihedral angles are now measured in a metric
depending on the diffusion matrix
of the underlying problem. Several variants of the new condition are obtained.
Based on one of them, two metric tensors for use in anisotropic
mesh generation are developed to account for DMP satisfaction and the combination of
DMP satisfaction and mesh adaptivity. Numerical examples are given to demonstrate
the features of the linear finite element method for anisotropic meshes generated with  the metric tensors.
\end{abstract}



\vspace{10pt}

\section{Introduction}
\label{Sec-intro}

We are concerned with the numerical solution of the diffusion equation
\beq
\label{bvp-pde}
 - \nabla \cdot (\mathbb{D} \, \nabla u)  =  f, \qquad \mbox{ in } \quad \Omega
\eeq
subject to the Dirichlet boundary condition
\beq
\label{bvp-bc}
u =  g, \qquad \mbox{ on } \quad \partial \Omega
\eeq
where $\Omega \subset \mathbb{R}^d \; (d=1,2, \text{ or } 3)$ is the physical domain,
$f$ and $g$ are given functions, and $\mathbb{D}=\mathbb{D}(\V{x})$ is the diffusion matrix
assumed to be symmetric and strictly positive definite on $\Omega$.
The boundary value problem (BVP) (\ref{bvp-pde}) and (\ref{bvp-bc})
becomes a heterogeneous anisotropic diffusion problem when $\mathbb{D}$ changes from place to place
(heterogeneous) and its eigenvalues are not all equal (anisotropic)
at least on a portion of $\Omega$.
When a standard numerical method, such as a finite element, a finite difference, or a finite volume method,
is used to solve this problem,
spurious oscillations can occur in the computed solution.
A challenge is then to design a proper numerical scheme and/or
a proper mesh so that the computed solution is free of spurious oscillations.
In some applications such as plasma physics, it is further desired
that the mesh be aligned with the fast diffusion direction so that no excessive numerical dissipation
is introduced in slow diffusion directions.
Moreover, mesh adaptation is often necessary for improving
computational efficiency and accuracy when
the physical solution and/or the diffusion matrix have sharp jumps.

Anisotropic diffusion problems arise in the various areas of science and engineering including
plasma physics in fusion experiments and astrophysics
\cite{GL09, GLT07, GYK05, NW00, SH07, Sti92},
petroleum reservoir simulation \cite{ABBM98a,ABBM98b,CSW95,EAK01,MD06},
and image processing \cite{CS00, CSV03, KM09, MS89, PM90, Wei98}.
In plasma physics, magnetized plasmas are constrained to move primarily along magnetic field lines.
Their heat conductivity in the direction parallel to the magnetic field is much higher than
those perpendicular to it, and the ratio of the conduction coefficients can easily exceed $10^{10}$ in fusion experiments.
The numerical simulation of the heat conduction of plasmas must not only produce
a physically meaningful temperature distribution but also avoid excessive numerical dissipation in the directions
perpendicular to the magnetic field.
In petroleum engineering, fluids such as water, crude oil, and natural gas are stored in reservoir rocks filled
by interconnected networks of pores.
The diffusion and flow of those fluids depend crucially on the rocks' permeability which
changes with location and flow direction and has much large values in horizontal directions
than in the vertical direction.
Finally, PDE-based anisotropic diffusion filters have been successfully used for shape recognition
and edge detection in image processing.

The BVP (\ref{bvp-pde}) and (\ref{bvp-bc}) is a representative example of anisotropic diffusion problems arising
in those areas. As typical for diffusion problems, it satisfies the maximum principle
\beq
\max_{\V{x} \in \Omega \cup \p \Omega} u(\V{x}) \le \max\{ 0, \max_{\V{s}\in \p \Omega} g(\V{s})\}
\label{MP}
\eeq
provided that $f (\V{x}) \le 0$ holds for all $\V{x} \in \Omega$. The BVP has been studied extensively in the past,
and a major effort has been made to avoid spurious oscillations in the numerical solution.
A common strategy is to develop numerical schemes satisfying the discrete counterpart of (\ref{MP}) --
the so-called discrete maximum principle (DMP), which are known to produce numerical solutions free of
spurious oscillations \cite{Cia70,Var66}.
The studies can be traced back to early works
by Varga \cite{Var66}, Ciarlet \cite{Cia70}, Ciarlet and Raviart \cite{CR73}, and Stoyan \cite{Sto82,Sto86}
where a number of sufficient conditions in a general and abstract setting are obtained for a class of linear elliptic
partial differential equations (PDEs). For example, denote by
$
A \V{u} = \V{f}
$
the linear algebraic system resulting
from the application of a numerical scheme to a linear elliptic PDE supplemented with a Dirichlet boundary
condition, where $A$ is the $n\times n$ stiffness matrix, $\V{u}$ is the unknown vector, and $\V{f}$
the right-hand-side vector.
Then, a sufficient condition is given as follows.

\begin{lem}
\label{lem1.1}
(\cite{Sto86}) If the stiffness matrix $A$ satisfies
\bey
\mbox{(a)} && \mbox{that $A$ is monotone with $A^{(-)}$ being either nonsingular, or singular and irreducible; and}\qquad
\label{cond-a}
\\
\mbox{(b)} && \mbox{that $A^{(-)} \V{e}^{(n)} \ge 0$,}
\label{cond-b}
\eey
then the numerical scheme satisfies DMP.
\end{lem}

Here, matrix $A$ is said to be monotone if $A$ is nonsingular and $A^{-1} \ge 0$
(i.e., all entries of $A^{-1}$ are non-negative),
and $A^{(-)}$ and $\V{e}^{(n)}$ are defined as
\beq
a^{(-)}_{ij} = \begin{cases} a_{ii}, & \mbox{ for } i = j\\ a_{ij}, & \mbox{ for } i \ne j, \; a_{ij} \le 0 \\
0, & \mbox{ for } i \ne j, \; a_{ij} > 0\end{cases} ,
\quad \V{e}^{(n)} = \left [\begin{array}{c} 1 \\ \vdots \\ 1\end{array} \right ].
\eeq
Note that condition (\ref{cond-b}) is equivalent to that $A^{(-)}$ has nonnegative
row sums. Moreover,  $A = A^{(-)}$ and the condition (\ref{cond-a}) holds when $A$ is an $M$-matrix \cite{Var62}.
From Lemma \ref{lem1.1} we have the following lemma.

\begin{lem}
\label{lem1.2}
If the stiffness matrix $A$ is an $M$-matrix and has nonnegative row sums,
then the numerical scheme satisfies DMP.
\end{lem}

Numerical schemes satisfying DMP have been developed
along the line of those sufficient conditions by either designing a proper discretization for
the underlying PDE or employing a suitable mesh.
To date most success has been made for the isotropic diffusion case where $\mathbb{D}$
is in the scalar matrix form, $\mathbb{D} = a(\V{x}) I$, with $a(\V{x})$ being a scalar function;
e.g., see \cite{BKK08,BE04,CR73,KK05,KK06,KKK07,KL95,Let92,SF73,XZ99}.
In particular, it is shown in \cite{BKK08,CR73} that the linear finite element method (FEM)
satisfies DMP when the mesh is simplicial and satisfies the so-called non-obtuse angle
condition requiring that the dihedral angles of all mesh elements be non-obtuse.
In two dimensions this condition can be replaced by a weaker condition (the Delaunay condition) that
the sum of any pair of angles opposite a common edge is less than or equal to $\pi$ \cite{Let92,SF73}.
Similar mesh conditions are developed in \cite{KK05,KK06,KKK07,KL95}
for elliptic problems with a nonlinear diffusion coefficient in the form
$\mathbb{D} = a(\V{x}, u, \nabla u) I$ and with mixed boundary conditions.
Burman and Ern \cite{BE04} propose a nonlinear stabilized Galerkin approximation
for the Laplace operator and prove that it satisfies DMP on arbitrary meshes and for
arbitrary space dimension without resorting to the non-obtuse angle condition.

On the other hand, the anisotropic diffusion case is more difficult and
only limited success has been made
\cite{CSW95,DDS04,GL09,GLT07,GYK05,KSS09,LePot05,LePot09,LePot09b,LSS07, LSSV07,LS08,MD06,SH07}.
For example, Dr\v{a}g\v{a}nescu et al. \cite{DDS04} show that the non-obtuse angle
condition fails to guarantee DMP satisfaction in the anisotropic diffusion case.
The techniques proposed by Liska and Shashkov \cite{LS08} and Kuzmin et al.
\cite{KSS09} to locally modify (or repair) the underlying numerical scheme, by Sharma and Hammett \cite{SH07}
to employ slope limiters in the discretization of the PDE, by Mlacnik and Durlofsky \cite{MD06} to
optimize the mesh for a multipoint flux approximation (MPFA) finite volume method
(e.g., see \cite{ABBM98a,ABBM98b} for the method),
and by Li et al. \cite{LSS07} to optimize a triangular mesh for the finite element solution,
help reduce spurious oscillations.
A nonlinear, first order finite volume method developed by Le Potier \cite{LePot05,LePot09} and
further improved by Lipnikov et al. \cite{LSSV07} gives rise to a stiffness $M$-matrix
on arbitrary meshes when applied to parabolic PDEs but fails to satisfy DMP
when applied to steady-state elliptic problems.
A first order finite difference method having similar features is proposed by Le Potier \cite{LePot09b}.

In this paper we study the linear finite element solution of BVP
(\ref{bvp-pde}) and (\ref{bvp-bc}) with a general diffusion matrix $\mathbb{D} = \mathbb{D}(\V{x})$.
The objective is threefold. The first is to develop a generalization of the well known non-obtuse angle condition,
the so-called anisotropic non-obtuse angle condition (cf. equation (\ref{g-nonobtuse})),
so that the linear FEM satisfies DMP when the mesh is simplicial and satisfies this condition.
The condition requires that the dihedral angles of all mesh elements,  measured in a metric
depending on $\mathbb{D}$, be non-obtuse. It reduces to the
non-obtuse angle condition for isotropic diffusion matrices. It also reproduces
several existing mesh conditions for homogeneous anisotropic media for which $\mathbb{D}$ is
a full, constant matrix (see Remark \ref{rem2.2}).
The second objective is to derive a metric tensor for use in mesh generation based on the anisotropic non-obtuse
angle condition. This is done by adopting the so-called $M$-uniform mesh approach \cite{Hua05b}
where an anisotropic mesh is generated as an $M$-uniform mesh or a uniform mesh
in the metric specified by a tensor. $M$-uniform meshes generated with the metric tensor
satisfy the anisotropic non-obtuse angle condition
and are aligned with the diffusion matrix $\mathbb{D}$ (cf. \S\ref{Sec-metric}).
The final objective is to combine both mesh adaptivity and DMP satisfaction
in the numerical solution of anisotropic diffusion problems.
An optimal metric tensor (see (\ref{M-DMP+adap})) accounting for both considerations is obtained
by minimizing an interpolation error bound, and advantages of using adaptive, DMP-bound meshes
are demonstrated in numerical examples.
To the authors' best knowledge, this is the first effort that
mesh adaptivity and DMP satisfaction are considered simultaneously
in the numerical solution of anisotropic diffusion problems.

The outline of this paper is as follows. In section \ref{Sec-fem}, the linear finite element solution of (\ref{bvp-pde})
and (\ref{bvp-bc})
is described and the anisotropic non-obtuse angle condition and several variants are derived.
Section \ref{Sec-metric} is devoted to the derivation of the metric tensor based on the anisotropic
non-obtuse angle condition.
In section \ref{Sec-adapt}, the combination of mesh adaptation and DMP satisfaction is addressed,
and an optimal metric tensor is obtained by minimizing an interpolation error bound.
Numerical examples are presented in section \ref{Sec-results}. Finally, section \ref{Sec-con} contains
conclusions and comments.

\section{Anisotropic non-obtuse angle conditions for linear finite element approximation}
\label{Sec-fem}

Consider the linear finite element solution of BVP (\ref{bvp-pde}) and (\ref{bvp-bc}). 
Assume that $\Omega$ is a connected polygon or polyhedron and an affine family of simplicial
triangulations $\{ \mathcal{T}_h \}$ is given thereon.
Let
\[
U_g = \{ v \in H^1(\Omega) \; | \: v|_{\p \Omega} = g\}.
\]
Denote by
$U_g^h \subset U_g$ the linear finite element space associated with mesh $\mathcal{T}_h$.
Then a linear finite element solution $\tilde{u}^h \in U_g^h$ to BVP (\ref{bvp-pde})
and (\ref{bvp-bc}) is defined by
\beq
\label{fem-form}
\int_{\Omega} (\nabla v^h)^T \; \mathbb{D} \nabla \tilde{u}^h d\V{x} =
 \int_{\Omega} f \, v^h d\V{x}, \quad \forall v^h \in U_0^h
\eeq
where $U_0^h = U_g^h$ with $g = 0$.
This equation can be rewritten as
\beq
\label{disc-1}
\sum_{K \in \mathcal{T}_h} \int_{K} (\nabla v^h)^{T} \, \mathbb{D} \,
\nabla \tilde{u}^h d\V{x} =
 \sum_{K \in \mathcal{T}_h} \int_{K} f \, v^h d\V{x}, \quad \forall v^h \in U_0^h.
\eeq
Generally speaking, the integrals in (\ref{disc-1}) cannot be carried out analytically, and numerical quadrature
is needed. We assume that a quadrature rule has been chosen on the reference element $\hat{K}$ for this purpose,
\beq
\label{quadra-rule}
\int_{\hat{K}} v(\xi) d \xi \approx |\hat{K}| \sum_{k=1}^m \hat{w}_k v(\hat{b}_k),
\quad \sum_{k=1}^m \hat{w}_k =1,
\eeq
where $\hat{w}_k$'s are the weights and $\hat{b}_k$'s the quadrature nodes.
A 2D example of such quadrature rules is given by $\hat{w}_k=\frac{1}{3}$ ($k=1,2,3$)
and the barycentric coordinates of the nodes
($\frac{1}{6},\frac{1}{6},\frac{2}{3}$), ($\frac{1}{6},\frac{2}{3},\frac{1}{6}$), and ($\frac{2}{3},\frac{1}{6},\frac{1}{6}$);
and a 3D example is $\hat{w}_i=\frac{1}{4}$ ($ i=1,2,3,4$) and the barycentric coordinates of the nodes
($a,a,a,1-3a$), ($a,a,1-3a,a$), ($a,1-3a,a,a$), and ($1-3a,a,a,a$) with $a=\frac{5-\sqrt{5}}{20}$;
e.g., see \cite{EG04}.

Let $F_K$ be the affine mapping from $\hat{K}$ to $K$ such that $K = F_K(\hat{K})$, and denote
$b_k^K = F_K(\hat{b}_k)$, $k=1,\cdots,m$. Upon applying (\ref{quadra-rule}) to the integrals in (\ref{disc-1})
and changing variables, the finite element approximation problem becomes seeking $u^h \in U_g^h$ such that
\beq
\label{disc-2}
\sum_{K \in \mathcal{T}_h} |K| \sum_{k=1}^m \hat{w}_k
\; (\nabla v^h|_K)^T \; \mathbb{D}(b_k^K) \; \nabla u^h|_K
=  \sum_{K \in \mathcal{T}_h} |K| \sum_{k=1}^m \hat{w}_k f(b_k^K) \; v^h (b_k^K),
\quad \forall v^h \in U_0^h
\eeq
where $\nabla v^h|_K$ and $\nabla u^h|_K$ denote the restriction of $\nabla v^h$
and $\nabla u^h$
on $K$, respectively. Note that we have used in (\ref{disc-2}) the fact that $\nabla v^h|_K$ and
$\nabla u^h|_K$ are constant. Letting
\beq
\label{def-DK}
\mathbb{D}_K = \sum_{k=1}^m \hat{w}_k \mathbb{D}(b_k^K) ,
\eeq
we can rewrite (\ref{disc-2}) into
\beq
\label{disc-3}
\sum_{K \in \mathcal{T}_h} |K| \; (\nabla v^h|_K)^T \; \mathbb{D}_K \; \nabla u^h|_K
=  \sum_{K \in \mathcal{T}_h} |K| \sum_{k=1}^m \hat{w}_k f(b_k^K) \; v^h(b_k^K),
\quad \forall v^h \in U_0^h.
\eeq

We now express (\ref{disc-3}) in a matrix form.
Denote the numbers of the elements, vertices, and interior vertices of $\mathcal{T}_h$
by $N$, $N_v$, and $N_{vi}$, respectively. Assume that the vertices are ordered in such a way that
the first $N_{vi}$ vertices are the interior vertices. Then $U_0^h$ and $u^h$ can be expressed as
\beq
U_0^h = \text{span} \{ \phi_1, \cdots, \phi_{N_{vi}} \}
\eeq
and
\beq
\label{soln-approx}
u^h = \sum_{j=1}^{N_{vi}} u_j \phi_j + \sum_{j=N_{vi}+1}^{N_{v}} u_j \phi_j ,
\eeq
where $\phi_j$ is the linear basis function associated with the $j$-th vertex, $\V{a}_j$.
Note that the boundary condition (\ref{bvp-bc}) can be approximated by
\beq
u_j = g_j \equiv g(\V{a}_j), \quad j = N_{vi}+1, ..., N_v .
\label{fem-bc}
\eeq
Substituting (\ref{soln-approx}) into and taking $v^h = \phi_i$ ($i=1, ..., N_{vi}$) in (\ref{disc-3})
and combining the resulting equations with (\ref{fem-bc}), we obtain the linear algebraic system
\beq
\label{fem-linsys}
A \, \V{u} = \V{f},
\eeq
where
\beq
A = \left [\begin{array}{cc} A_{11} & A_{12} \\ 0 & I \end{array} \right ],
\label{fem-matrix}
\eeq
$I$ is the identity matrix of size $(N_{v} - N_{vi})$, and
\[
\V{u} = (u_1,..., u_{N_{vi}}, u_{N_{vi}+1},..., u_{N_v})^T,
\]
\[
\V{f} = ( f_1, ..., f_{N_{vi}}, g_{N_{vi}+1}, ..., g_{N_v} )^T.
\]
The entries of the stiffness matrix $A$ and the right-hand-side
vector $\V{f}$ are given by
\beq
\label{stiffness-matrix}
a_{ij}= \sum_{K \in \mathcal{T}_h} |K| \;(\nabla \phi_i|_K)^T \; \mathbb{D}_K \; \nabla \phi_j|_K,
\quad i=1, ..., N_{vi},\; j=1, ..., N_{v},
\eeq
\beq
f_i = \sum_{K \in \mathcal{T}_h} |K| \sum_{k=1}^m \hat{w}_k f(b_k^K)\; \phi_i (b_k^K),
\quad i=1, ..., N_{vi}.
\eeq
We recall that (\ref{fem-linsys}) and (\ref{fem-matrix}) have been obtained under
the Dirichlet boundary condition (\ref{bvp-bc}). It is not difficult to show that a linear system in the same form
can be obtained for mixed boundary conditions provided that
$\Gamma_D \ne \emptyset$, with $\Gamma_D$ being the part of the boundary
where the Dirichlet condition is imposed. Therefore,
the mesh conditions developed below also work for mixed boundary conditions with $\Gamma_D \ne \emptyset$.

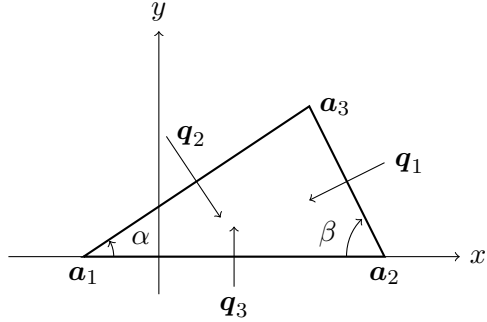
\begin{figure}[thb]
\centering
\begin{tikzpicture}[scale = 1]
\draw [->] (-2,0) -- (4,0);
\draw [right] (4,0) node {$x$};
\draw [->] (0,-0.5) -- (0,3);
\draw [above] (0,3) node {$y$};
\draw [thick] (-1,0) -- (3,0) -- (2, 2) -- cycle;
\draw [below] (-1,0) node {$\V{a}_1$};
\draw [below] (3,0) node {$\V{a}_2$};
\draw [right] (2,2) node {$\V{a}_3$};
\draw [->] (-0.6,0) arc (0:30:0.42);
\draw [right] (-0.5,0.25) node {$\alpha$};
\draw [<-] (2.7,0.5) arc (135:180:0.7);
\draw [left] (2.5,0.3) node {$\beta$};
\draw [->] (1,-0.4) -- (1,0.4);
\draw [below] (1,-0.4) node {$\V{q}_3$};
\draw [->] (3.0,1.25) -- (2,0.75);
\draw [right] (3.0,1.25) node {$\V{q}_1$};
\draw [->] (0.1,1.6) -- (0.82,0.52);
\draw [right] (0.1,1.6) node {$\V{q}_2$};
\end{tikzpicture}
\caption{A sketch of the $\V{q}$ vectors for an arbitrary element. The angles sharing the edge connecting vertices
$\V{a}_1$ and $\V{a}_2$ are $\alpha$ and $\beta$.}
\label{f1}
\end{figure}

We now study under what mesh conditions the scheme (\ref{fem-linsys}) satisfies DMP.
Our basic tool is Lemma \ref{lem1.2}, i.e., we show that $A$ is an $M$-matrix
and has non-negative row sums
when the mesh satisfies the condition (\ref{g-nonobtuse}) below.
To this end, we first introduce some notation.
Denote the vertices of $K$ by
$\V{a}_1^K, \V{a}_2^K, \cdots, \V{a}_{d+1}^K$. The edge matrix of $K$ is defined as
\beq
\nn
E_K = [\V{a}_2^K-\V{a}_1^K, \,\V{a}_3^K-\V{a}_1^K, \, \cdots, \, \V{a}_{d+1}^K-\V{a}_1^K ].
\eeq
From the definition of simplices, $E_K$ is nonsingular \cite{Som29}.
Then, a set of $\V{q}$-vectors (cf. Fig. \ref{f1}) can be defined as
\beq
\label{def-q}
[\V{q}_2^K, \, \V{q}_3^K, \, \cdots, \, \V{q}_{d+1}^K] = E_K^{-T},
\quad \V{q}_1^K = - \sum \limits_{i=2}^{d+1} \V{q}_i^K .
\eeq
This set of vectors has the following properties.
\begin{enumerate}

\item[(i)] By definition, it follows that
\beq
\label{def-q2}
\begin{array}{ll}
\V{q}_i^K \cdot (\V{a}_j^K-\V{a}_1^K) = \delta_{ij}, \\
\V{q}_1^K \cdot (\V{a}_j^K-\V{a}_i^K) = \delta_{1j},
\end{array}
\quad i=2, \cdots,d+1;\; j=1, \cdots, d+1
\eeq
where $\delta_{ij}$ is the Kronecker delta function.

\item[(ii)]
Denote by $\V{S}_i^K$ the face opposite to vertex $\V{a}_i^K$ (i.e., the face not having $\V{a}_i$
as a vertex). Then (\ref{def-q2}) implies that $\V{q}_i^K$ is the inward normal to the face $\V{S}_i^K$;
see Fig. \ref{f1}.

\item[(iii)]
The dihedral angle, $\alpha_{ij}$, between any two faces $\V{S}_i^K$ and $\V{S}_j^K$ ($i \ne j$)
is defined as the supplement of the angle between the inward normals to the faces.
It can be calculated by
\beq
\label{def-dihedral}
\cos(\alpha_{ij}) = - \frac{\V{q}_i^K \cdot \V{q}_j^K }{ \|\V{q}_i^K\|\; \|\V{q}_j^K\| }, \quad i \ne j.
\eeq

\item[(iv)] It is known \cite{BKK07,KL95} that, for any vertex of $K$ with the global and local indices $i$
and $i_K$, respectively, there holds
\beq
\label{q-phi}
\nabla \phi_{i}|_K = \V{q}_{i_K}^K.
\eeq

\end{enumerate}

The main result of this section is stated in the following theorem.

\begin{thm}
\label{thm-g-nonobtuse}
If the mesh satisfies the anisotropic non-obtuse angle condition
\beq
\label{g-nonobtuse}
(\V{q}_i^K)^T \, \mathbb{D}_K \, \V{q}_j^K \le 0,
\quad \forall i \ne j,\; i,j = 1,2, ..., d+1,\; \forall K \in \mathcal{T}_h
\eeq
then the linear finite element scheme (\ref{disc-3}) for solving BVP (\ref{bvp-pde}) and (\ref{bvp-bc})
satisfies DMP.
\end{thm}

{\bf Proof.} We prove this theorem using Lemma \ref{lem1.2}. That is, we show that
the stiffness matrix $A$ has non-negative row sums and is an $M$-matrix when the mesh
satisfies condition (\ref{g-nonobtuse}).

(i) We first show that $A$ has non-negative row sums. From (\ref{fem-matrix}) we only need to show
$\sum_{j=1}^{N_v} a_{ij} \ge 0$ for $i = 1, ..., N_{vi}$. From (\ref{stiffness-matrix}) we have
\bey
\sum_{j=1}^{N_v} a_{ij} & = &
\sum_{j=1}^{N_v} \sum_{K \in \mathcal{T}_h} |K| \;(\nabla \phi_i|_K)^T \; \mathbb{D}_K \; \nabla \phi_j|_K
\nn \\
& = &
\sum_{K \in \mathcal{T}_h} |K| \;(\nabla \phi_i|_K)^T \; \mathbb{D}_K \; \nabla
\left. \left ( \sum_{j=1}^{N_v} \phi_j\right )\right |_K
\nn \\
& = & 0,
\eey
where we have used the fact that $\sum_{j=1}^{N_v} \phi_j(\V{x}) \equiv 1$ for any $\V{x} \in K$.

(ii) Next we show that
\bey
a_{ij} \le 0, && \quad\forall \; i \neq j,\; i, j = 1, ..., N_{v}
\label{thm-g-nonobtuse-1}
\\
a_{ii} \ge 0, && \quad\forall \; i = 1, ..., N_{v} .
\label{thm-g-nonobtuse-2}
\eey
Let $\omega_i$ (or $\omega_j$) be the patch of the elements containing $\V{a}_i$ (or $\V{a}_j$)
as a vertex. Notice that $\nabla \phi_i|_K = 0$ when $K \notin \omega_i$.
Denote the local indices of vertices $\V{a}_i$ and $\V{a}_j$ on $K$ by $i_K$ and $j_K$, respectively.
Then from (\ref{stiffness-matrix}), (\ref{q-phi}), and (\ref{g-nonobtuse}), we have,
for $i\ne j,\; i=1, ..., N_{vi},\; j = 1, ..., N_{v}$,
\bey
a_{ij} & = & \sum_{K \in \omega_i\cap \omega_j}
|K| \; (\nabla \phi_i|_K)^T \; \mathbb{D}_K \; \nabla \phi_j|_K
\nn \\
& = & \sum_{K \in \omega_i\cap \omega_j}
 |K| (\V{q}_{i_K}^K)^T \; \mathbb{D}_K \; \V{q}_{j_K}^K
 \label{thm-g-nonobtuse-3}
 \\
 & \le & 0.
\label{thm-g-nonobtuse-4}
\eey
From (\ref{fem-matrix}) it is obvious that $a_{ij} = 0$ for $i\ne j,\; i = N_{vi}+1,\cdots N_{v},\; j = 1, ..., N_{v}$.
Hence, the off-diagonal entries of $A$ are non-positive.

The inequality (\ref{thm-g-nonobtuse-2}) follows immediately from (\ref{fem-matrix}), (\ref{stiffness-matrix}),
and the positive definiteness of $\mathbb{D}_K$.

(iii) We now show that $A_{11}$ defined in (\ref{fem-matrix}) is an $M$-matrix.
Notice that the non-negativeness of the row sums of $A$ and the properties (\ref{thm-g-nonobtuse-1})
and (\ref{thm-g-nonobtuse-2}) imply that $A_{11}$ is diagonally dominant.
In theory, we can show that $A_{11}$ is an $M$-matrix by proving it is irreducible \cite{Var62}.
However, we will need to assume that any pair of interior
vertices is connected at least by an interior edge path \cite{DDS04}. To avoid
this additional restriction on the mesh, we instead opt to show $A_{11}$ is symmetric and positive definite,
which together with (\ref{thm-g-nonobtuse-1}) and (\ref{thm-g-nonobtuse-2}) implies
that $A_{11}$ is an $M$-matrix \cite{Var62}.

From (\ref{stiffness-matrix}) it is obvious that $A_{11}$ is symmetric.
It suffices to show $A_{11}$ is positive definite.  From the strictly positive definiteness
of the diffusion matrix $\mathbb{D}$,
there exists a positive constant $\beta$ such that
\[
 \mathbb{D}_K \ge \beta \V{I}, \quad \forall K \in \mathcal{T}_h.
\]
For any vector $\V{v} = (v_1, ..., v_{N_{vi}})^T$, we define $v^h = \sum \limits_{i=1}^{N_{vi}} v_i \phi_i \in U_0^h$.
From the definition of $A_{11}$ and the fact that $\nabla v^h|_K$ is constant on $K$, we have
\bey
\nn
\V{v}^T A_{11} \V{v} & = & \sum_{K \in \mathcal{T}_h} |K| \; (\nabla v^h|_K)^T
\; \mathbb{D}_K \; \nabla v^h|_K \\
\nn
& \ge & \beta \sum_{K \in \mathcal{T}_h} |K| \; (\nabla v^h|_K)^T
\; \nabla v^h|_K \\
\nn
& = & \beta \sum_{K \in \mathcal{T}_h} \int_K (\nabla v^h)^T
\nabla v^h d\V{x}
\\
\nn
& = & \beta \int_{\Omega} (\nabla v^h)^T \nabla v^h d\V{x} \\
\nn
& \ge & \beta C_{p} \int_{\Omega} |v^h|^2 d \V{x},
\eey
where in the last step we have used Poincare's inequality and $C_p > 0$ is the associated constant.
For any nonzero vector $\V{v}$, $v^h=\sum \limits_{i=1}^{N_{vi}} v_i \phi_i \not \equiv 0$ and is piecewise linear
and continuous on $\Omega$.
Consequently,
\[
\V{v}^T A_{11} \V{v} \ge \beta C_p \int_{\Omega} |v^h|^2 d \V{x} > 0,
\quad \forall \V{v} \ne 0
\]
which implies that $A_{11}$ is positive definite. Hence, $A_{11}$ is an $M$-matrix.

(iv) From (\ref{fem-matrix}) it is easy to verify that the inverse of $A$ is given by
\[
A^{-1}  = \left [\begin{array}{rr} A_{11}^{-1} & - A_{11}^{-1} A_{12} \\ 0 & I \end{array}\right ] .
\]
Then (\ref{thm-g-nonobtuse-1}) and the fact $A_{11}^{-1} \ge 0$ imply that $A^{-1} \ge 0$
and therefore $A$ is an $M$-matrix.

We have shown above that $A$ is an $M$-matrix and has non-negative row sums.
By Lemma \ref{lem1.2} we conclude that the linear FEM satisfies DMP when the simplicial mesh
satisfies (\ref{g-nonobtuse}).
\proofend


\begin{rem}{\em
\label{rem2.1}
For the isotropic case where $\mathbb{D} = a(\V{x}) \V{I}$ for some scalar function
$a(\V{x})$, condition (\ref{g-nonobtuse}) reduces to the well known non-obtuse angle condition \cite{BKK07,CR73}
\beq
\label{non-obtuse}
 \V{q}_i^K \cdot \V{q}_j^K \le 0,\quad \forall i \ne j, \; \forall K \in \mathcal{T}_h,
\eeq
which requires the dihedral angles $\alpha_{ij}$ (cf. (\ref{def-dihedral}))
of all mesh elements be non-obtuse.
Thus, condition (\ref{g-nonobtuse}) is a generalization of the non-obtuse angle condition.
An alternative interpretation of (\ref{g-nonobtuse}) is that the dihedral angles
of element $K$, measured in the Riemannian metric $\mathbb{D}_K$ (piecewise constant),
are non-obtuse.
}\proofend\end{rem}

\begin{rem}{\em
\label{rem2.2}
It is interesting to point out that an explicit mesh condition similar to (\ref{g-nonobtuse}) is obtained
by Eigestad et al. \cite{EAE02} for a multipoint flux approximation (MPFA) finite volume method
on triangular meshes for anisotropic homogeneous media (i.e., $\mathbb{D}$ is constant).
Moreover, (\ref{g-nonobtuse}) reduces to
a mesh condition obtained by Li et al. \cite{LSS07} for a similar situation with 
constant $\mathbb{D}$ and triangular meshes.
	To see this,
	let the eigen-decomposition of the constant diffusion matrix $\mathbb{D}$ be
	\beq
	\mathbb{D} = \left [\begin{array}{rr} \cos \theta & - \sin \theta \\ \sin \theta & \cos \theta \end{array} \right ]
	\left [\begin{array}{rr} k_1 & 0 \\ 0 & k_2 \end{array} \right ]
	\left [\begin{array}{rr} \cos \theta &  \sin \theta \\ - \sin \theta & \cos \theta \end{array} \right ] .
	\label{D-eigen-d}
	\eeq
	For an arbitrary triangular element $K$, denote the angles sharing the edge connecting vertices
	$\V{a}_1$ and $\V{a}_2$ by $\alpha$ and $\beta$; see Fig. \ref{f1}. Then, a mesh condition of \cite{LSS07}
	is given by
	\beq
	\begin{cases}
	- k_1 \sin \beta \sin \alpha + k_2 \cos \beta \cos \alpha \le 0,& \\
	- k_2 \cos \beta \le 0,& \\
	- k_2 \cos \alpha \le 0,&
	\end{cases}
	\label{LSS}
	\eeq
	provided that the edge connecting
	$\V{a}_1$ and $\V{a}_2$ is parallel to the primary diffusion direction $(\cos \theta, \sin\theta)^T$
	(the eigenvector corresponding to the first eigenvalue of $\mathbb{D}$, $k_1$).
	We now show that (\ref{g-nonobtuse}) reduces to (\ref{LSS}) for the current situation.
	Without loss of generality we assume that the primary diffusion direction
	and the edge connecting $\V{a}_1$ and $\V{a}_2$  are in the direction of the $x$-axis;
	cf. Fig. \ref{f1}. (In this case we have $\theta = 0$.) It is not difficult to obtain
	\[
	\V{q}_1 = c_1 \left [\begin{array}{r} - \sin\beta \\ - \cos\beta \end{array}\right ],\quad
	\V{q}_2 = c_2 \left [\begin{array}{r}  \sin\alpha \\ - \cos\alpha \end{array}\right ],\quad
	\V{q}_3 = c_3 \left [\begin{array}{r} 0 \\ 1 \end{array}\right ],
	\]
	where $c_1$, $c_2$, and $c_3$ are positive constants. From these and (\ref{D-eigen-d}),
	(\ref{g-nonobtuse}) reduces to
	\beq
	\begin{cases}
	\V{q}_1^T \mathbb{D}_K \V{q}_2 =
	\V{q}_1^T \mathbb{D} \V{q}_2 = c_1 c_2 (-k_1 \sin\alpha \sin\beta + k_2 \cos \alpha \cos \beta) \le 0,\\
	\V{q}_1^T \mathbb{D}_K \V{q}_3
	= \V{q}_1^T \mathbb{D} \V{q}_3 = c_1 c_3 (- k_2 \cos \beta ) \le 0,\\
	\V{q}_2^T \mathbb{D}_K \V{q}_3 =
	\V{q}_2^T \mathbb{D} \V{q}_3 = c_2 c_3 (- k_2 \cos \alpha ) \le 0,
	\end{cases}
	\nn
	\eeq
	which gives (\ref{LSS}).
}\proofend\end{rem}

\vspace{15pt}

It is often more convenient to express the anisotropic non-obtuse angle condition (\ref{g-nonobtuse})
in terms of mapping $F_K$ from $\hat K$
to $K$. Denote the Jacobian matrix of $F_K$ by $F_K'$.
We define the vectors $\hat{\V{q}}_k$, $k=1, ..., d+1$ for the reference element $\hat{K}$
as in (\ref{def-q}). The chain rule of differentiation implies
\[
\nabla \phi_i = (F_K')^{-T} \nabla_{\xi} \hat{\phi}_i,
\]
where $\hat{\phi}_i(\V{\xi}) = \phi_i(F_K(\V{\xi}))$. From (\ref{q-phi}), we have
\[
\V{q}_i = (F_K')^{-T} \hat{\V{q}}_i .
\]
Inserting this into (\ref{g-nonobtuse}) we obtain the following theorem.

\begin{thm}
\label{thm-mesh-map}
If the mesh satisfies
\beq
\label{g-nonobtuse-2}
\hat{\V{q}}_i^T (F_K')^{-1} \mathbb{D}_K (F_K')^{-T} \hat{\V{q}}_j \le 0,
\quad \forall i \ne j,\; i,j=1, ..., d+1,\; \forall K \in \mathcal{T}_h
\eeq
then the linear finite element scheme (\ref{disc-3}) for solving BVP (\ref{bvp-pde}) and (\ref{bvp-bc})
satisfies DMP.
\end{thm}

\begin{co}
\label{co2.1}
Suppose that the reference element $\hat{K}$ is taken as a simplex with non-obtuse
dihedral angles. If the mesh satisfies
\beq
\label{g-nonobtuse-3}
(F_K')^{-1} \mathbb{D}_K (F_K')^{-T} = C_K I,\quad \forall K \in \mathcal{T}_h
\eeq
where $C_K$ is a positive constant on $K$ and $I$ is the $d\times d$ identity matrix, then
the linear finite element scheme (\ref{disc-3}) for solving BVP (\ref{bvp-pde}) and (\ref{bvp-bc})
satisfies DMP.
\end{co}

{\bf Proof.} Since $\hat{K}$ is a simplex with non-obtuse
dihedral angles, we have
\[
\hat{\V{q}}_i^T \hat{\V{q}}_j \le 0,\; i \ne j,\quad i,j=1,..., d+1.
\]
From this it is easy to see that (\ref{g-nonobtuse-3}) is sufficient for (\ref{g-nonobtuse-2}) to hold.
\proofend

\vspace{15pt}

In the next two sections mesh condition (\ref{g-nonobtuse-3}) will be used to develop
metric tensors accounting for DMP satisfaction and mesh adaptivity. These metric tensors
are needed in anisotropic mesh generation.
It is emphasized that (\ref{g-nonobtuse-3}), as well as mesh conditions (\ref{g-nonobtuse}) and
(\ref{g-nonobtuse-2}), can also be used more directly via direct minimization \cite{LSS07,MD06}
or variational formulation \cite{Hua01b} for optimizing the current mesh
to improve DMP satisfaction.

\section{Anisotropic mesh generation: metric tensor based on DMP satisfaction}
\label{Sec-metric}

In this section we develop a metric tensor for use in anisotropic mesh generation based on
mesh condition (\ref{g-nonobtuse-3}). To this end, we adopt the so-called $M$-uniform
mesh approach \cite{Hua05b, Hua06}
where an anisotropic mesh is viewed as an $M$-uniform mesh or a uniform
one in the metric specified by a tensor $M=M(\V{x})$.
The tensor, chosen to be symmetric and positive definite,
provides the information on the size, shape, and orientation of mesh elements over $\Omega$
necessary for the actual implementation of mesh generation.
Various formulations of the metric tensor
have been developed in the past for anisotropic mesh adaptation; e.g., see
\cite{BGHLS97, CHMP97, FG00, Hua05b,HL08}.
Once a metric tensor has been determined, the corresponding anisotropic meshes can be generated
using a variety of techniques including
blue refinement \cite{KR90, Lan91},
directional refinement \cite{Rac93,Rac97},
Delaunay-type triangulation \cite{BGHLS97,BGM97,CHMP97,PVMZ87},
front advancing \cite{GS98},
bubble packing \cite{YS00},
local refinement and modification \cite{ABHFDV02,BH96,DVBFH02,HDBAFV00,RLSF03},
and variational mesh generation \cite{BS82,Dvi91,Hua01b,Jac88,Knu96,KMS02,Win81}.
In this paper we restrict our attention to the determination of a metric tensor for DMP satisfaction,
and refer the interested reader to the above mentioned references for meshing strategies.

It is shown in \cite{Hua06} that when the reference element $\hat{K}$ is taken to be equilateral and
unitary in volume,
a simplicial $M$-uniform mesh $\mathcal{T}_h$ for a given
$M=M(\V{x})$ satisfies
\bey
\label{cond-equi}
\rho_K |K| & = & \frac{\sigma_h}{N}, \quad \forall K \in \mathcal{T}_h \\
\label{cond-align}
\frac{1}{d} \mbox{tr} \left ( (F_K')^T M_K F_K' \right ) & = &
\mbox{det} \left ( (F_K')^T M_K F_K' \right )^{\frac{1}{d}}, \quad \forall K \in \mathcal{T}_h
\eey
where $N$ is the number of mesh elements, $F_K$ is the affine mapping from $\hat{K}$ to $K$,
$F_K'$ is the Jacobian matrix of $F_K$,  and
\beq
M_K = \frac{1}{|K|} \int_K M(\V{x}) d \V{x},\quad
\rho_K=\sqrt{\mbox{det}(M_K)}, \quad
\sigma_h = \sum_{K \in \mathcal{T}_h} \rho_K |K|.
\eeq
Condition (\ref{cond-equi}), referred to as {\em the equidistribution condition}, determines the size of $K$
from $\rho_K$. The larger $\rho_K$ is,  the smaller $|K|$ is.
On the other hand, (\ref{cond-align}), called {\em the alignment condition},
characterizes the shape and orientation of
$K$ in the sense that the principal axes of the circumscribed ellipsoid of $K$ are parallel to
the eigenvectors of $M_K$ while their lengths are reciprocally proportional to
the square roots of the respective eigenvalues \cite{Hua06}.

To determine $M$ from mesh condition
(\ref{g-nonobtuse-3}), we first notice that the left and right sides of (\ref{cond-align}) represents
the arithmetic and geometric means of the eigenvalues
of matrix $(F_K')^T M_K F_K'$, respectively.
From the arithmetic-mean geometric-mean inequality, (\ref{cond-align})
implies that all of the eigenvalues are equal to each other. In other words, $(F_K')^T M_K F_K'$
is a scalar matrix, i.e.,
\beq
\label{M-unif-cond}
(F_K')^T M_K F_K' = \tilde{C}_K I
\quad \text{ or } \quad (F_K')^{-1} M_K^{-1} (F_K')^{-T} = \tilde{C}_K^{-1} I
\eeq
for some constant $\tilde{C}_K$. A direct comparison of
(\ref{M-unif-cond}) with (\ref{g-nonobtuse-3}) suggests that the metric tensor $M$ be chosen in the form
\beq
\label{M-DMP}
M_{DMP,K} = \theta_K \mathbb{D}_K^{-1}, \quad \forall K \in \mathcal{T}_h
\eeq
where $\theta = \theta_K > 0$ is an arbitrary piecewise constant function.
Thus, any $M$-uniform mesh associated with a metric tensor in the form (\ref{M-DMP})
satisfies condition (\ref{g-nonobtuse-3}). The following theorem follows from Corollary \ref{co2.1}.

\begin{thm}
\label{stiff-M-thm}
Suppose that the reference element $\hat{K}$ is taken to be equilateral and unitary in volume.
For an $M$-uniform mesh associated with any metric tensor in the form (\ref{M-DMP}),
the linear finite element scheme (\ref{disc-3}) for solving BVP (\ref{bvp-pde}) and (\ref{bvp-bc})
satisfies DMP.
\end{thm}

\begin{rem}
\label{rem3.1}
{\rm
Since an $M$-uniform mesh is aligned with the metric tensor $M$
as characterized by the alignment condition (\ref{cond-align}),
we can conclude that when $M$ is chosen in the form (\ref{M-DMP}),
a corresponding $M$-uniform mesh is aligned with
the diffusion matrix $\mathbb{D}$ in the sense that the principal axes of the circumscribed ellipsoid of
element $K$ are parallel to the eigenvectors of $\mathbb{D}_K$ while their lengths are
proportional to the square roots of the respective eigenvalues.
As a consequence, the length of $K$ is greater in a faster diffusion direction and
smaller in a slower diffusion direction. A small length scale of mesh elements
in slow diffusion directions helps reduce numerical dissipation in those directions.
}\proofend\end{rem}

\begin{rem}
\label{rem3.2}
{\rm
Note that $\theta = \theta_K$ in (\ref{M-DMP}) is arbitrary. Thus,
in addition to satisfying DMP, there is a degree of freedom for the mesh
to account for other considerations. In the next section we shall consider
mesh adaptation and choose $\theta_K$ to minimize a certain error bound.
}\proofend\end{rem}

\section{Metric tensors based on DMP satisfaction and mesh adaptivity}
\label{Sec-adapt}

In this section we develop a metric tensor taking both the satisfaction of DMP
and mesh adaptivity into consideration. The metric tensor takes the form (\ref{M-DMP}),
with the scalar function $\theta = \theta_K$ being determined to minimize an interpolation error bound.
For simplicity, we consider here an error bound for linear Lagrange interpolation. Other interpolation
error bounds (e.g., see \cite{Hua06}) can be considered without major modification.

\begin{lem} (\cite{Hua06})
Let $K \subset \mathbb{R}^d $ be a simplicial element and
$\Pi_h$ be the linear Lagrange interpolation operator. Then,
\beq
\label{bound-iso-1}
|v - \Pi_h v|_{H^1(K)} \le C \|(F_K')^{-1}\| 
\left [ \int_K \left [\mbox{tr}\left ( (F_K')^T |H(v)| F_K'\right )\right ]^2  d \V{x}\right ]^{\frac 1 2},
\quad \forall v \in H^2(K)
\eeq
where $\| \cdot \|$ denotes the $l_2$ matrix norm, $H(v)$ is the Hessian of $v$, and $|H(v)| = \sqrt{H(v)^2}$.
\end{lem}

\begin{lem}
\label{trace-lem}
For any given $d\times d$ symmetric matrix $S$, there holds that
\beq
| \mbox{tr}(A^T S A) | \le \mbox{tr}(A^T  A) \; \| S \|, \quad \forall A \in \mathbb{R}^{d\times d} .
\label{lem4.2-2-2}
\eeq
If $S$ is further positive definite, then
\beq
\| S \|^{-1}\; \mbox{tr}(A^T S A) \le \mbox{tr}(A^T  A) \le \mbox{tr}(A^T S A)\; \| S^{-1} \| .
\label{lem4.2-2-1}
\eeq
\end{lem}

{\bf Proof.} Denote the eigen-decomposition of $S$ by
\[
S = Q \Sigma Q^T,
\]
where $Q$ is an orthogonal matrix, $\Sigma = \mbox{ diag}(\lambda_1, ..., \lambda_d)$,
and $\lambda_i,\; i = 1, ..., d$ are the eigenvalues of $S$. Write
\[
A^T Q = [\V{v}_1, ..., \V{v}_d] .
\]
Then
\[
A^T S A =  (A^T Q) \Sigma (Q^T A) = [\V{v}_1, ..., \V{v}_d] \Sigma [\V{v}_1, ..., \V{v}_d]^T
= \sum_i \lambda_i \V{v}_i \V{v}_i^T .
\]
It follows that
\bey
| \mbox{tr}(A^T S A)| & = & | \sum_i \lambda_i \mbox{tr}(\V{v}_i \V{v}_i^T) |
\nn \\
& = & | \sum_i \lambda_i \| \V{v}_i\|^2 |
\nn \\
& \le & \sum_i  \| \V{v}_i\|^2 \cdot | \lambda |_{max}
\nn \\
& = & \mbox{tr}(A^T A)  \, \| S \| ,
\nn
\eey
which gives (\ref{lem4.2-2-2}). Inequality (\ref{lem4.2-2-1}) follows from (\ref{lem4.2-2-2}) and that
\beq
\mbox{tr}(A^T  A) = \mbox{tr}(A^T S^{\frac{1}{2}} S^{-1} S^{\frac{1}{2}} A) \le \mbox{tr}(A^T S A)\; \| S^{-1} \| .
\label{lem4.2-2-3}
\eeq
\proofend

The scalar function $\theta = \theta_K$ in (\ref{M-DMP})
is determined based on interpolation error bound (\ref{bound-iso-1}).
From the definition of the Frobenius matrix norm, we have
\[
\| A \| \le \|A\|_{F} = \sqrt{\mbox{tr}(A^T A) } = \sqrt{\mbox{tr}(A A^T) },\quad \forall A \in \mathbb{R}^{d\times d}.
\]
Using this, taking squares of both sides of (\ref{bound-iso-1}), and summing the result over all elements
of $\mathcal{T}_h$, we have
\bey
\nn
& & |u - \Pi_h u|^2_{H^1 (\Omega)}  = \sum_{K \in \mathcal{T}_h} |u - \Pi_h u|^2_{H^1 (K)} \\
\nn
& & \le C \sum_{K \in \mathcal{T}_h} \|(F_K')^{-1}\|^{2} 
\int_K \left [\mbox{tr}\left ( (F_K')^T |H(u)| F_K'\right )\right ]^2  d \V{x}
\\
\nn
& & \le C \sum_{K \in \mathcal{T}_h}  \|(F_K')^{-1}\|_F^{2}
\int_K \left [\mbox{tr}\left ( (F_K')^T |H(u)| F_K'\right )\right ]^2  d \V{x}
\\
\nn
& & = C \sum_{K \in \mathcal{T}_h} \left [ \mbox{tr} ( (F_K')^{-1} (F_K')^{-T} ) \right ]
\int_K \left [\mbox{tr}\left ( (F_K')^T |H(u)| F_K'\right )\right ]^2  d \V{x}.
\eey
From Lemma~\ref{trace-lem} it follows that
\bey
&& |u - \Pi_h u|^2_{H^1 (\Omega)}
\nn \\
& & \le C \sum_{K \in \mathcal{T}_h}  \left [ \mbox{tr}( (F_K')^{-1} \mathbb{D}_K (F_K')^{-T} ) \right ]
\cdot \|\mathbb{D}_K^{-1} \| \cdot 
\int_K \left [ \mbox{tr}( (F_K')^T \mathbb{D}_K^{-1} (F_K') ) \right ]^2   \|\mathbb{D}_K |H(u)| \|^2 d \V{x}
\nn \\
&  & = C \sum_{K \in \mathcal{T}_h}  |K| \cdot \left [ \mbox{tr}( (F_K')^{-1} \mathbb{D}_K (F_K')^{-T} ) \right ]
\cdot \left [ \mbox{tr}( (F_K')^T \mathbb{D}_K^{-1} (F_K') ) \right ]^2 
\nn \\
&& \qquad \mbox{ } \qquad \times  \|\mathbb{D}_K^{-1} \| \cdot 
\frac{1}{|K|} \int_K  \|\mathbb{D}_K |H(u)| \|^2 d \V{x} .
\label{bound-iso-2}
\eey

Consider an $M$-uniform mesh $\mathcal{T}_h$ corresponding to a metric tensor $M_K$ in the form
(\ref{M-DMP}). Then, alignment condition (\ref{cond-align}) reduces to
\beq
\frac{1}{d} \mbox{tr} \left ( (F_K')^T \mathbb{D}_K^{-1} F_K' \right ) =
\mbox{det} \left ( (F_K')^T \mathbb{D}_K^{-1} F_K' \right )^{\frac{1}{d}} .
\label{e4-1}
\eeq
From the arithmetic-mean geometric-mean inequality, (\ref{e4-1}) implies that
all of the eigenvalues of matrix $(F_K')^T \mathbb{D}_K^{-1} F_K'$ are equal to each other.
As a consequence, all of the eigenvalues of the inverse of $(F_K')^T \mathbb{D}_K^{-1} F_K'$
are equal to each other, which in turn implies
\beq
 \frac{1}{d} \mbox{tr} \left ( (F_K')^{-1} \mathbb{D}_K (F_K')^{-T} \right ) =
\mbox{det} \left ( (F_K')^{-1} \mathbb{D}_K (F_K')^{-T} \right )^{\frac{1}{d}} .
\label{e4-2}
\eeq
Inserting (\ref{e4-1}) and (\ref{e4-2}) into (\ref{bound-iso-2}) and noticing
\[
\mbox{det} \left ( (F_K')^T \mathbb{D}_K^{-1} F_K' \right )
= |K|^2 \mbox{det} \left ( \mathbb{D}_K \right )^{-1} ,
\quad
\mbox{det} \left ( (F_K')^{-1} \mathbb{D}_K (F_K')^{-T} \right )
= |K|^{-2} \mbox{det} \left ( \mathbb{D}_K \right ),
\]
we have
\bey
|u - \Pi_h u|^2_{H^1 (\Omega)} & \le & C \sum_{K \in \mathcal{T}_h} |K|^{\frac{d+2}{d}}
\mbox{det} \left ( \mathbb{D}_K \right )^{-\frac{1}{d}}
 \|\mathbb{D}_K^{-1} \|\cdot \frac{1}{|K|} \int_K  \|\mathbb{D}_K |H(u)| \|^2 d \V{x} .
\eey
Rewrite this bound as
\beq
|u - \Pi_h u|^2_{H^1 (\Omega)} \le
C \sum_{K \in \mathcal{T}_h} |K|^{\frac{d+2}{d}} B_K ,
\label{bound-iso-3}
\eeq
where
\beq
B_K = \mbox{det} \left ( \mathbb{D}_K \right )^{-\frac{1}{d}}
 \|\mathbb{D}_K^{-1} \|\cdot  \frac{1}{|K|} \int_K  \|\mathbb{D}_K |H(u)| \|^2 d \V{x} .
\label{B_K}
\eeq
Notice that $\int_K  \|\mathbb{D}_K |H(u)| \|^2 d \V{x}$ and therefore $B_K$ can vanish locally.
To ensure the positive definiteness of the metric tensor to be defined, we regularize the above bound
with a parameter $\alpha_h>0$ as
\beq
|u - \Pi_h u|^2_{H^1 (\Omega)}  \le  C \sum_{K \in \mathcal{T}_h} |K|^{\frac{d+2}{d}} \left [\alpha_h + B_K\right ]
 = C \alpha_h \sum_{K \in \mathcal{T}_h} |K|^{\frac{d+2}{d}} \left [1+\frac{1}{\alpha_h} B_K\right ].
 \label{bound-iso-4}
\eeq
From H\"older's inequality, we have
\bey
\sum_{K \in \mathcal{T}_h} |K|^{\frac{d+2}{d}}
\left [1+\frac{1}{\alpha_h} B_K\right ]
& = &
\sum_{K \in \mathcal{T}_h} \left ( |K|
\left [1+\frac{1}{\alpha_h} B_K\right ]^{\frac{d}{d+2}}
 \right )^{\frac{d+2}{d}}
\nn \\
& \ge &
N^{-\frac{2}{d}} \left ( \sum_{K \in \mathcal{T}_h} |K|
\left [1+\frac{1}{\alpha_h} B_K\right ]^{\frac{d}{d+2}}
\right )^{\frac{d+2}{d}} ,
\label{e4-3}
\eey
with equality in the last step if and only if
\beq
|K| \left [1+\frac{1}{\alpha_h} B_K\right ]^{\frac{d}{d+2}}
= \mbox{ constant} ,\quad \forall K \in \mathcal{T}_h .
\label{e4-4}
\eeq
A direct comparison of this with equidistribution condition (\ref{cond-equi}) suggests that
the optimal $\rho_K$ be defined as
\beq
\label{rho-1}
\rho_K = \left [1+\frac{1}{\alpha_h} B_K\right ]^{\frac{d}{d+2}}.
\eeq
From the relation $\rho_K =\sqrt{\mbox{det}(M_K)}$,  we find the optimal $\theta_K$ and $M_K$ as
\beq
\label{theta-1}
\theta_K = \rho_K^{\frac{2}{d}} \mbox{det} \left ( \mathbb{D}_K \right )^{\frac{1}{d}}
= \left [ 1+ \frac{1}{\alpha_h}  B_K \right ]^{\frac{2}{d+2}}
\mbox{det} \left ( \mathbb{D}_K \right )^{\frac{1}{d}},
\eeq
\beq
\label{M-DMP+adap}
M_{DMP+adap, K} = \left [ 1+ \frac{1}{\alpha_h}  B_K \right ]^{\frac{2}{d+2}}
\mbox{det} \left ( \mathbb{D}_K \right )^{\frac{1}{d}} \mathbb{D}_K^{-1},
\eeq
where $B_K$ is defined in (\ref{B_K}).
With the so-defined metric tensor, the error bound can be obtained from (\ref{bound-iso-4}) and (\ref{e4-3})
for a corresponding $M$-uniform mesh as
\beq
\label{bound-iso-5}
| u - \Pi_h u |_{H^1(\Omega)} \le C N^{-\frac{1}{d}} \sqrt{\alpha_h}
\sigma_h^{\frac{d+2}{2 d}} .
\eeq

To complete the definition, we need to determine the regularization parameter $\alpha_h$.
We follow \cite{Hua05b} to define $\alpha_h$ such that
\beq
\sigma_h \equiv \sum_{K \in \mathcal{T}_h} \rho_K |K| \le 2|\Omega|,
\label{sigma-1}
\eeq
with which roughly 50\% of the mesh points are concentrated in regions of large $\rho_K$.
From (\ref{rho-1}) and Jensen's inequality, we have
\bey
\sigma_h & = & \sum_{K \in \mathcal{T}_h}  |K| \left [1+\frac{1}{\alpha_h} B_K\right ]^{\frac{d}{d+2}}
\nn \\
& \le & \sum_{K \in \mathcal{T}_h}  |K| \left [1+\alpha_h^{-\frac{d}{d+2}} B_K^{\frac{d}{d+2}}\right ]
\nn \\
& = & |\Omega| + \alpha_h^{-\frac{d}{d+2}} \sum_{K \in \mathcal{T}_h}  |K| B_K^{\frac{d}{d+2}} .
\eey
By requiring the above bound to be less than or equal to $2 |\Omega |$, we obtain
\beq
\alpha_h = \left ( \frac{1}{|\Omega|} \sum_{K \in \mathcal{T}_h} |K| B_K^{\frac{d}{d+2}} \right )^{\frac{d+2}{d}}.
\label{alpha-1}
\eeq
Combining (\ref{bound-iso-5}) with (\ref{sigma-1}) and (\ref{alpha-1}) and summarizing
the above derivation, we have the following theorem.

\begin{thm}
\label{thm4.1}
Suppose that the reference element $\hat{K}$ is chosen to be equilateral and unitary in volume.
For any $M$-uniform simplicial mesh corresponding to the metric tensor (\ref{M-DMP+adap}),
the linear finite element scheme (\ref{disc-3}) for solving BVP (\ref{bvp-pde}) and (\ref{bvp-bc})
satisfies DMP and the interpolation error for the exact solution $u$ is bounded by
\beq
| u - \Pi_h u |_{H^1(\Omega)} \le C N^{-\frac{1}{d}}
\left ( \sum_{K \in \mathcal{T}_h} |K| B_K^{\frac{d}{d+2}} \right )^{\frac{d+2}{2 d}},
\label{thm4.2-1}
\eeq
where $B_K$ is defined in (\ref{B_K}).
\end{thm}

\vspace{10pt}

It is remarked that the metric tensor (\ref{M-DMP+adap}) (cf. (\ref{B_K}))
depends on the second derivatives of the exact solution $u$ which is what
we are seeking/approximating. In actual computation, the second derivatives are replaced
with approximations obtained with a Hessian recovery technique such as
the one of using piecewise quadratic polynomials fitting in least-squares sense to
nodal values of the computed solution (e.g., see \cite{Hua05b}).
A hierarchical basis error estimator can also be used to approximate the Hessian of the exact solution.
It is shown in \cite{HKL10} that the least-squares fitting and the hierarchical basis methods
work comparably for all considered cases except for one where the diffusion coefficient is discontinuous and
the interfaces are predefined in the mesh. In this case, the latter works better than the former
since hierarchical basis estimation does not over-concentrate mesh elements near the interfaces.
Since our main goal is to study DMP satisfaction instead of the discontinuity of the diffusion coefficient,
we choose to use the least squares fitting method for Hessian recovery in our computation
due to its simplicity and problem independent feature.

It is interesting to note that the term in the bracket in (\ref{thm4.2-1}) can be viewed as a Riemann sum
of an integral, i.e.,
\[
\sum_{K \in \mathcal{T}_h} |K| B_K^{\frac{d}{d+2}}
 \sim \int_\Omega
\mbox{det} \left ( \mathbb{D} \right )^{-\frac{1}{d+2}}
 \|\mathbb{D}^{-1} \|^{\frac{d}{d+2}} \cdot  \|\mathbb{D}|H(u)| \|^{\frac{2d}{d+2}}d \V{x} .
\]
Thus, the interpolation error has an asymptotic bound as
\bey
| u - \Pi_h u |_{H^1(\Omega)} & \le & C N^{-\frac{1}{d}}
\left ( \sum_{K \in \mathcal{T}_h} |K| B_K^{\frac{d}{d+2}} \right )^{\frac{d+2}{2 d}}
\nn \\
& \sim & C   N^{-\frac{1}{d}} \left (  \int_\Omega
\mbox{det} \left ( \mathbb{D} \right )^{-\frac{1}{d+2}}
 \|\mathbb{D}^{-1} \|^{\frac{d}{d+2}} \cdot  \|\mathbb{D} |H(u)| \|^{\frac{2d}{d+2}}
 d \V{x}  \right )^{\frac{d+2}{2 d}}.
\label{bound-iso-6}
\eey

We emphasize that both the satisfaction of DMP and mesh adaptation (through minimization of an error bound)
are taken into account in the definition of metric tensor (\ref{M-DMP+adap}).
An interesting question is what the interpolation error bound looks like
if mesh adaptation is not taken into consideration. For example, we consider a case $\theta_K = 1$ in
(\ref{M-DMP}). This gives the metric tensor
\beq
M_K = \mathbb{D}_K^{-1} .
\label{M-DMP-2}
\eeq
Recall that the interpolation error is bounded in (\ref{bound-iso-3}), i.e.,
\beq
|u - \Pi_h u|_{H^1 (\Omega)} \le
C \left (\sum_{K \in \mathcal{T}_h} |K|^{\frac{d+2}{d}} B_K\right )^{\frac 1 2}  ,
\label{bound-iso-8}
\eeq
where $B_K$ is defined in (\ref{B_K}). Moreover,
for an $M$-uniform mesh corresponding to this metric tensor the equidistribution condition (\ref{cond-equi})
reduces to
\beq
\mbox{det}(\mathbb{D}_K)^{-\frac{1}{2}} |K| = \frac{\sigma_h}{N},
\label{cond-equi-2}
\eeq
where $\sigma_h = \sum\limits_{K \in \mathcal{T}_h} \mbox{det}(\mathbb{D}_K)^{-\frac{1}{2}} |K| $.
Inserting (\ref{cond-equi-2}) into (\ref{bound-iso-8}), we have
\bey
|u - \Pi_h u|_{H^1 (\Omega)} & \le &
C \left (\sum_{K \in \mathcal{T}_h} |K|
\left (\mbox{det}(\mathbb{D}_K)^{\frac{1}{2}} \frac{\sigma_h}{N} \right )^{\frac{2}{d}} B_K\right )^{\frac 1 2}
\nn \\
& = & C N^{-\frac{1}{d}}\sigma_h^{\frac{1}{d}}
\left (\sum_{K \in \mathcal{T}_h} |K| \mbox{det}(\mathbb{D}_K)^{\frac{1}{d}}  B_K\right )^{\frac 1 2}
\nn \\
& = & C N^{-\frac{1}{d}}\left ( \sum_{K \in \mathcal{T}_h} \mbox{det}(\mathbb{D}_K)^{-\frac{1}{2}} |K|
\right )^{\frac{1}{d}}
\left (\sum_{K \in \mathcal{T}_h} |K| \mbox{det}(\mathbb{D}_K)^{\frac{1}{d}}  B_K\right )^{\frac 1 2} .
\nn
\eey
Thus,
\bey
|u - \Pi_h u|_{H^1 (\Omega)} & \le &
C N^{-\frac{1}{d}}\left ( \sum_{K \in \mathcal{T}_h} \mbox{det}(\mathbb{D}_K)^{-\frac{1}{2}} |K|
\right )^{\frac{1}{d}}
\left (\sum_{K \in \mathcal{T}_h} |K| \mbox{det}(\mathbb{D}_K)^{\frac{1}{d}}  B_K\right )^{\frac 1 2} 
\label{bound-iso-9}
\\
& \sim & C N^{-\frac{1}{d}} \left ( \int_\Omega \mbox{det}(\mathbb{D})^{-\frac{1}{2}} d \V{x} \right )^{\frac{1}{d}}
\left ( \int_\Omega \|\mathbb{D}^{-1} \|\cdot  \|\mathbb{D} |H(u)| \|^2 d \V{x} \right )^{\frac 1 2} .
\label{bound-iso-10}
\eey
This is the interpolation error bound for an $M$-uniform mesh corresponding to metric tensor (\ref{M-DMP-2}).

From H\"older's inequality, it follows that
\[
\left ( \sum_{K \in \mathcal{T}_h} |K| B_K^{\frac{d}{d+2}} \right )^{\frac{d+2}{2d}}
\le \left ( \sum_{K \in \mathcal{T}_h} \mbox{det}(\mathbb{D}_K)^{-\frac{1}{2}} |K|
\right )^{\frac{1}{d}}
\left (\sum_{K \in \mathcal{T}_h} |K| \mbox{det}(\mathbb{D}_K)^{\frac{1}{d}}  B_K\right )^{\frac 1 2} .
\]
Thus, the solution-dependent factor of bound (\ref{thm4.2-1})
is small than or equal to that of bound (\ref{bound-iso-9}).
In this sense, $M_{DMP+adap}$ defined in (\ref{M-DMP+adap}) leads to a more accurate interpolant
than $M_{DMP}$ defined in (\ref{M-DMP-2}) (or  (\ref{M-DMP}) with $\theta_K = 1$).

Moreover, from the standard interpolation theory we recall that the interpolation error for a uniform mesh
is bounded by
\beq
|u - \Pi_h u|_{H^1 (\Omega)} \le C N^{-\frac{1}{d}}
\left ( \int_\Omega \| \nabla^2 u \|^2 d \V{x} \right )^{\frac 1 2} .
\label{bound-iso-11}
\eeq
It is easy to see that the solution dependent factor in error bound (\ref{bound-iso-6}) for $M_{DMP+adap}$
is in the order of $|\nabla^2 u|_{L^{\frac{2d}{d+2}}(\Omega)}$ and those in (\ref{bound-iso-10}) for $M_{DMP}$ and (\ref{bound-iso-11}) for a uniform mesh are in the order of $|\nabla^2 u|_{L^{2}(\Omega)}$. Thus,
(\ref{bound-iso-6}) has the smallest solution dependent factor, an indication of the advantage of using adaptive
meshes. On the other hand, the error bounds (\ref{bound-iso-6}) and (\ref{bound-iso-10}) depend on
the determinant and norm of the diffusion matrix $\mathbb{D}$ and its inverse. This indicates that
DMP satisfaction may sacrifice accuracy. Indeed, as we shall see in the next section,  the solution error
for DMP-bound meshes can sometimes be larger than that for a uniform mesh.

\section{Numerical results}
\label{Sec-results}

In this section we present three two-dimensional examples to demonstrate the performance of metric tensors
$M_{DMP}$ in (\ref{M-DMP}) with $\theta_K = 1$ based on DMP satisfaction and $M_{DMP+adap}$ in
(\ref{M-DMP+adap}) combining DMP satisfaction and mesh adaptivity.
For comparison purpose, we also include numerical results obtained with almost uniform meshes (labelled
with $M_{unif}$) and with a metric tensor $M_{adap}$ based on minimization of a bound
on the $H^1$ semi-norm of linear interpolation error \cite{Hua05b}:
\beq
\label{M-adap}
M_{adap,K}=\rho_K^{\frac{2}{d}}\det
\left( I+\frac{1}{\alpha _{h}}|H_K(u)|\right) ^{-\frac{1}{d}}
\left[ I+\frac{1}{\alpha _{h}}|H_K(u)|\right] ,
\eeq
where
\beq
\nn
\rho_K = \Big \| I + \frac{1}{\alpha_h} | H_K(u) | \Big \|_F ^{\frac{d}{d+2}
}\,\det \left( I+\frac{1}{\alpha _{h}}|H_K(u)|\right)^{\frac{1}{d+2}},
\eeq
and $\alpha_h$ is defined implicitly through
\[
\sum_{K\in \mathcal{T}_h} |K| \rho_K = 2 | \Omega | .
\]
Once again, the second derivatives of the exact solution are replaced 
in actual computation with approximations obtained with a Hessian recovery technique
(using piecewise quadratic polynomials fitting in least-squares sense to
nodal values of the computed solution \cite{Hua05b}).

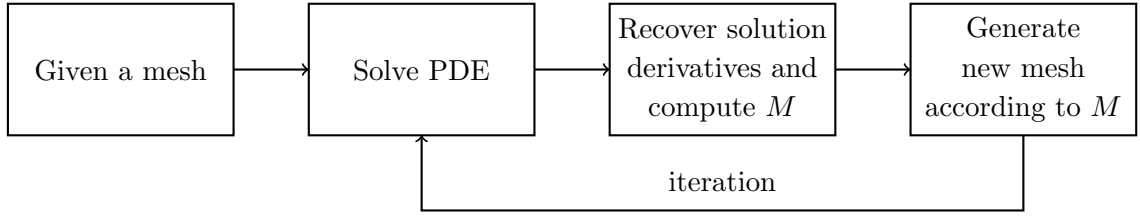
\begin{figure}
\centering
\tikzset{my node/.code=\ifpgfmatrix\else\tikzset{matrix of nodes}\fi}
\begin{tikzpicture}[every node/.style={my node},scale=0.5]
\draw[thick] (0,0) rectangle (6,3.5);
\draw[thick] (8,0) rectangle (14,3.5);
\draw[thick] (16,0) rectangle (22,3.5);
\draw[thick] (24,0) rectangle (30,3.5);
\draw[->,thick] (6,1.75)--(8,1.75);
\draw[->,thick] (14,1.75)--(16,1.75);
\draw[->,thick] (22,1.75)--(24,1.75);
\draw[->,thick] (27,0)--(27,-2)--(11,-2)--(11,0);
\node[above] at (19,-2) {iteration\\};
\node (node1) at (3,1.75) {Given a mesh\\};
\node (node2) at (11,1.75) {Solve PDE\\};
\node (node3) at (19,1.75) {Recover solution\\ derivatives and\\ compute $M$\\};
\node (node4) at (27,1.75) {Generate\\ new mesh\\ according to $M$\\};
\end{tikzpicture}
\caption{An iterative procedure for adaptive mesh solution of PDEs}
\label{f.1}
\end{figure}

An iterative procedure for solving PDEs is shown in Fig. \ref{f.1}. In the current computation, each run
is stopped after ten iterations.  We have found that there is very little improvement in the computed solution
after ten iterations for all the examples considered.
A new mesh is generated using the computer code BAMG (bidimensional
anisotropic mesh generator) developed by Hecht \cite{Hec97} based on
a Delaunay-type triangulation method \cite{CHMP97}.
The code allows the user to supply his/her own metric tensor defined on a background
mesh. In our computation, the background mesh has been taken as the most recent mesh available.

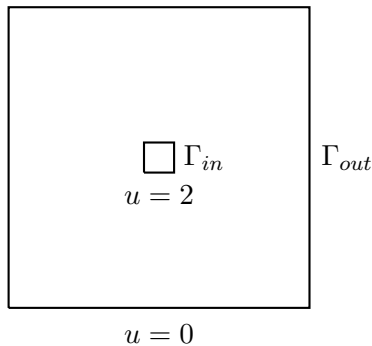
\begin{figure}[thb]
\begin{center}
\begin{tikzpicture}[scale = 1.0]
\draw [thick] (0,0) -- (0,4) -- (4, 4) -- (4, 0) -- (0, 0);
\draw [thick] (1.8,1.8) -- (1.8,2.2) -- (2.2, 2.2) -- (2.2, 1.8) -- (1.8, 1.8);
\draw (2, -0.35) node {$ u = 0 $};
\draw (4.5, 2) node {$ \Gamma_{out} $};
\draw (2, 1.5) node {$ u = 2 $};
\draw (2.6, 2) node {$ \Gamma_{in} $};
\end{tikzpicture}
\end{center}
\caption{The physical domain and boundary conditions for Example~\ref{ex1}.}
\label{ex1-domain}
\end{figure}

\begin{figure}[thb]
\centering
\hbox{
\hspace{10mm}
\begin{minipage}[b]{2.5in}
\includegraphics[width=2.5in]{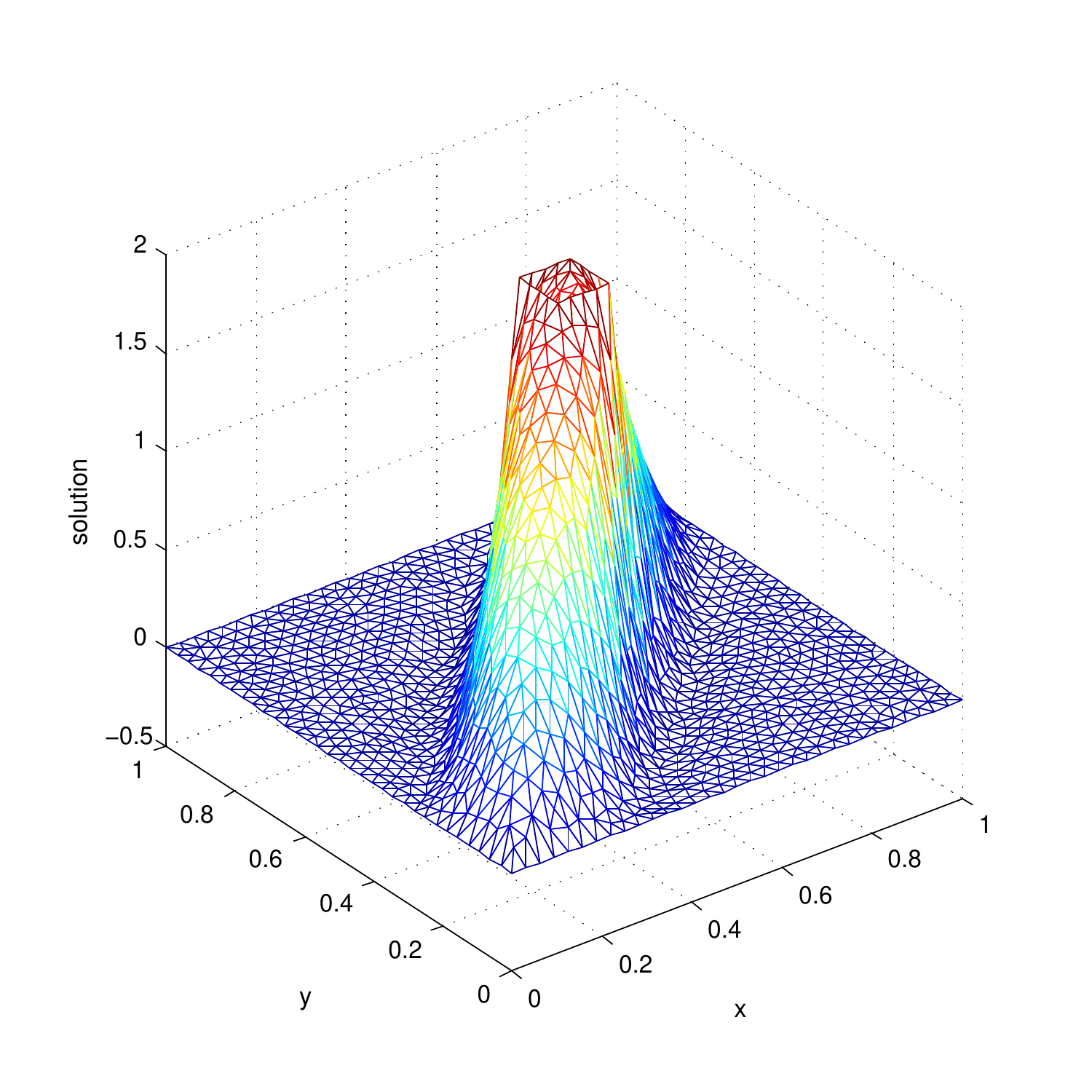}
\centerline{(a): $M_{unif}$}
\end{minipage}
\begin{minipage}[b]{2.5in}
\includegraphics[width=2.5in]{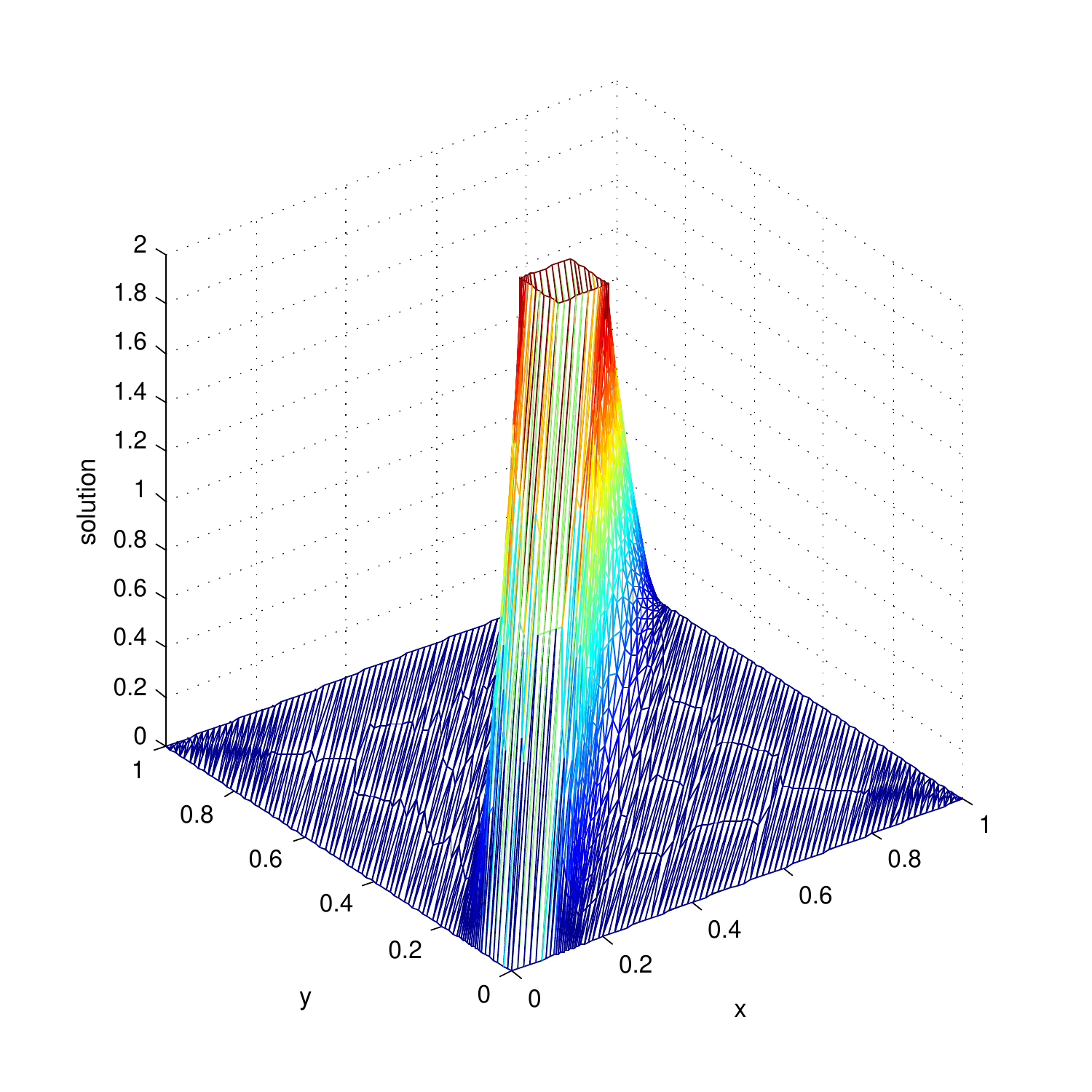}
\centerline{(b): $M_{DMP+adap}$}
\end{minipage}
}
\caption{Example~\ref{ex1} with constant $\mathbb{D}$. Finite element solutions obtained with
(a) $M_{unif}$ and (b) $M_{DMP+adap}$.}
\label{ex1-soln}
\end{figure}

\begin{figure}[thb]
\centering
\hbox{
\hspace{10mm}
\begin{minipage}[t]{2.5in}
\includegraphics[width=2.5in]{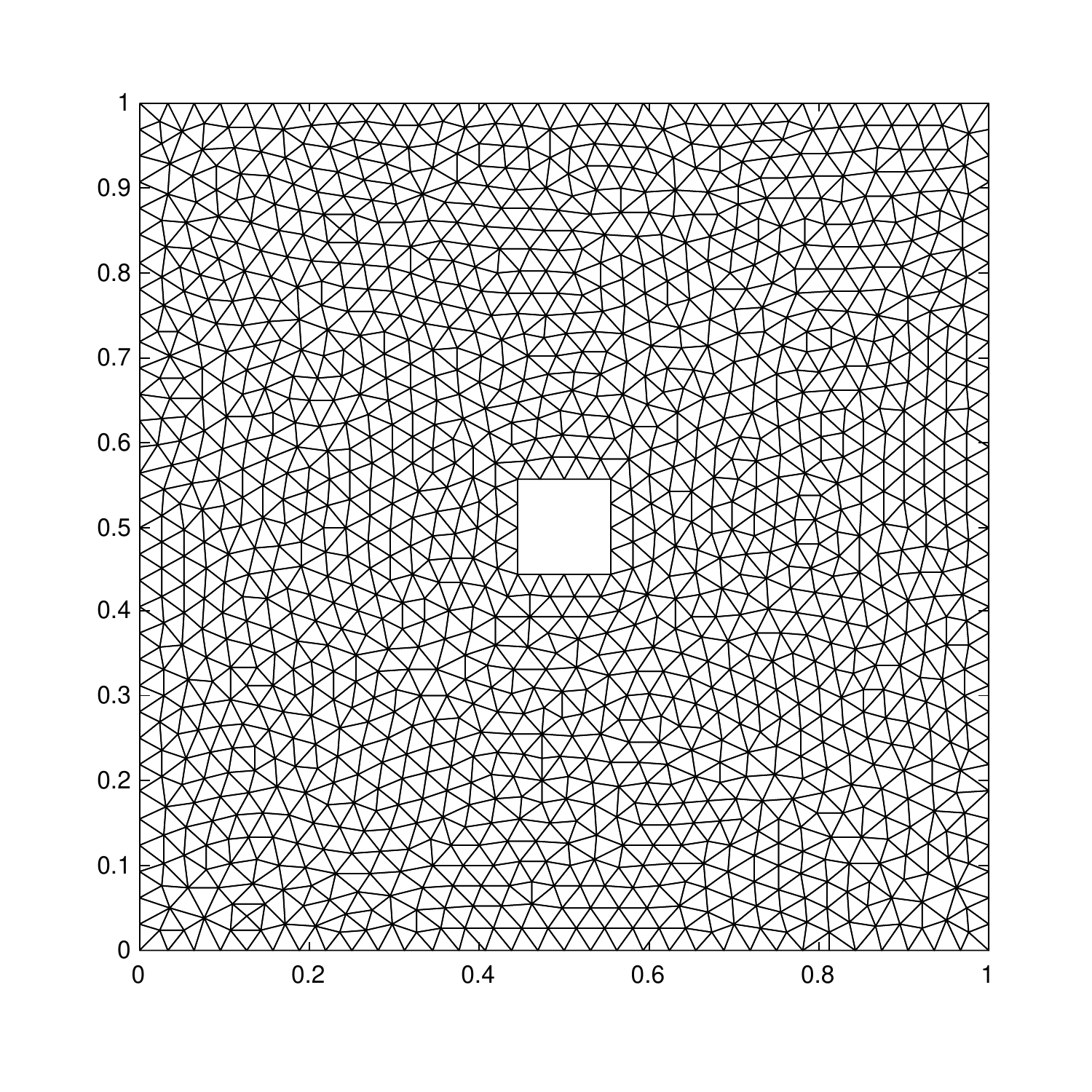}
\centerline{(a): $M_{unif}$, $N=2460$}
\end{minipage}
\hspace{10mm}
\begin{minipage}[t]{2.5in}
\includegraphics[width=2.5in]{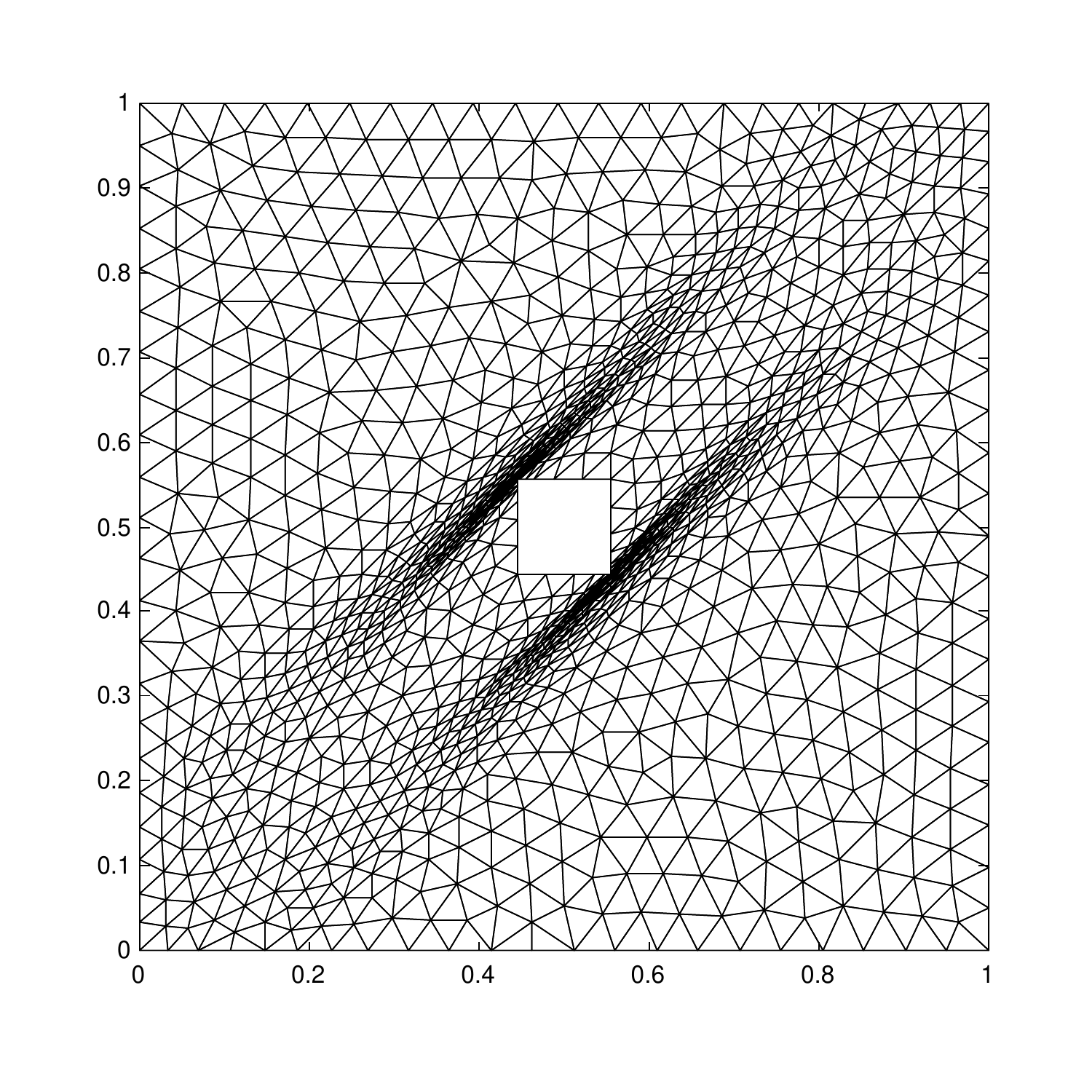}
\centerline{(b): $M_{adap}$, $N=2583$}
\end{minipage}
}
\hbox{
\hspace{10mm}
\begin{minipage}[t]{2.5in}
\includegraphics[width=2.5in]{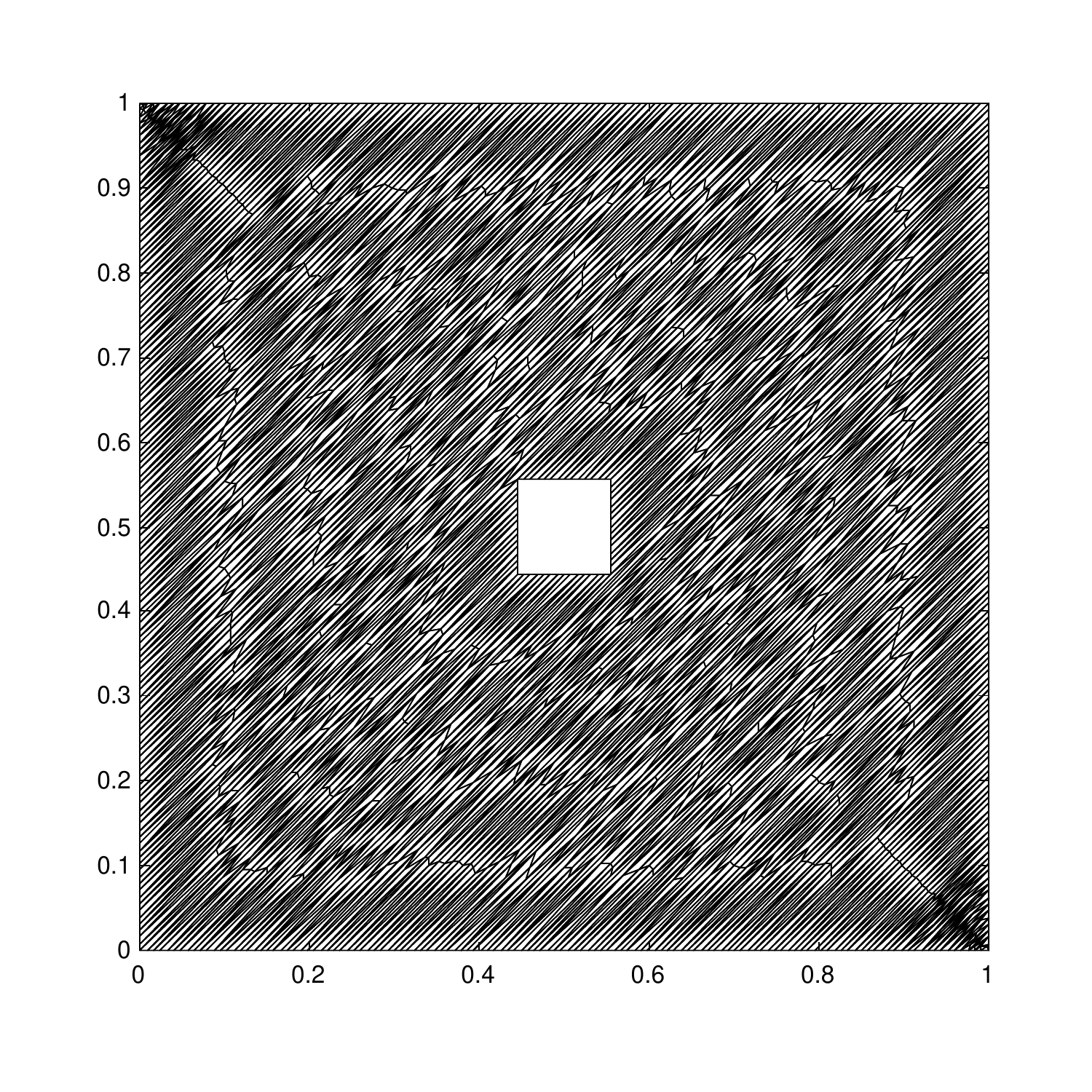}
\centerline{(c): $M_{DMP}$, $N=2530$}
\end{minipage}
\hspace{10mm}
\begin{minipage}[t]{2.5in}
\includegraphics[width=2.5in]{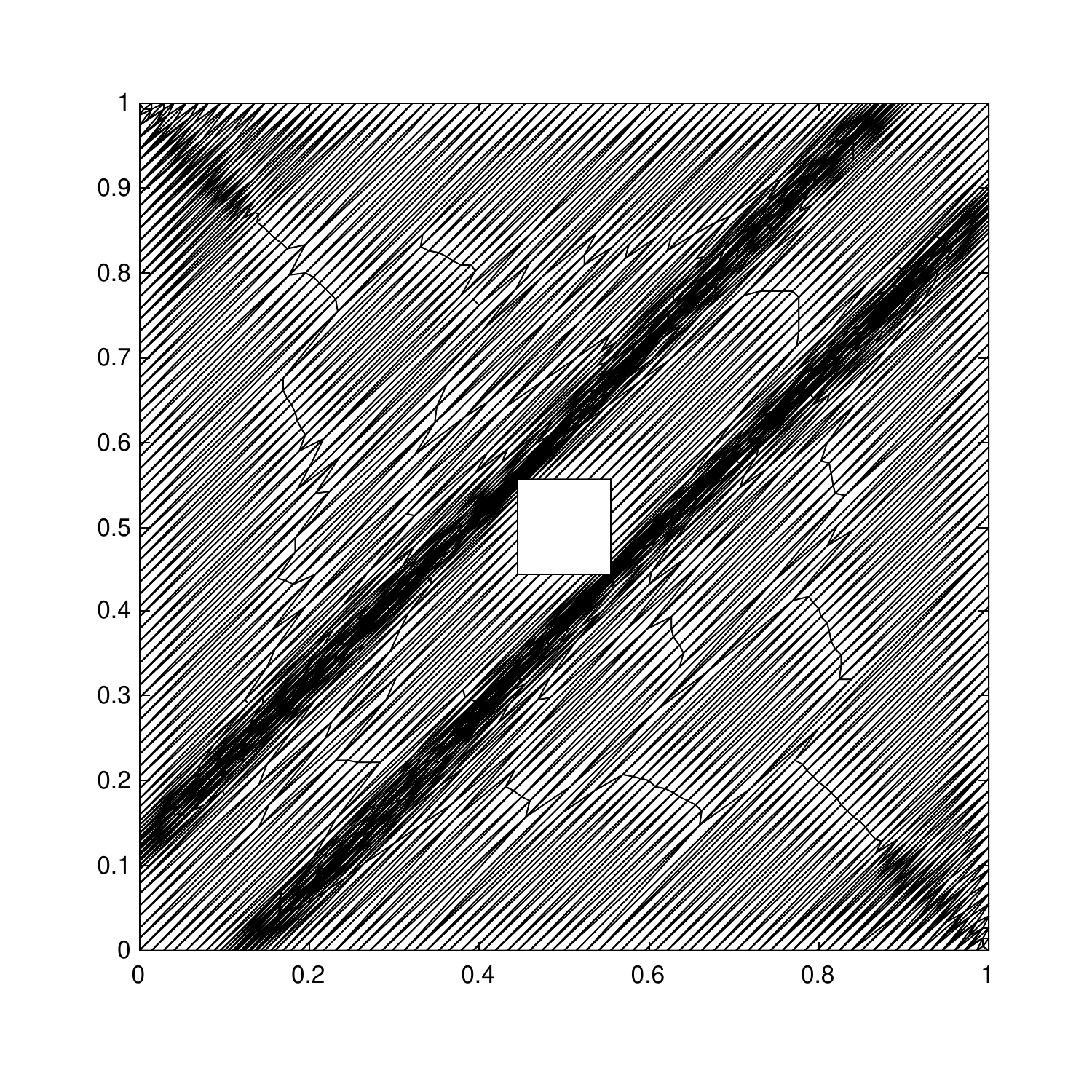}
\centerline{(d): $M_{DMP+adap}$, $N=2474$}
\end{minipage}
}
\caption{Example~\ref{ex1} with constant $\mathbb{D}$. Meshes obtained with different metric tensors.}
\label{ex1-mesh}
\end{figure}

\begin{figure}[thb]
\centering
\hbox{
\hspace{10mm}
\begin{minipage}[t]{2.5in}
\includegraphics[width=2.5in]{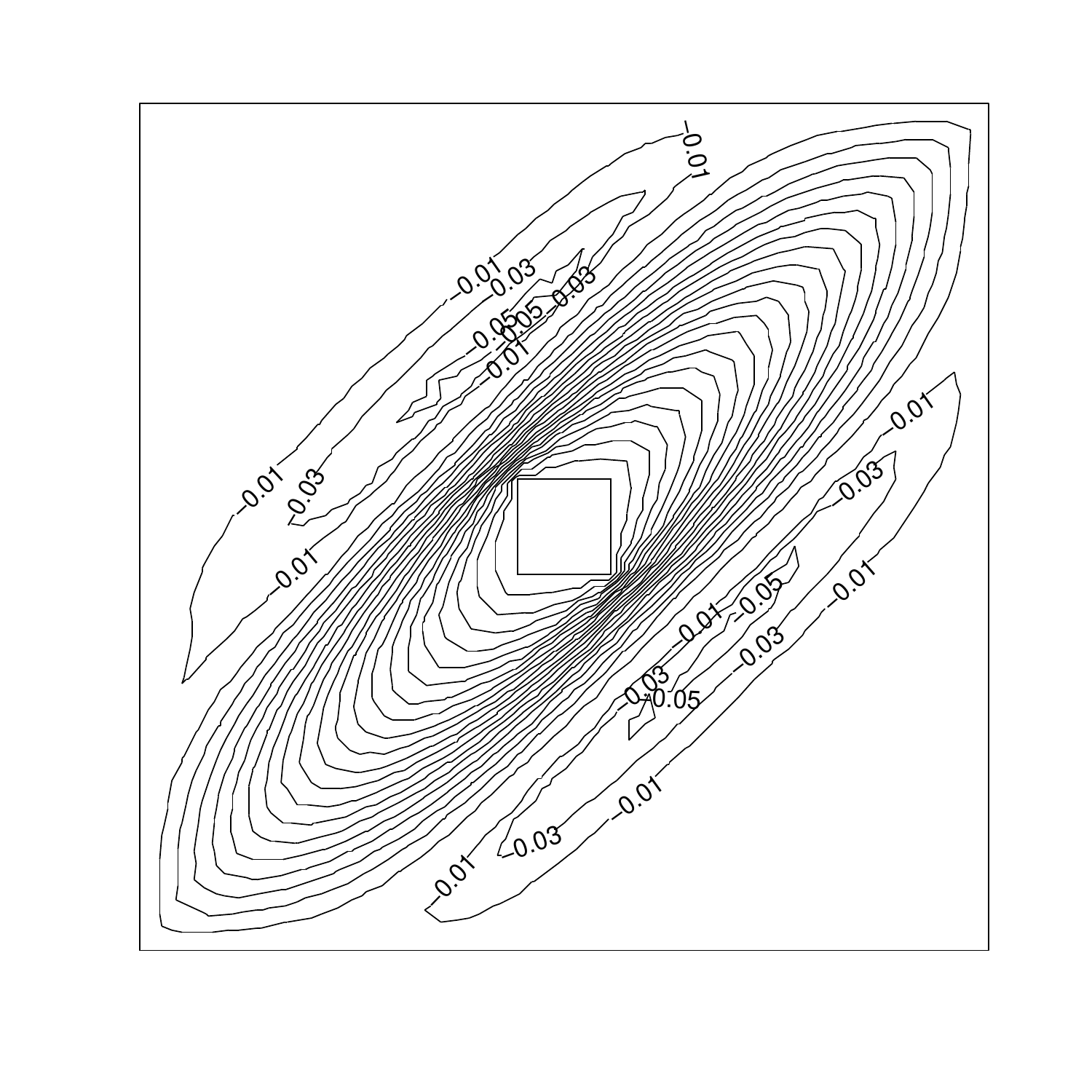}
\centerline{(a): $M_{unif}$, $u_{min}=-0.0602$}
\end{minipage}
\hspace{10mm}
\begin{minipage}[t]{2.5in}
\includegraphics[width=2.5in]{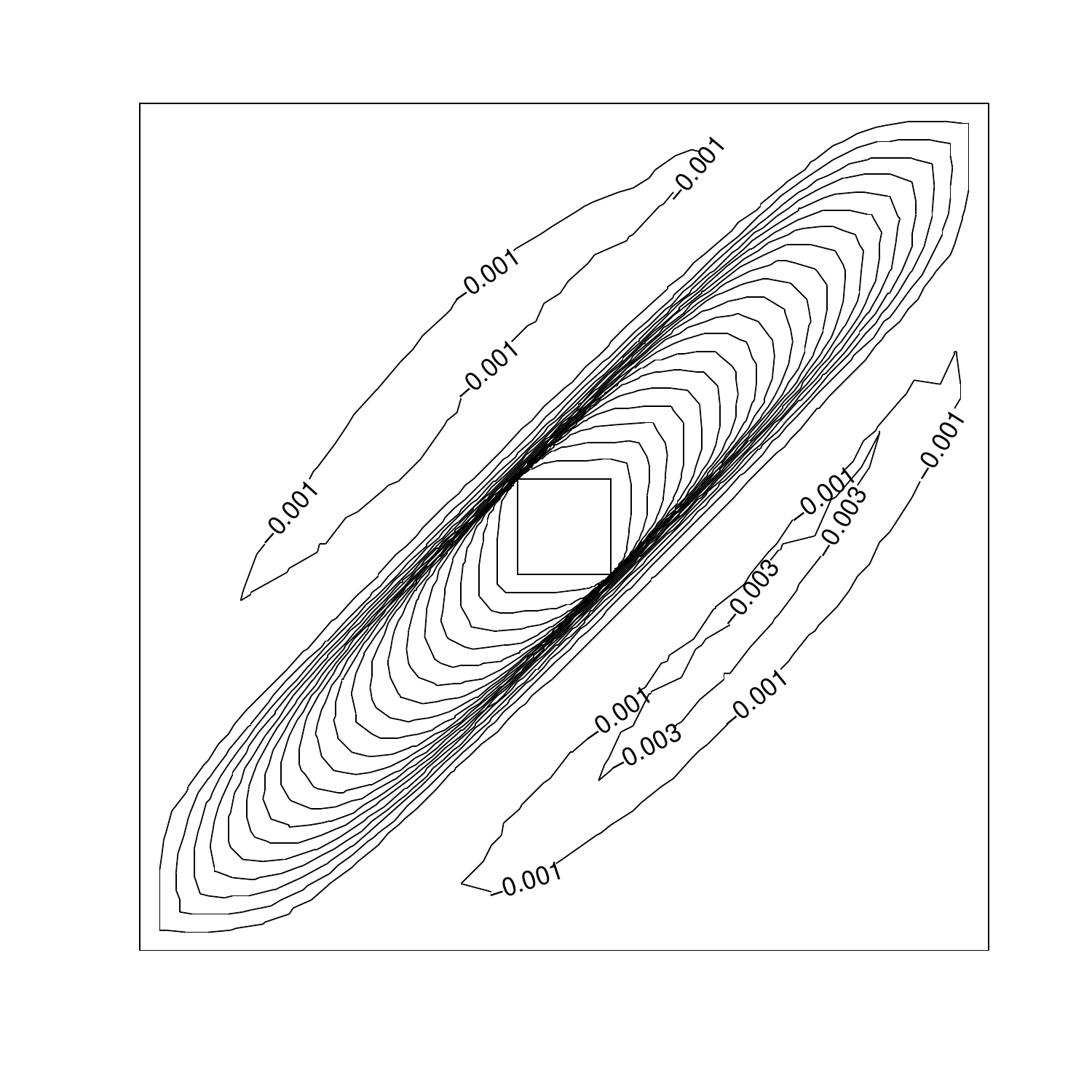}
\centerline{(b): $M_{adap}$, $u_{min}=-0.0039$}
\end{minipage}
}
\hbox{
\hspace{10mm}
\begin{minipage}[t]{2.5in}
\includegraphics[width=2.5in]{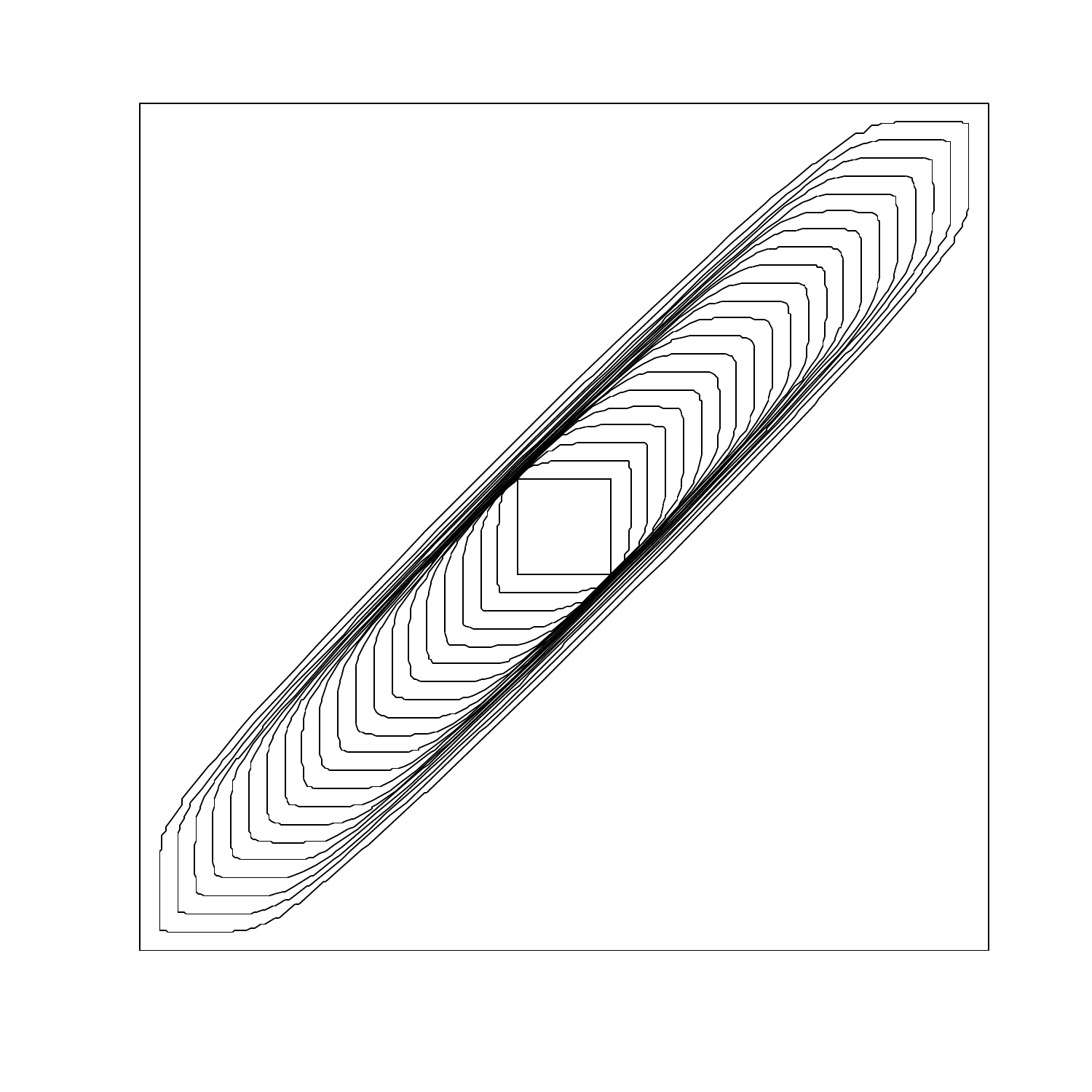}
\centerline{(c): $M_{DMP}$, $u_{min}=0$}
\end{minipage}
\hspace{10mm}
\begin{minipage}[t]{2.5in}
\includegraphics[width=2.5in]{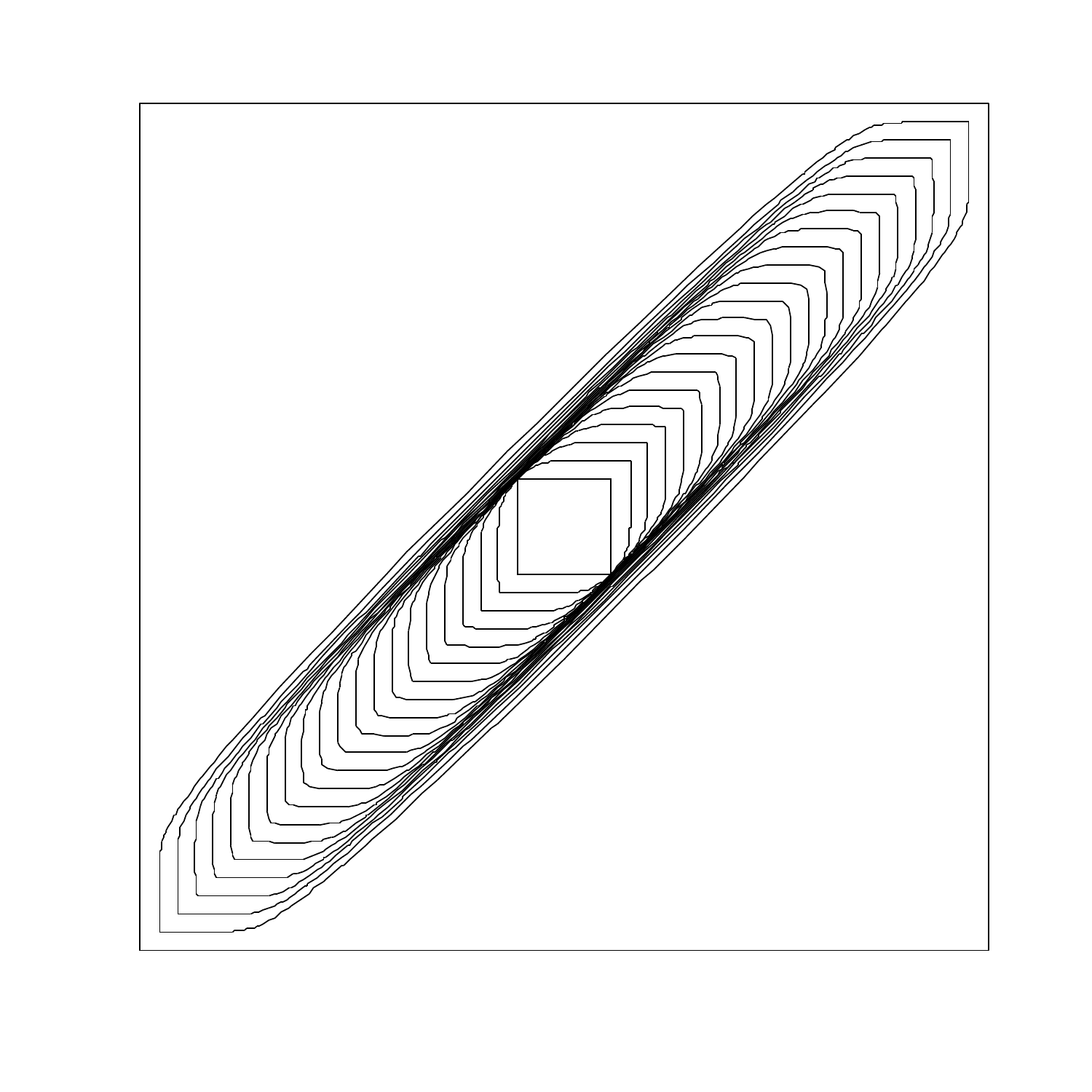}
\centerline{(d): $M_{DMP+adap}$, $u_{min}=0$}
\end{minipage}
}
\caption{Example~\ref{ex1} with constant $\mathbb{D}$. Contours of the finite element solutions obtained with
different metric tensors.}
\label{ex1-contour}
\end{figure}

\begin{figure}[thb]
\centering
\includegraphics[width=3in]{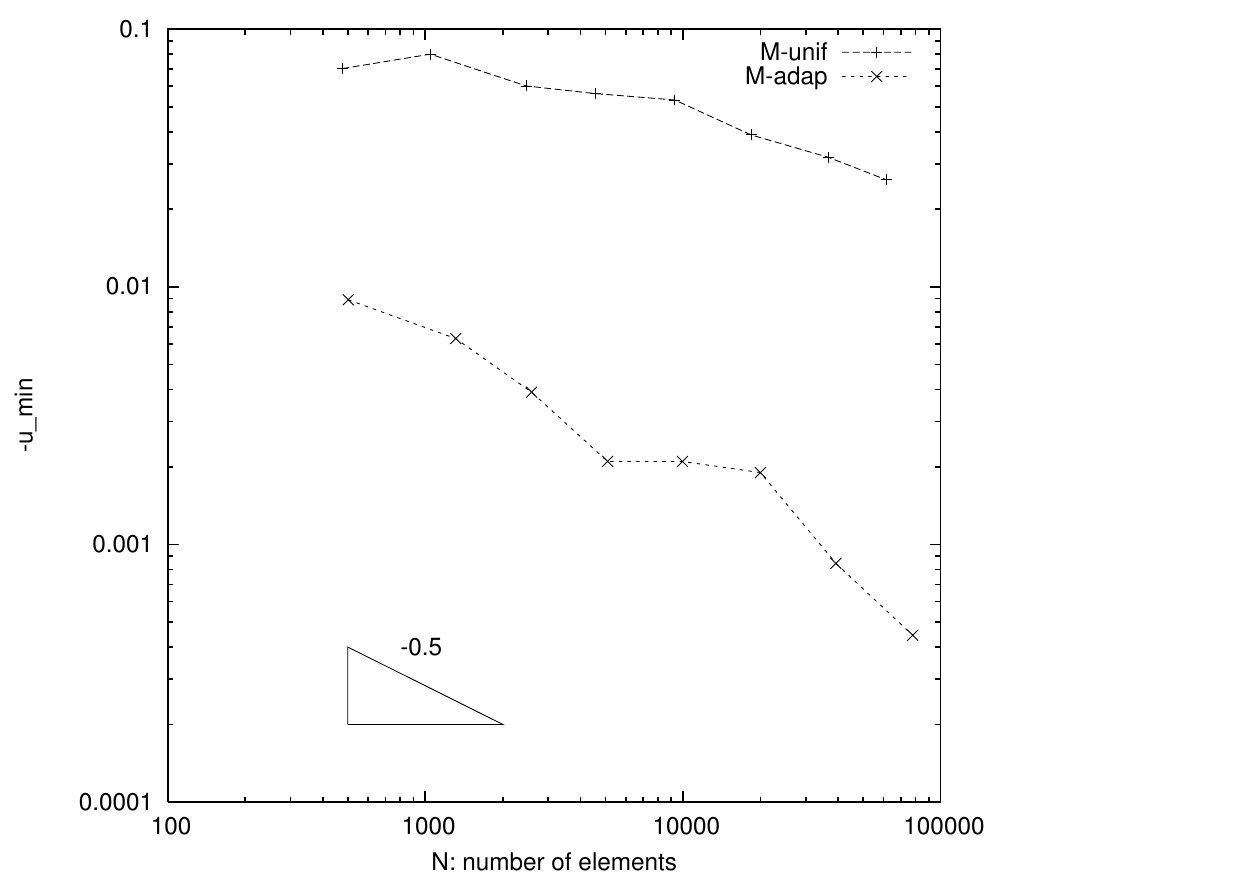}
\caption{Example~\ref{ex1} with constant $\mathbb{D}$. The undershoot, $-u_{min}$, is shown as functions
of the number of elements.}
\label{ex1-umin}
\end{figure}

\begin{figure}[thb]
\centering
\hbox{
\hspace{10mm}
\begin{minipage}[b]{2.5in}
\includegraphics[width=2.5in]{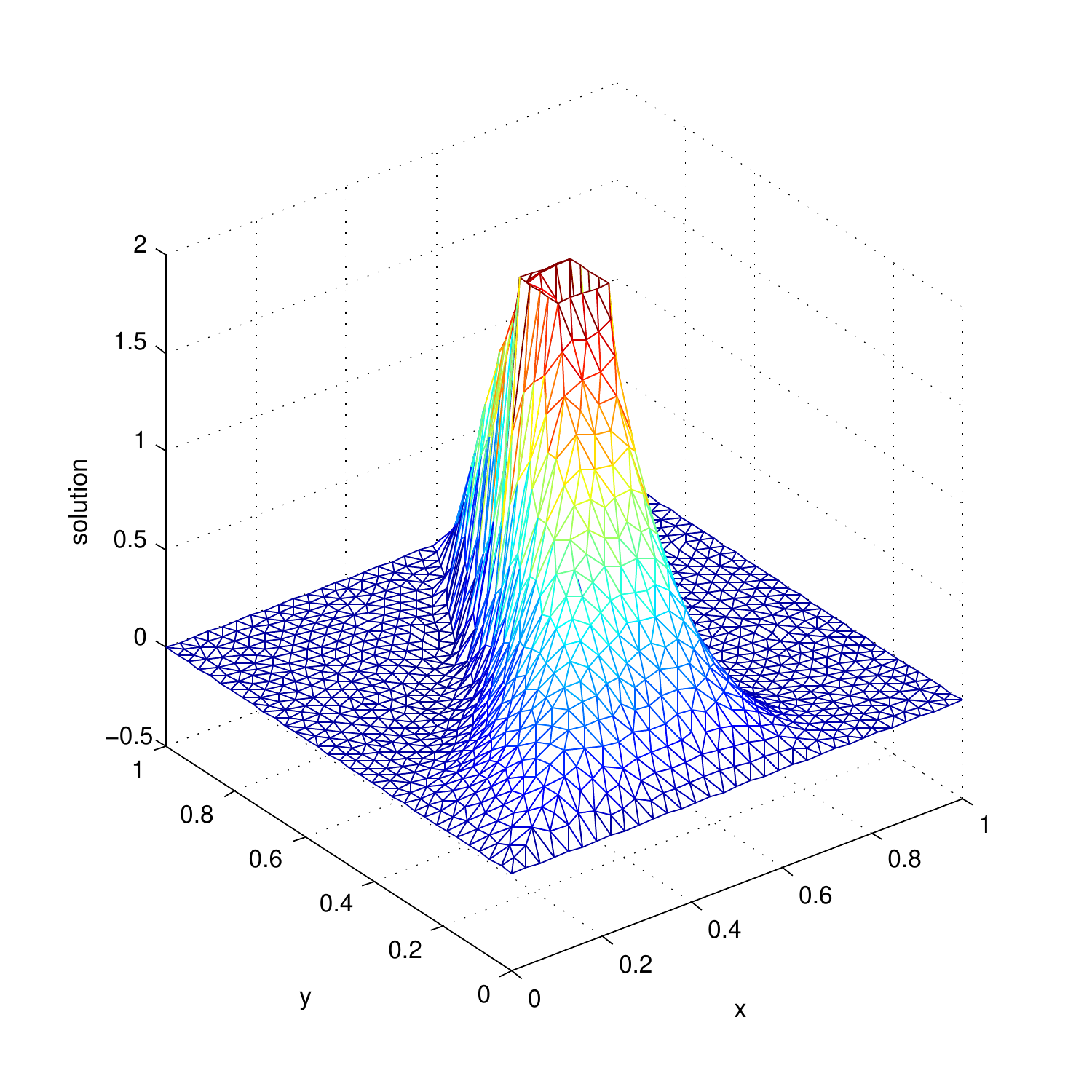}
\centerline{(a): $M_{unif}$}
\end{minipage}\begin{minipage}[b]{2.5in}
\includegraphics[width=2.5in]{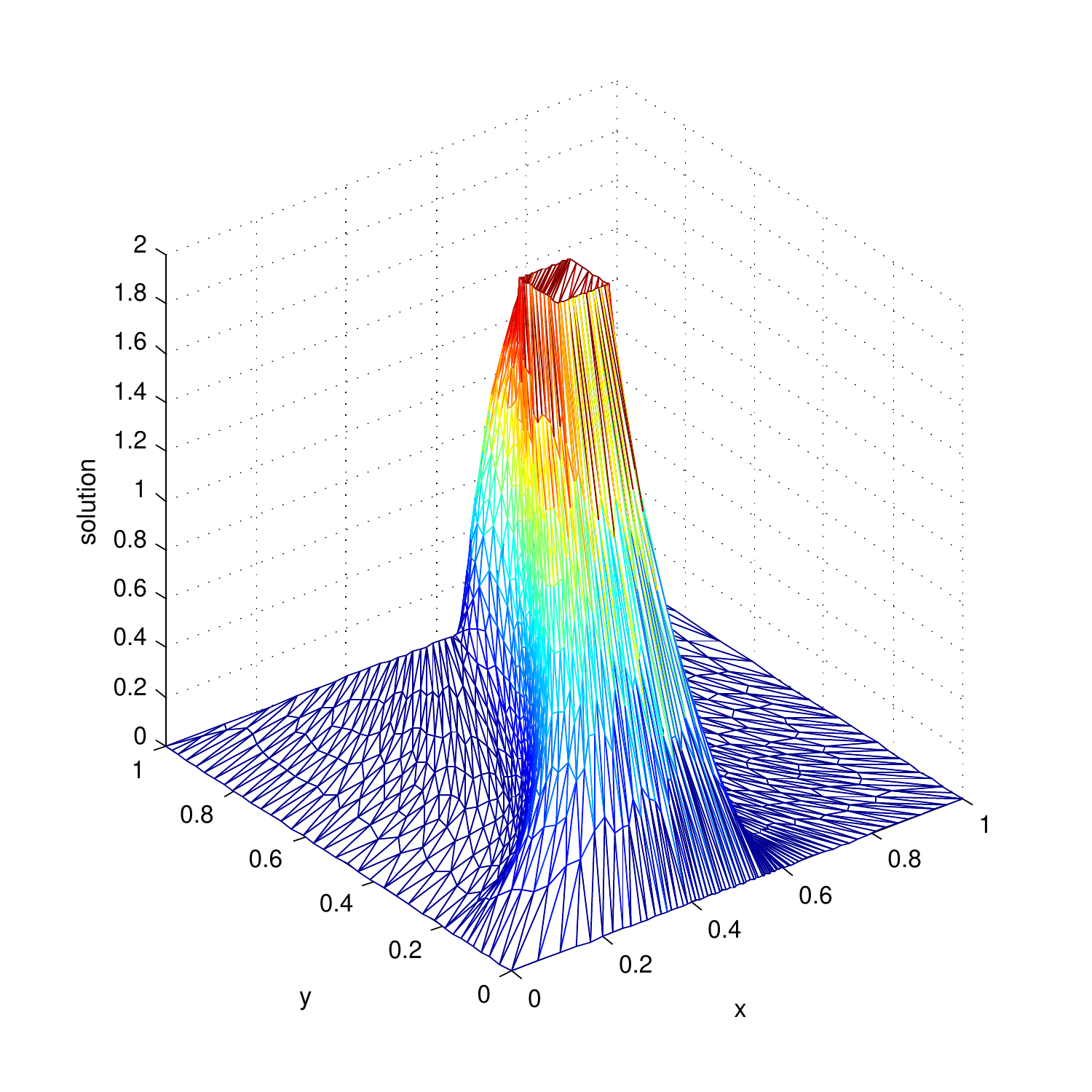}
\centerline{(b): $M_{DMP+adap}$}
\end{minipage}
}
\caption{Example~\ref{ex1} with variable $\mathbb{D}$. Finite element solutions obtained with
(a) $M_{unif}$ and (b) $M_{DMP+adap}$.}
\label{ex1b-soln}
\end{figure}

\begin{figure}[thb]
\centering
\hbox{
\begin{minipage}[t]{3in}
\includegraphics[width=2.5in]{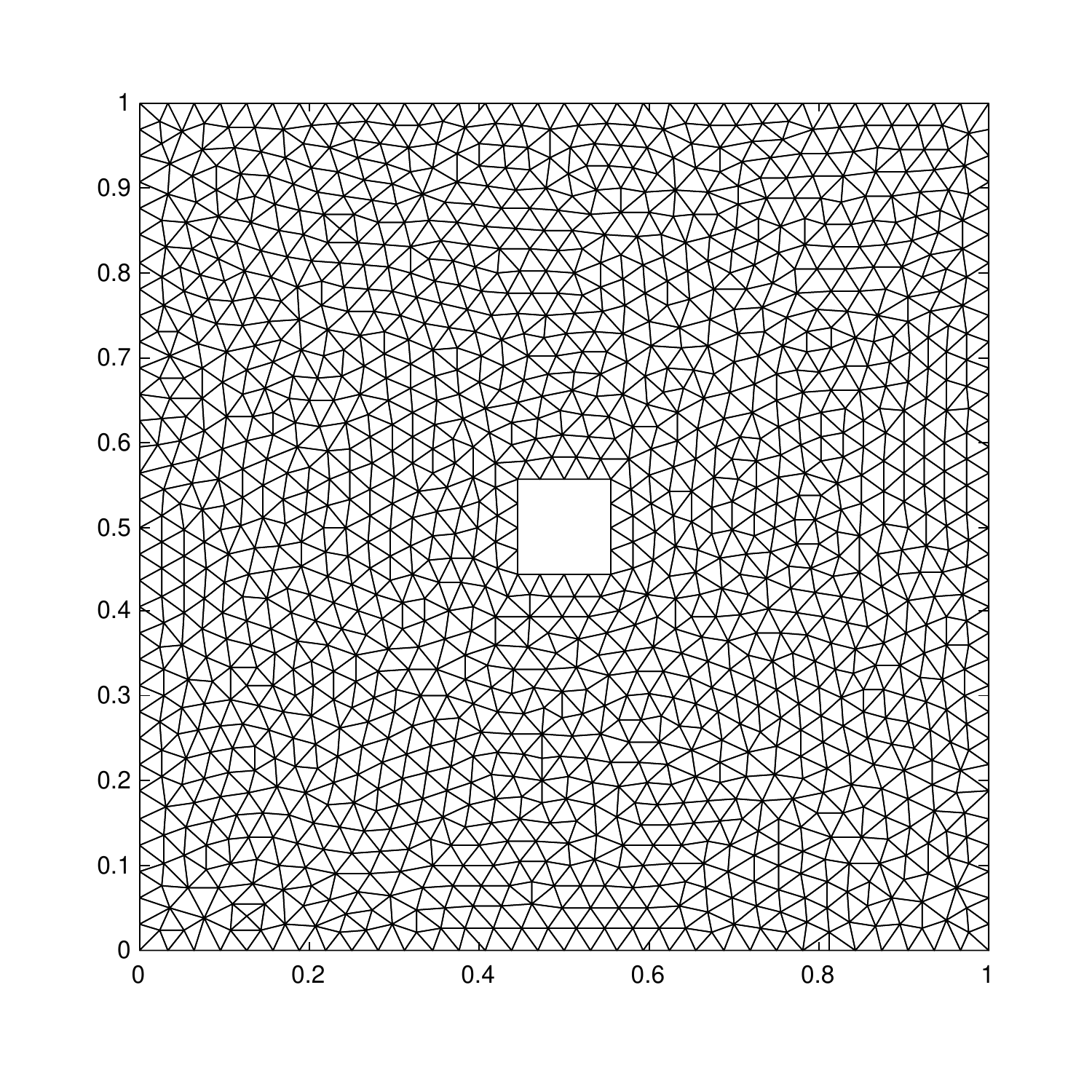}
\centerline{(a):  $M_{unif}$, $N=2460$}
\end{minipage}
\hspace{10mm}
\begin{minipage}[t]{3in}
\includegraphics[width=2.5in]{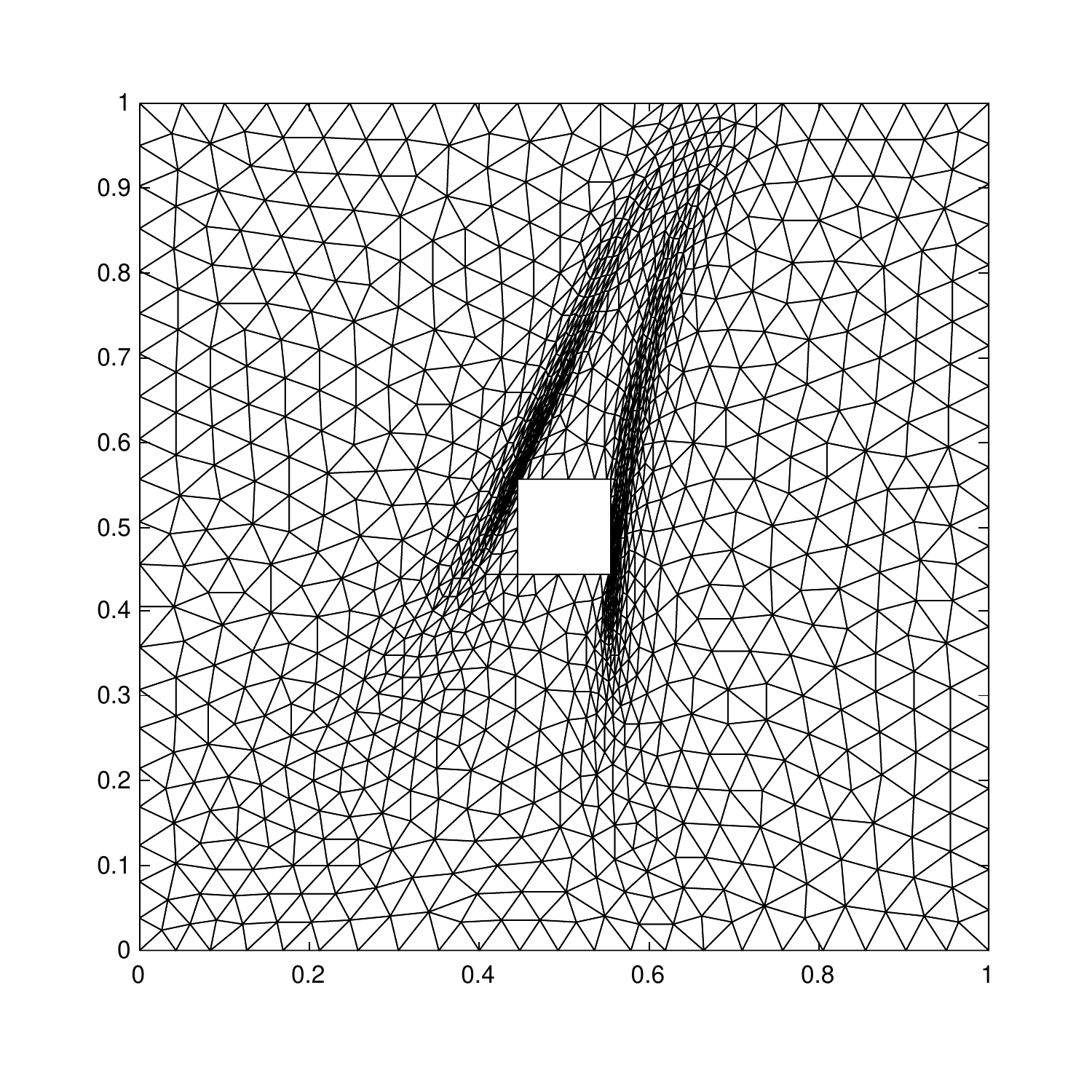}
\centerline{(b): $M_{adap}$, $N=2568$}
\end{minipage}
}
\hbox{
\begin{minipage}[t]{3in}
\includegraphics[width=2.5in]{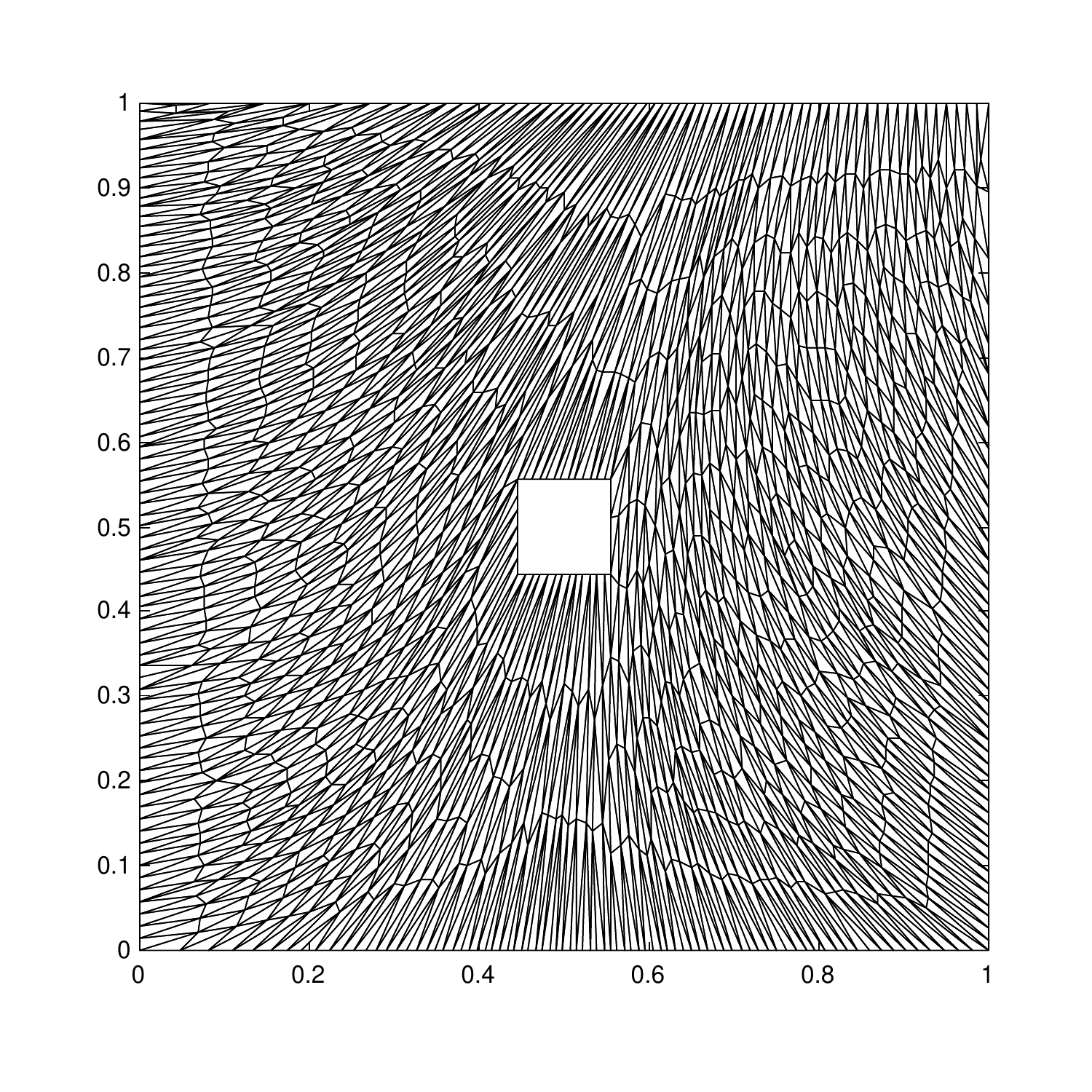}
\centerline{(c): $M_{DMP}$, $N=2510$}
\end{minipage}
\hspace{10mm}
\begin{minipage}[t]{3in}
\includegraphics[width=2.5in]{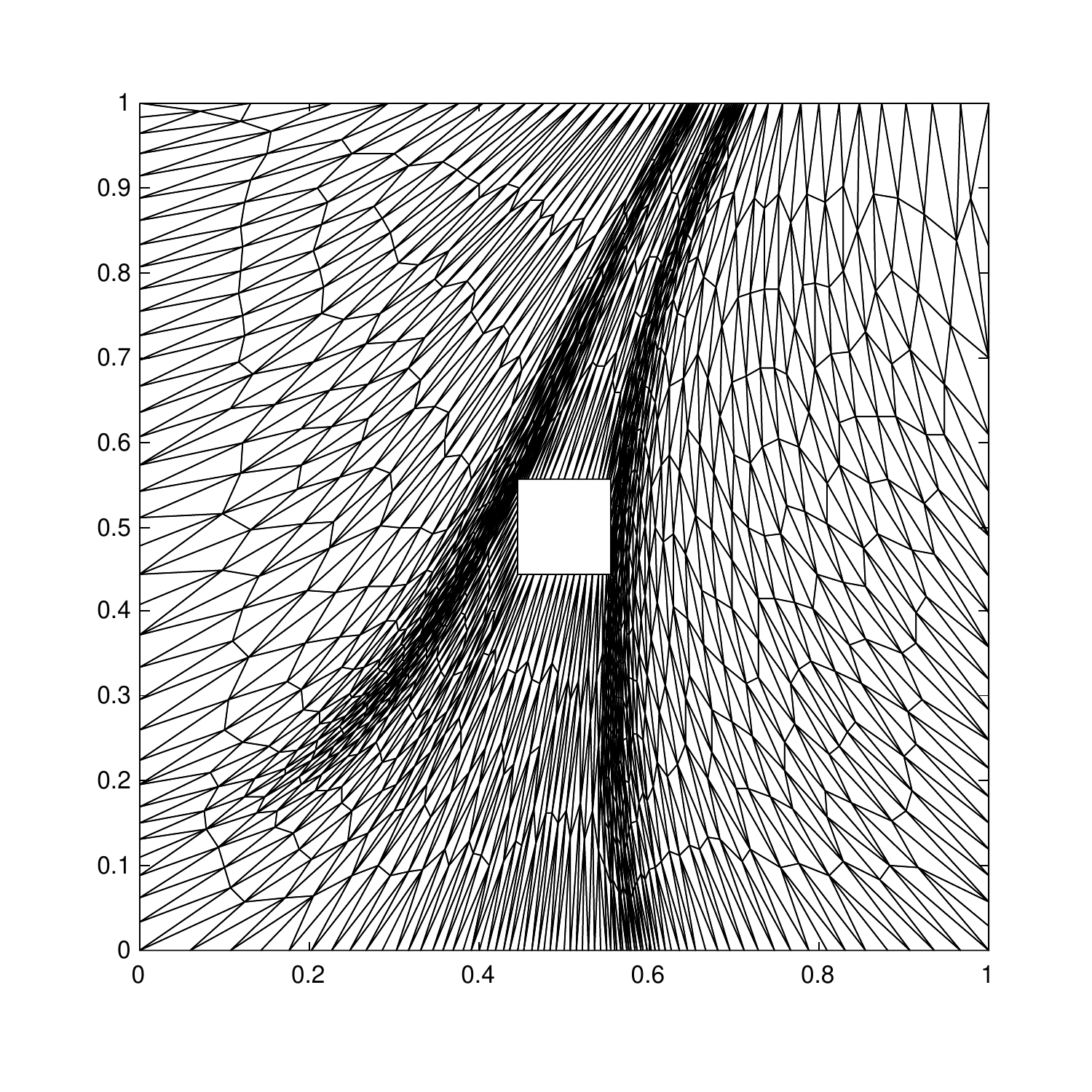}
\centerline{(d): $M_{DMP+adap}$, $N=2463$}
\end{minipage}
}
\caption{Example~\ref{ex1} with variable $\mathbb{D}$. Meshes obtained from different metric tensors.}
\label{ex1b-mesh}
\end{figure}

\begin{figure}[thb]
\centering
\hbox{
\hspace{10mm}
\begin{minipage}[t]{2.5in}
\includegraphics[width=2.5in]{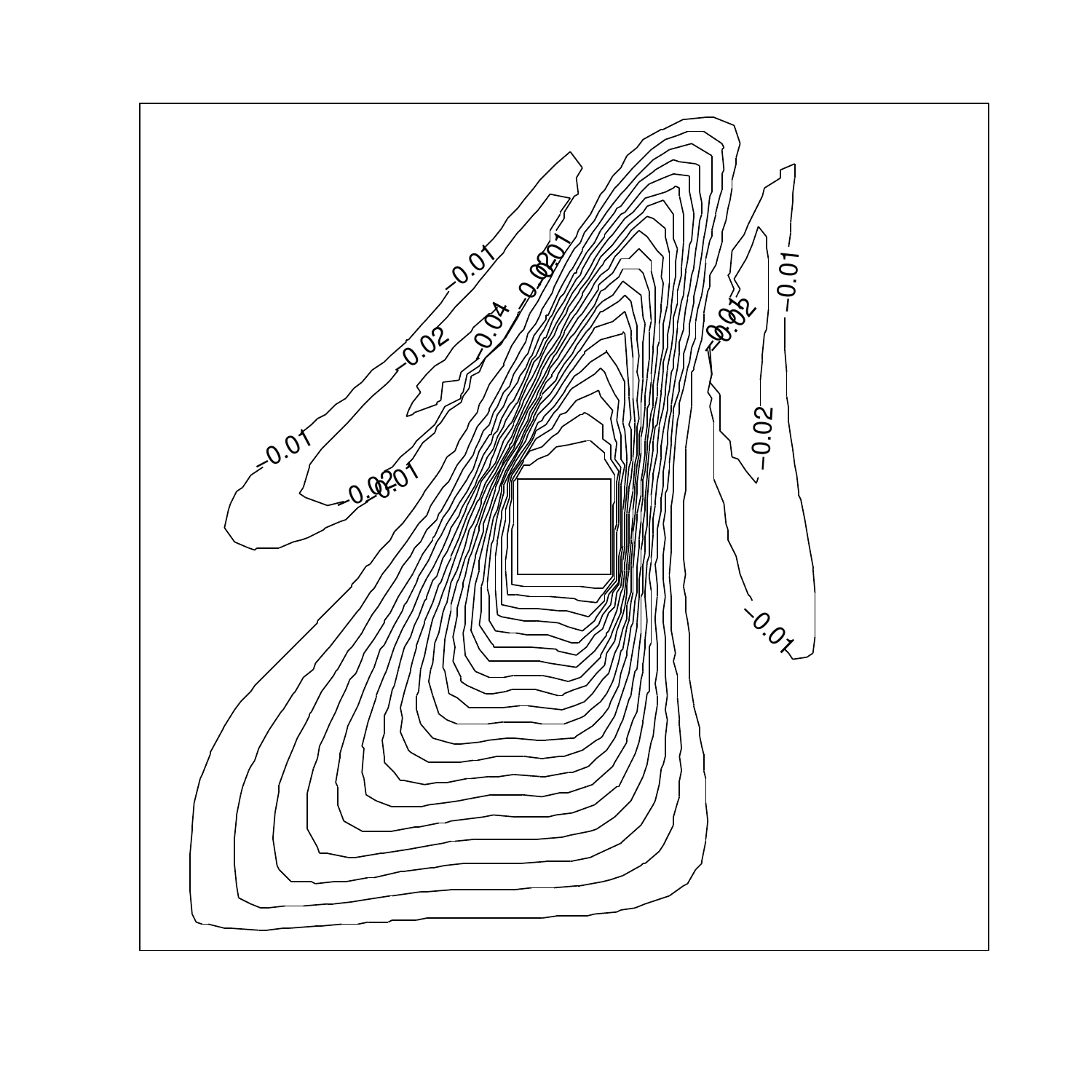}
\centerline{(a): $M_{unif}$, $u_{min}=-0.0506$}
\end{minipage}
\hspace{10mm}
\begin{minipage}[t]{2.5in}
\includegraphics[width=2.5in]{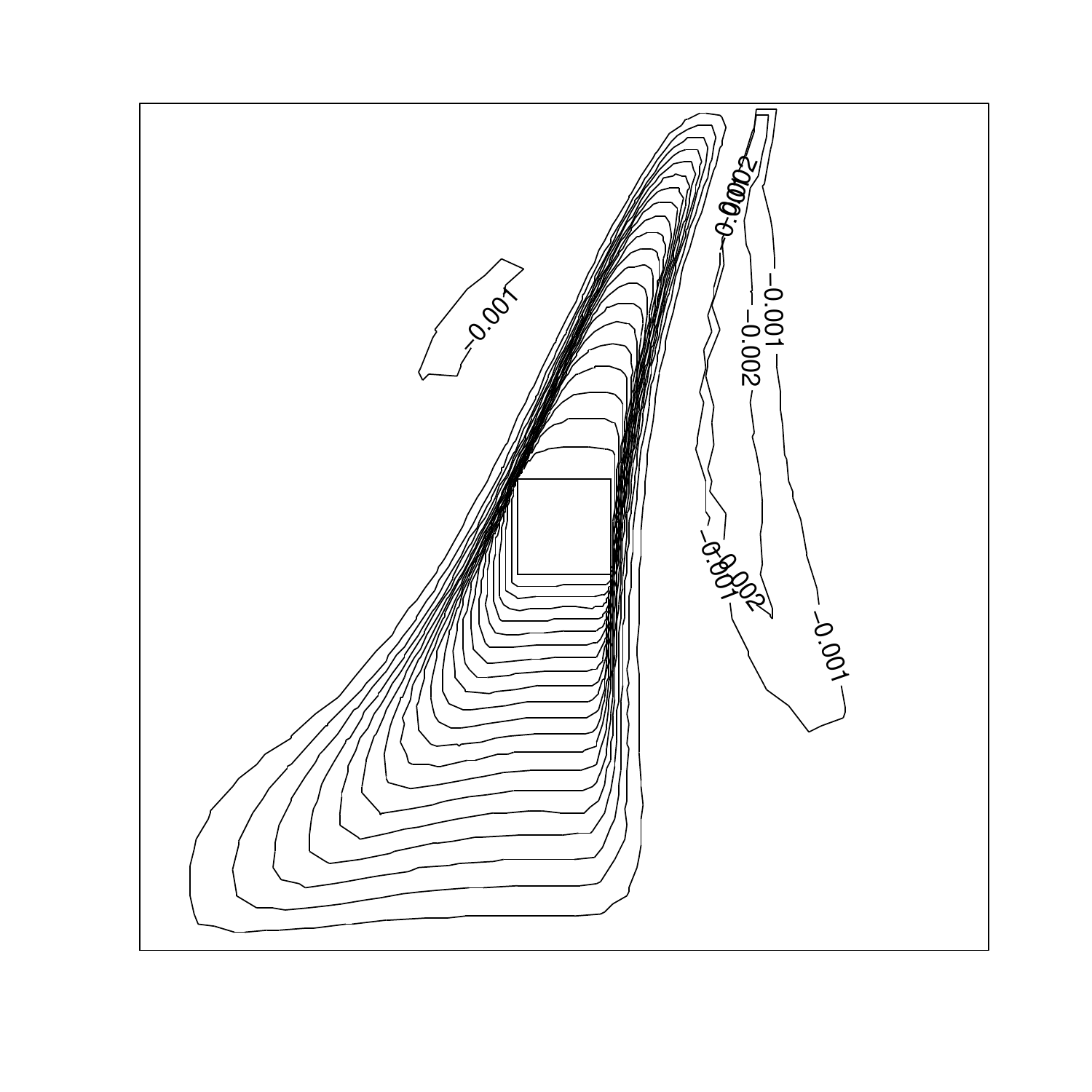}
\centerline{(b): $M_{adap}$, $u_{min}=-0.0039$}
\end{minipage}
}
\hbox{
\hspace{10mm}
\begin{minipage}[t]{2.5in}
\includegraphics[width=2.5in]{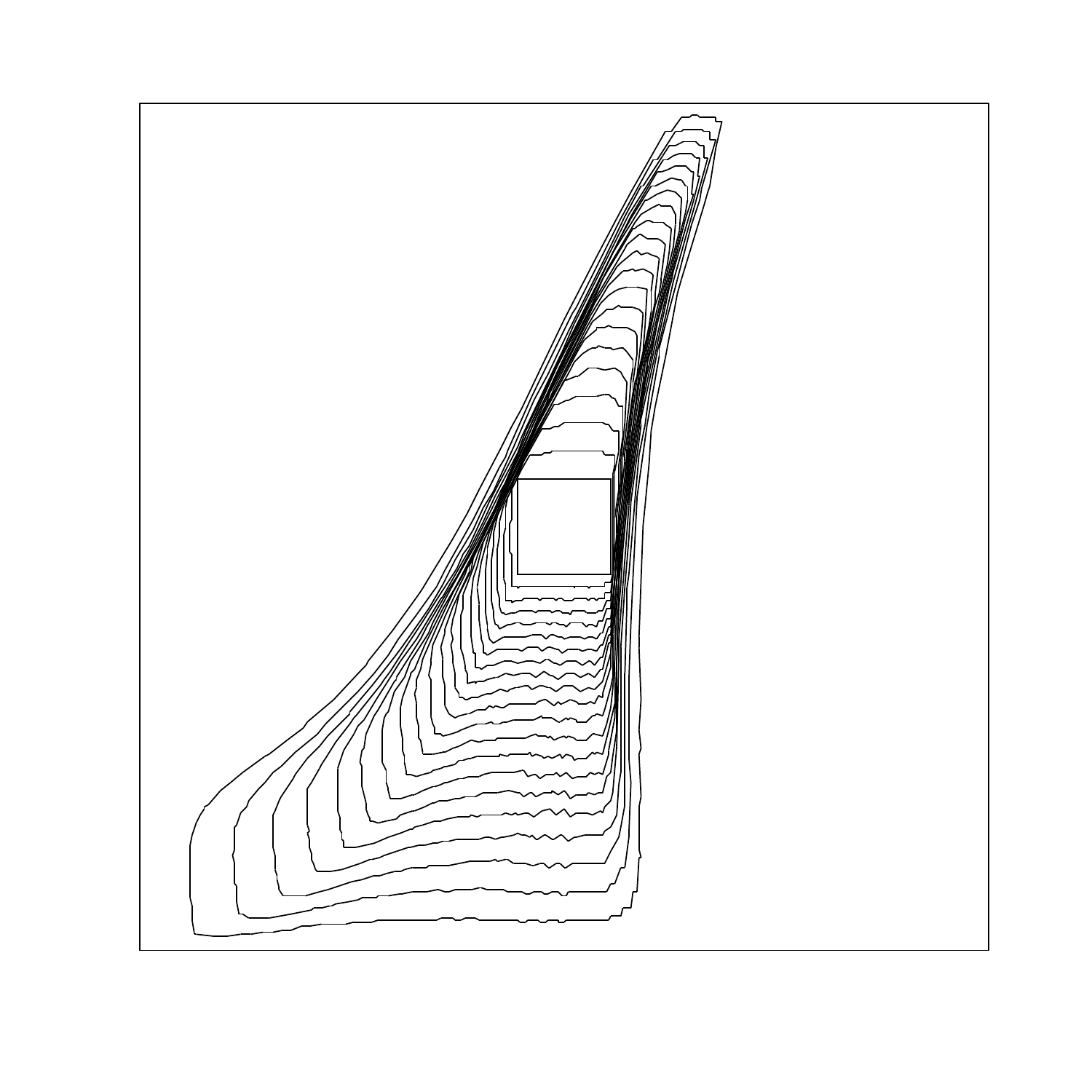}
\centerline{(c): $M_{DMP}$, $u_{min}=0$}
\end{minipage}
\hspace{10mm}
\begin{minipage}[t]{2.5in}
\includegraphics[width=2.5in]{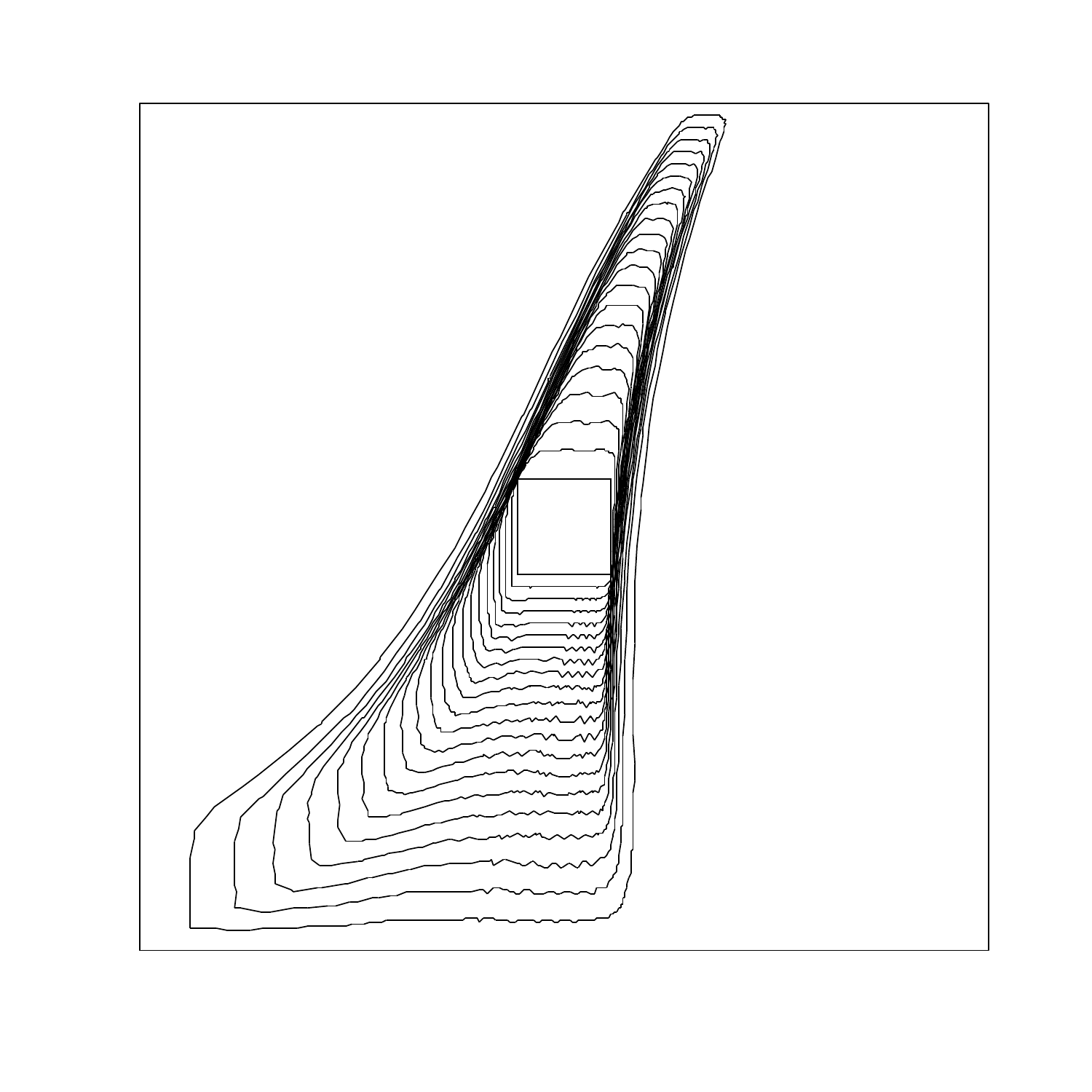}
\centerline{(d): $M_{DMP+adap}$, $u_{min}=0$}
\end{minipage}
}
\caption{Example~\ref{ex1} with variable $\mathbb{D}$. Contours of the finite element solutions obtained with
different metric tensors.}
\label{ex1b-contour}
\end{figure}

\begin{exam}
\label{ex1}
{\em
The first example is to consider BVP (\ref{bvp-pde}) and (\ref{bvp-bc}) with
\[
f\equiv 0,\quad \Omega = [0,1]^2\backslash\left [\frac{4}{9},\frac{5}{9}\right ]^2,
\quad g = 0 \mbox{ on } \Gamma_{out}, \quad  g = 2 \mbox{ on } \Gamma_{in},
\]
where $\Gamma_{out}$ and $\Gamma_{in}$ are the outer and inner boundaries of $\Omega$, respectively;
see Fig. \ref{ex1-domain}. The diffusion matrix is given by (\ref{D-eigen-d}) with $k_1 = 1000$, $k_2 = 1$,
and $\theta$ being the angle of the primary diffusion direction (parallel to the first eigenvector of $\mathbb{D}$).

This example satisfies the maximum principle and the solution (whose analytical expression is unavailable)
stays between 0 and 2. Our goal is to produce a numerical solution which also satisfies DMP and stays
between 0 and 2. Moreover, for both cases with a constant and a variable $\theta$ we consider,
the exact solution has sharp jumps near the inner boundary
(cf. Figs. \ref{ex1-soln} and \ref{ex1b-soln}) so mesh adaptation
is needed for a proper resolution of them.
This example has been studied in \cite{KSS09,LSS07}.

We first consider the case of constant $\mathbb{D}$ with $\theta = \pi/4$. Fig. \ref{ex1-soln}
shows finite element solutions obtained with $M_{unif}$ and $M_{DMP+adap}$.
Meshes and solution contours obtained with various
metric tensors are shown in Figs. \ref{ex1-mesh} and \ref{ex1-contour}, respectively.
No overshoots in the finite element solutions are observed for all cases. However,
undershoots and unphysical minima occur in the solutions obtained with $M_{unif}$ ($u_{min}=-0.0602$)
and $M_{adap}$ ($u_{min}=-0.0039$) (cf. Fig. \ref{ex1-contour})(a) and (b)).
Fig.~\ref{ex1-umin} shows the decrease of $-u_{min}$
as the mesh is refined. For the range of the number of mesh elements considered, the undershooting
improves at a rate of $-u_{min} = O( N^{-0.5})$ for both $M_{unif}$ and $M_{adap}$.
On the other hand, the results confirm the theoretical prediction that
the solutions obtained with $M_{DMP}$ and $M_{DMP+adap}$ satisfy DMP and no overshoot/undershoot
and no unphysical extremum occur. It should be pointed out that the solution contour obtained
with an almost uniform mesh is very smooth but the sharp jumps of the solution are smeared;
see Figs. \ref{ex1-soln}(a) and \ref{ex1-contour}(a). The solution contours obtained with $M_{DMP}$ and $M_{DMP+adap}$ are comparable to the one obtained with $M_{adap}$. 

Next we consider a case of variable $\mathbb{D}$ with $\theta= \pi \sin(x) \cos(y)$.
The finite element solutions, meshes, and solution contours are shown in Figs. \ref{ex1b-soln},
\ref{ex1b-mesh}, and \ref{ex1b-contour}, respectively. Similar observations as for the constant $\mathbb{D}$
case can be made. Especially, undershoots and unphysical extrema occur in the solutions obtained with
$M_{unif}$ and $M_{adap}$ but not with $M_{DMP}$ and $M_{DMP+adap}$. Once again, the results
confirm our theoretical predictions in the previous sections.
}\end{exam}

\begin{figure}[thb]
\centering
\hbox{
\hspace{10mm}
\begin{minipage}[b]{2.5in}
\includegraphics[width=2.5in]{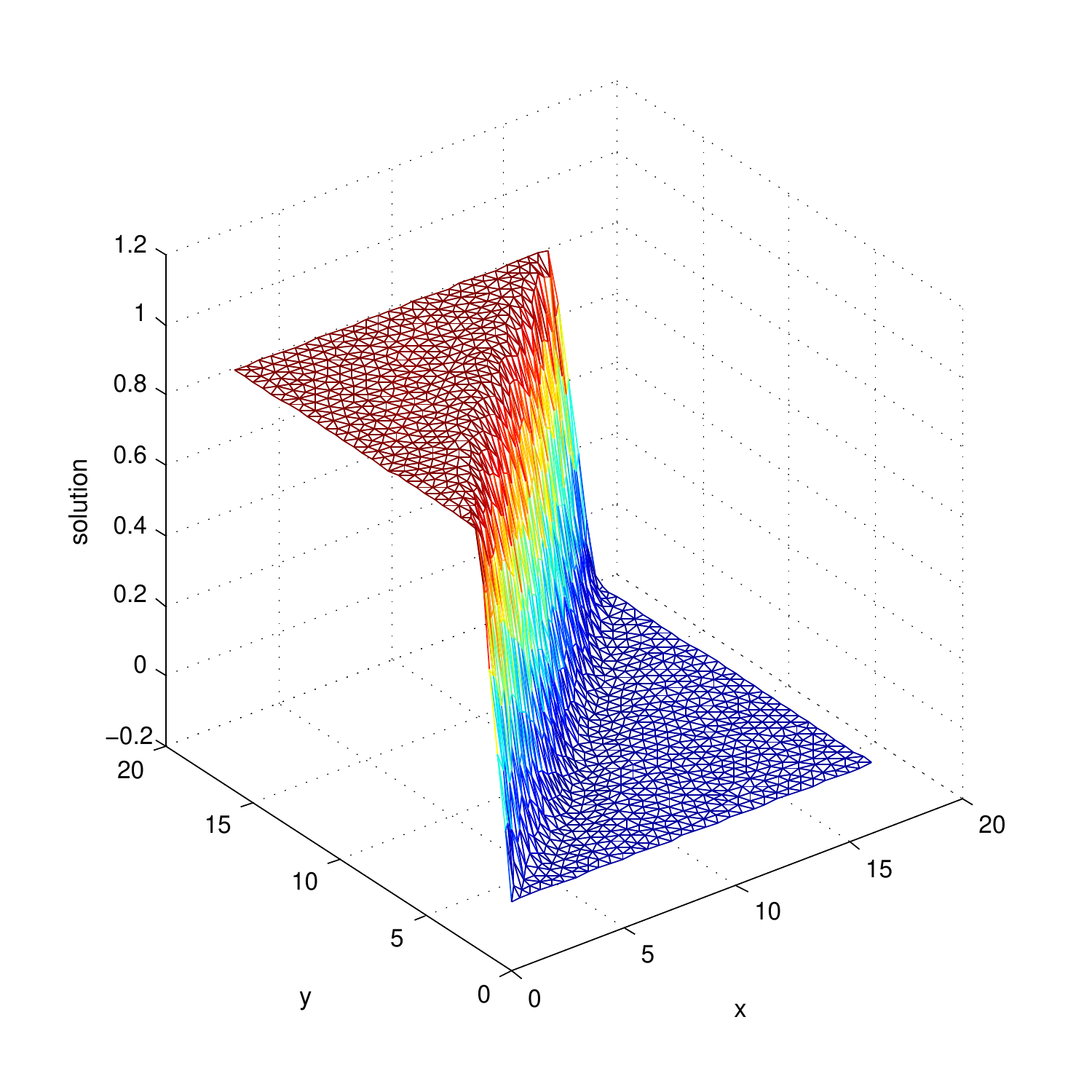}
\centerline{(a): $M_{adap}$, $u_{min}=-0.0195$}
\end{minipage}
\begin{minipage}[b]{2.5in}
\includegraphics[width=2.5in]{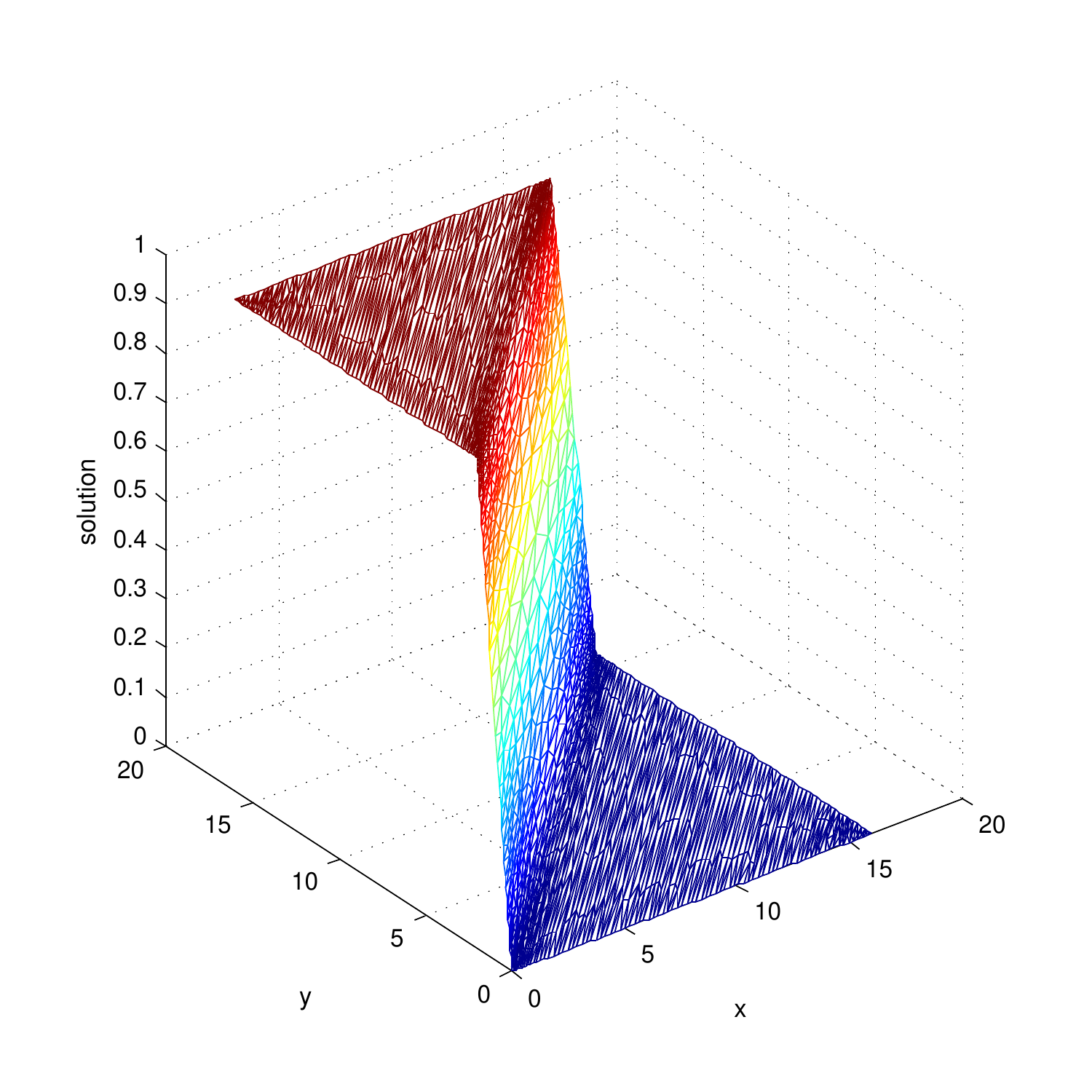}
\centerline{(b): $M_{DMP+adap}$, $u_{min}=0$}
\end{minipage}
}
\caption{Example~\ref{ex2}. Finite element solutions obtained with
(a) $M_{adap}$ and (b) $M_{DMP+adap}$.}
\label{ex2-soln}
\end{figure}

\begin{figure}[thb]
\centering
\hbox{
\hspace{10mm}
\begin{minipage}[t]{3in}
\includegraphics[width=2.5in]{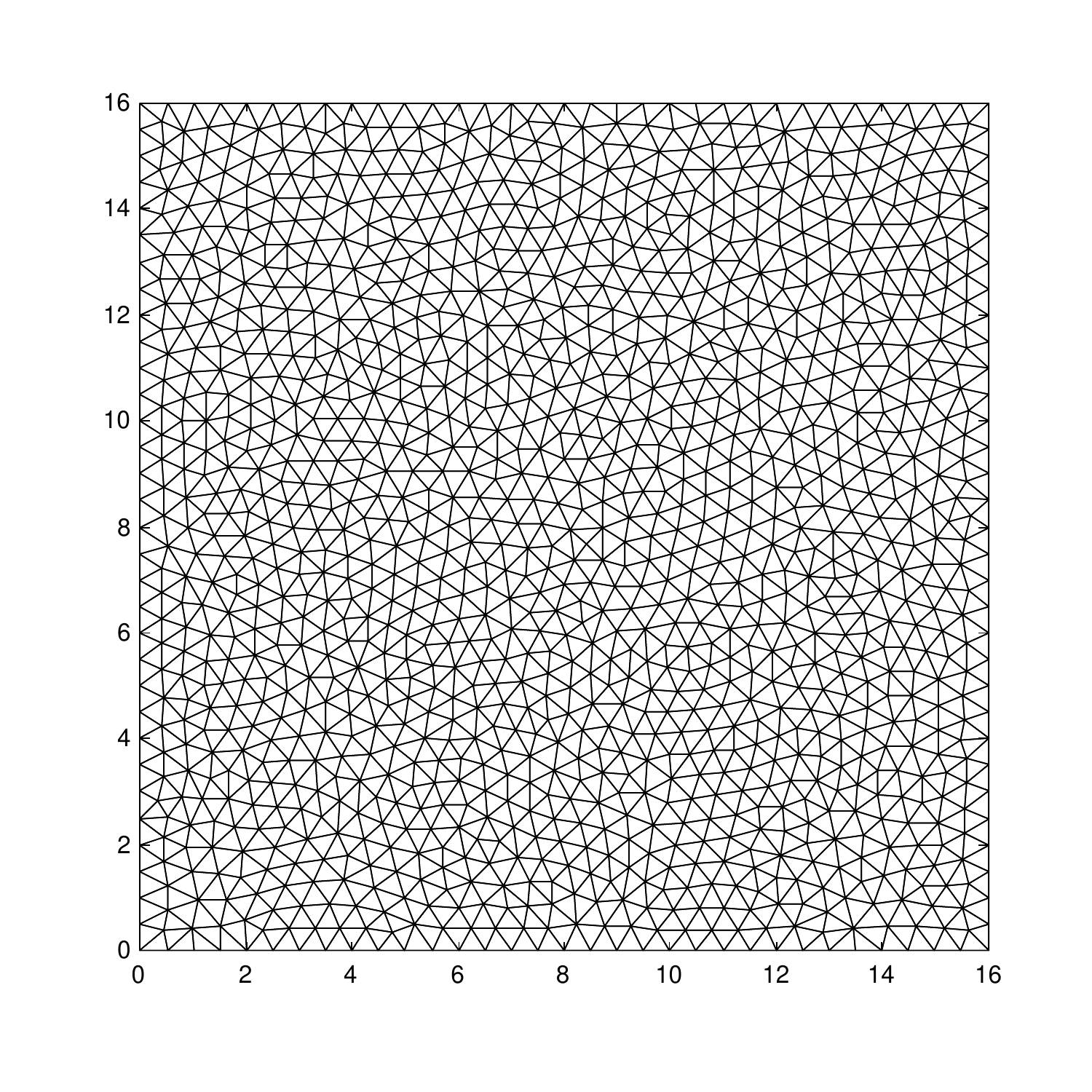}
\centerline{(a): $M_{unif}$, $N=2480$}
\end{minipage}
\hspace{10mm}
\begin{minipage}[t]{3in}
\includegraphics[width=2.5in]{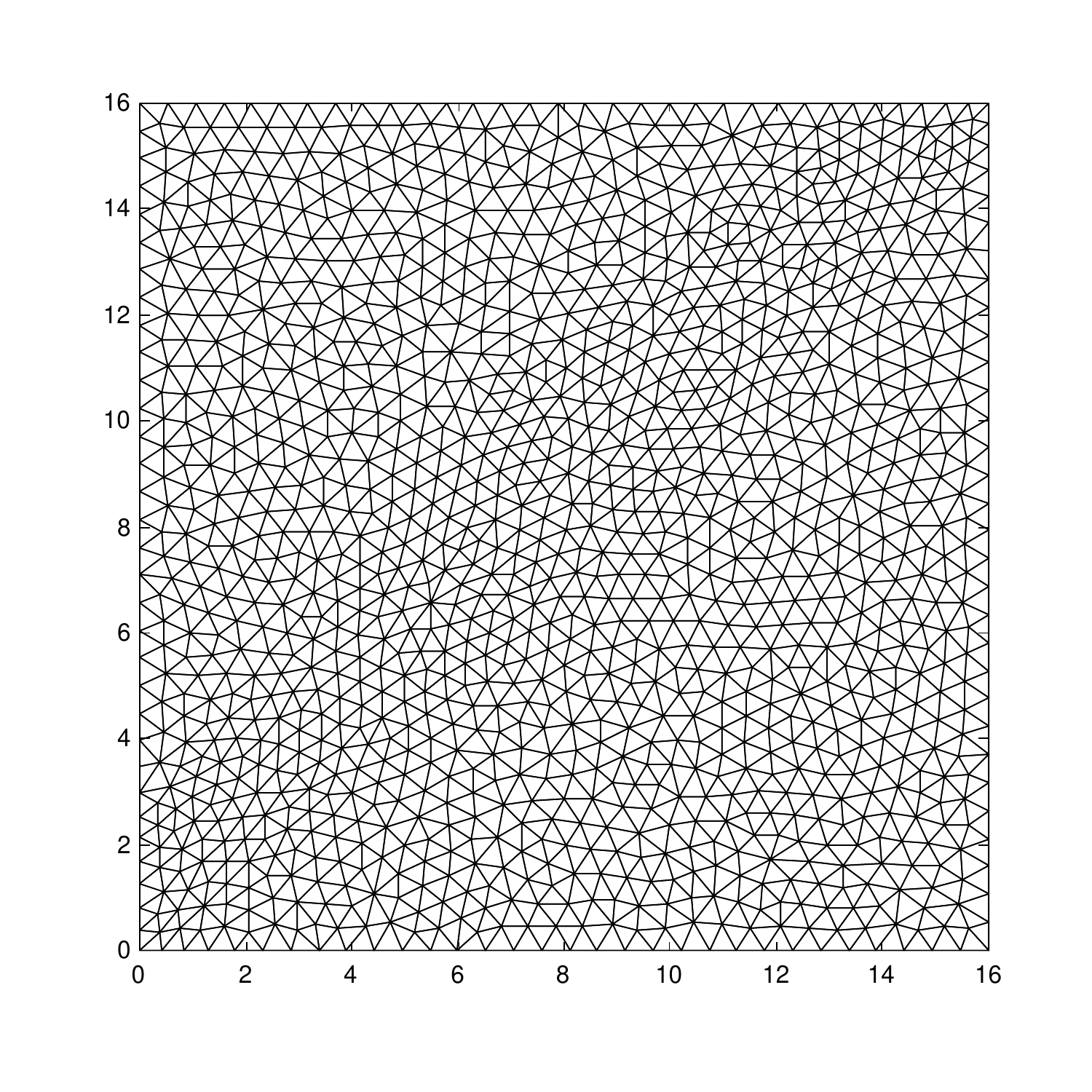}
\centerline{(b): $M_{adap}$, $N=2480$}
\end{minipage}
}
\hbox{
\hspace{10mm}
\begin{minipage}[t]{3in}
\includegraphics[width=2.5in]{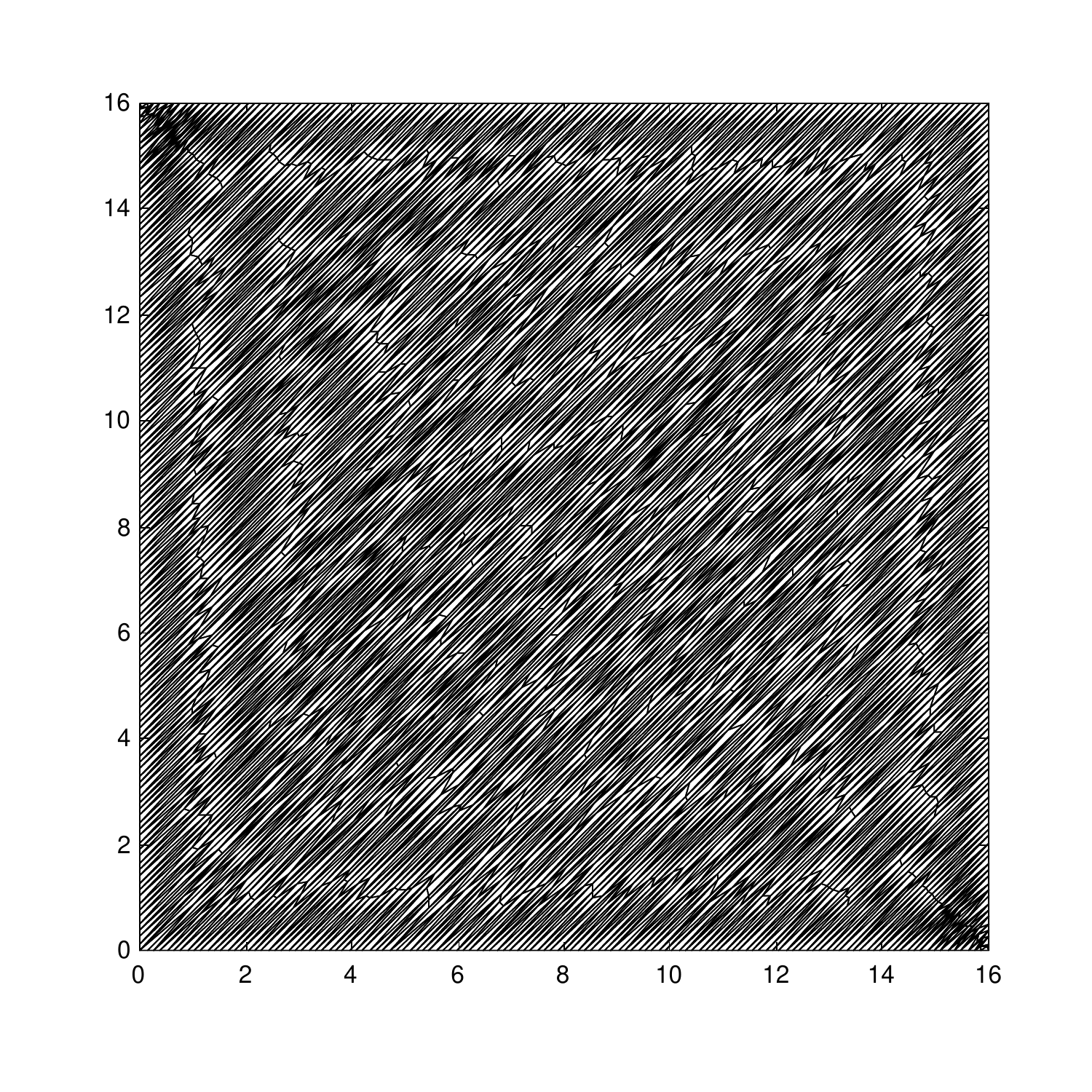}
\centerline{(c): $M_{DMP}$, $N=2754$}
\end{minipage}
\hspace{10mm}
\begin{minipage}[t]{3in}
\includegraphics[width=2.5in]{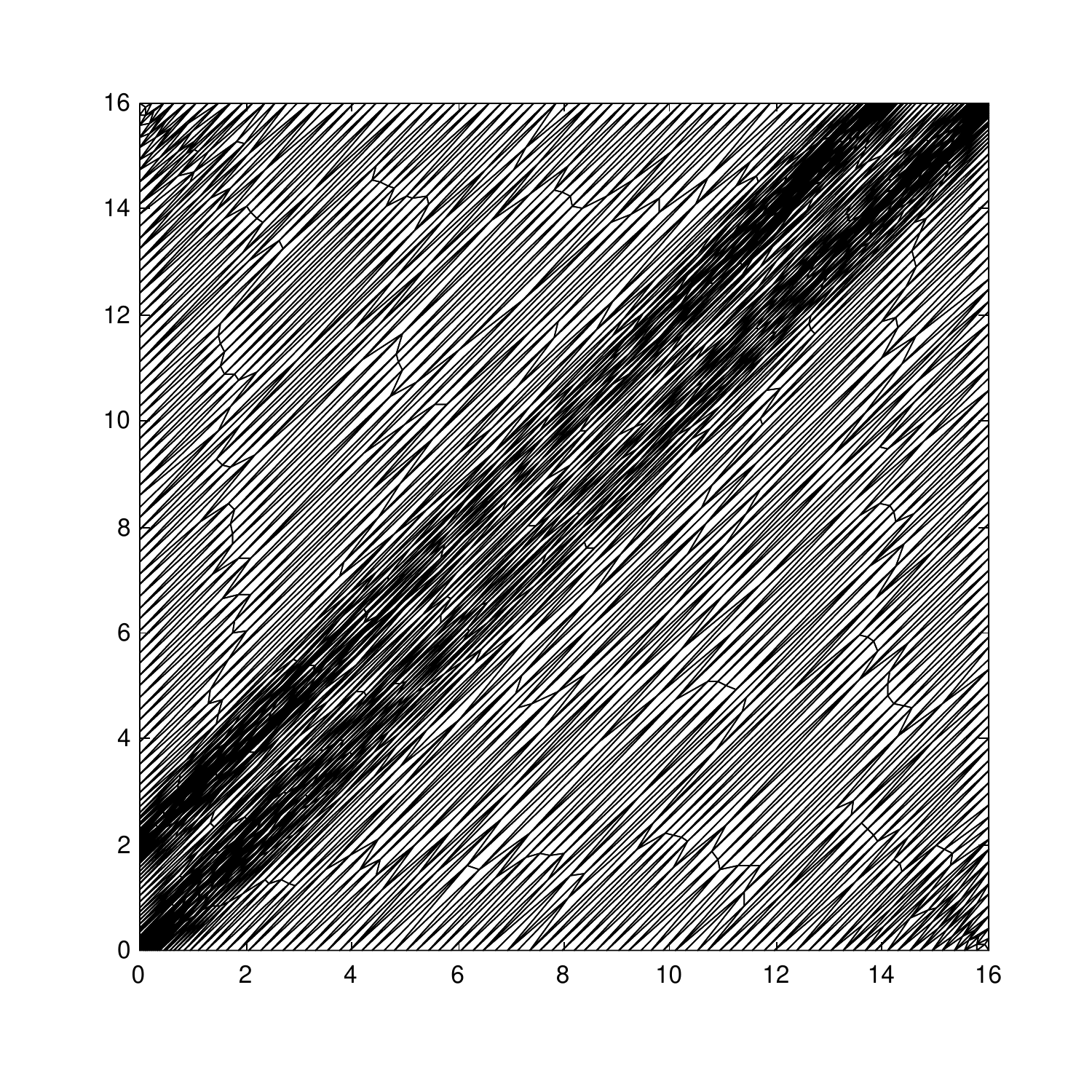}
\centerline{(d): $M_{DMP+adap}$, $N=2529$}
\end{minipage}
}
\caption{Example~\ref{ex2}. The adaptive meshes obtained with various metric tensors.}
\label{ex2-mesh}
\end{figure}

\begin{figure}[thb]
\centering
\hbox{
\hspace{10mm}
\begin{minipage}[t]{2.5in}
\includegraphics[width=2.5in]{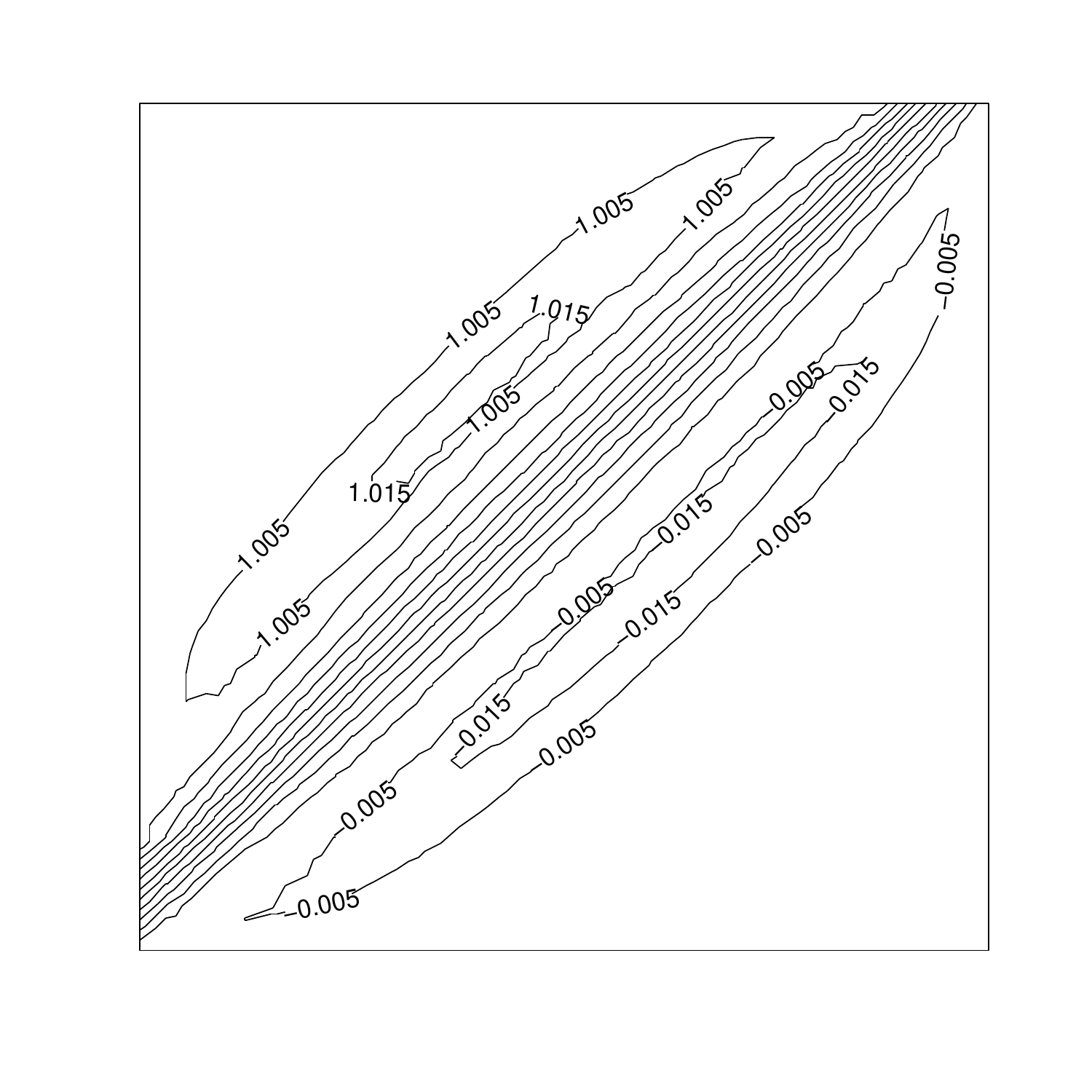}
\centerline{(a): $M_{unif}$, $u_{min}=-0.0280$, $u_{max}=1.0209$}
\end{minipage}
\hspace{10mm}
\begin{minipage}[t]{2.5in}
\includegraphics[width=2.5in]{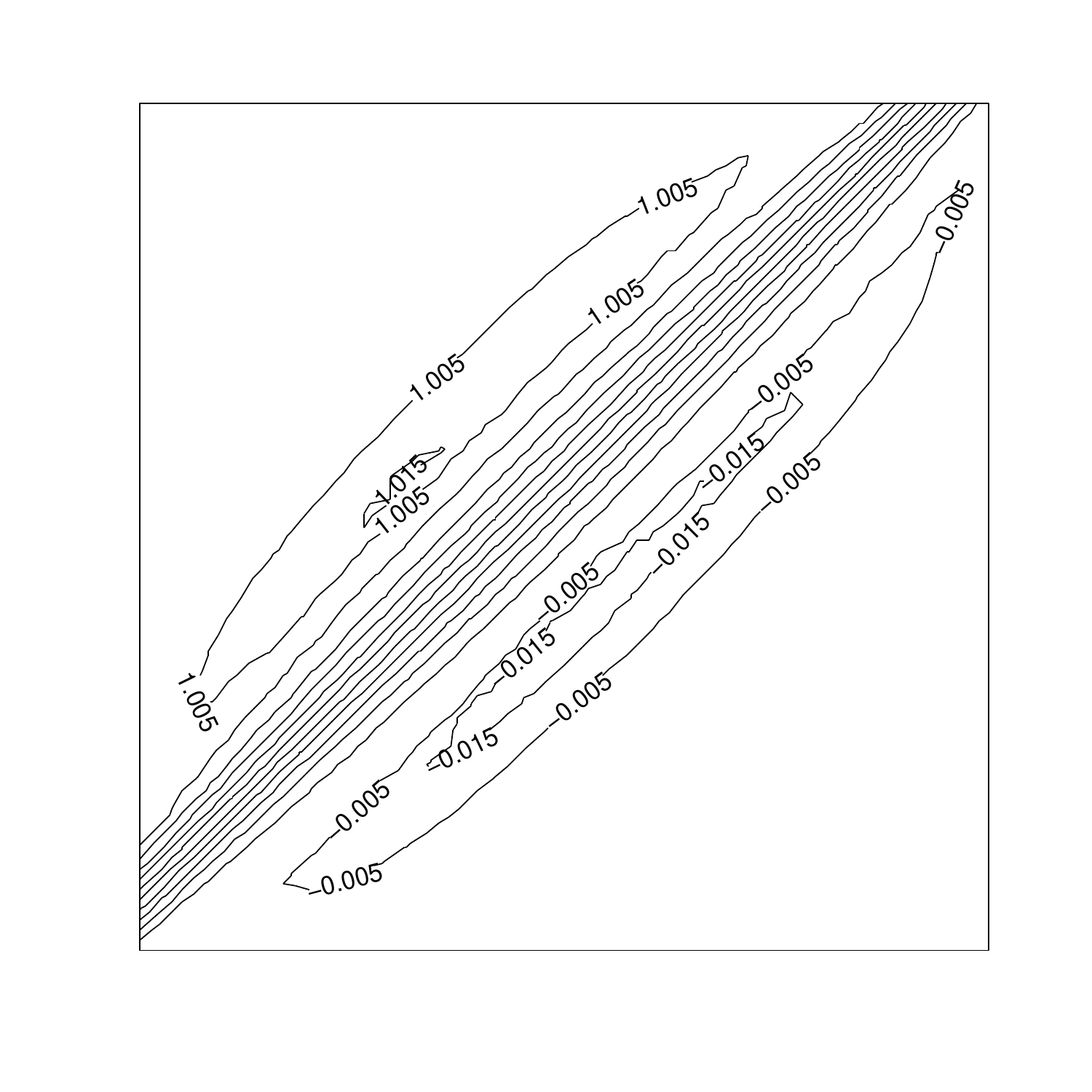}
\centerline{(b): $M_{adap}$, $u_{min}=-0.0195$, $u_{max}=1.0160$}
\end{minipage}
}
\hbox{
\hspace{10mm}
\begin{minipage}[t]{2.5in}
\includegraphics[width=2.5in]{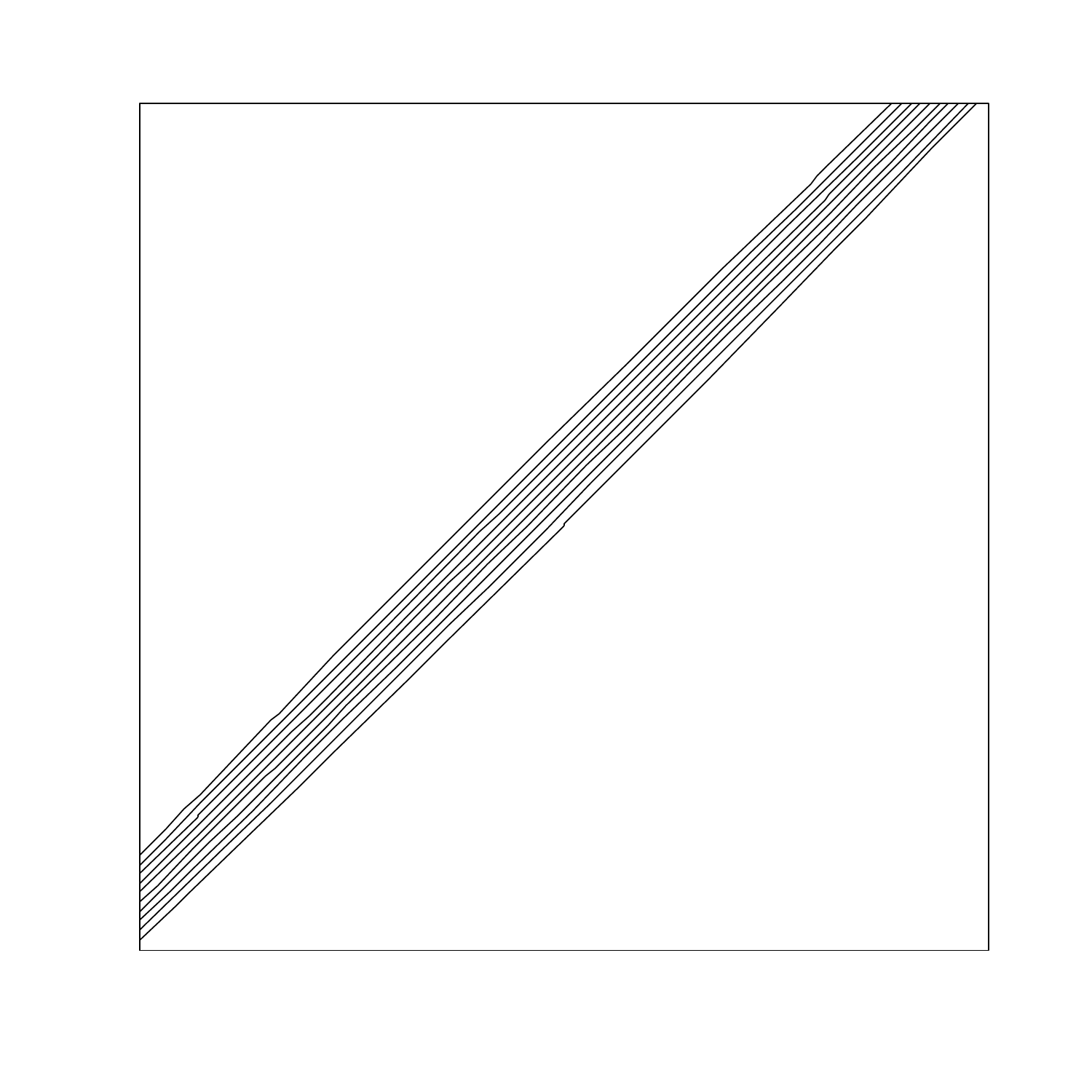}
\centerline{(c): $M_{DMP}$, $u_{min}=0$, $u_{max}=0$}
\end{minipage}
\hspace{10mm}
\begin{minipage}[t]{2.5in}
\includegraphics[width=2.5in]{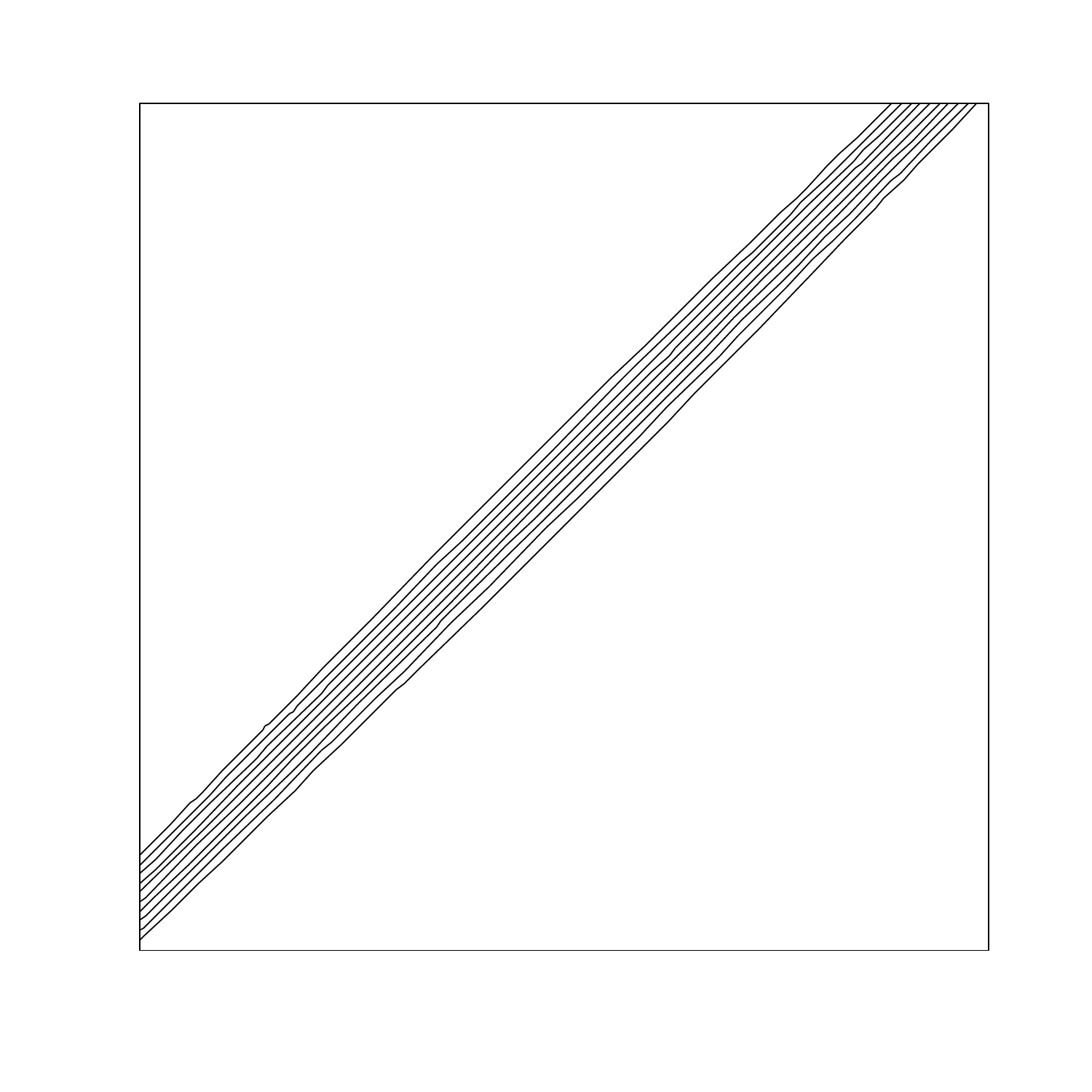}
\centerline{(d): $M_{DMP+adap}$, $u_{min}=0$, $u_{max}=0$}
\end{minipage}
}
\caption{Example~\ref{ex2}. Contours of the finite element solutions obtained with
different metric tensors.}
\label{ex2-contour}
\end{figure}

\begin{exam}
\label{ex2}
{\em
In this example, we consider BVP (\ref{bvp-pde}) and (\ref{bvp-bc}) with
\[
f \equiv 0,\quad g(x,0)=g(16,y)=0, 
\]
\beq
\nn 
g(0,y)= \left\{ \begin{array}{ll} 0.5y & \text{if  } 0 \le y < 2, \\
 	1 & \text{if  } 2 \le y \le 16, \end{array} \right. \text{ and  }
g(x,16)= \left\{ \begin{array}{ll} 1 & \text{if  } 0 \le x \le 14, \\
 	8-0.5x & \text{if  } 14 < x \le  16. \end{array} \right.
\eeq
The diffusion matrix is defined as
\beq
\nn
\mathbb{D}(x,y) = \left ( \begin{array}{cc} 500.5 & 499.5 \\ 499.5 & 500.5 \end{array} \right ).
\eeq
This is a simple example with a constant but anisotropic $\mathbb{D}$ and with 
a continuous boundary condition. It satisfies the maximum principle and its solution
stays between $0$ and $1$.

Numerical solutions, meshes, and solution contours are shown in Figs. \ref{ex2-soln}, \ref{ex2-mesh}, and
\ref{ex2-contour}, respectively. For this example, both undershoots and overshoots are observed
in the computed solutions with $M_{unif}$ and $M_{adap}$ but not with 
with $M_{DMP}$ and $M_{DMP+adap}$. 
This example demonstrates that a scheme violating DMP can produce unphysical extrema even for a simple
problem with constant diffusion, continuous boundary conditions, and a convex domain.
}\end{exam}

\begin{figure}[thb]
\centering
\hbox{
\begin{minipage}[t]{3in}
\includegraphics[width=2.5in]{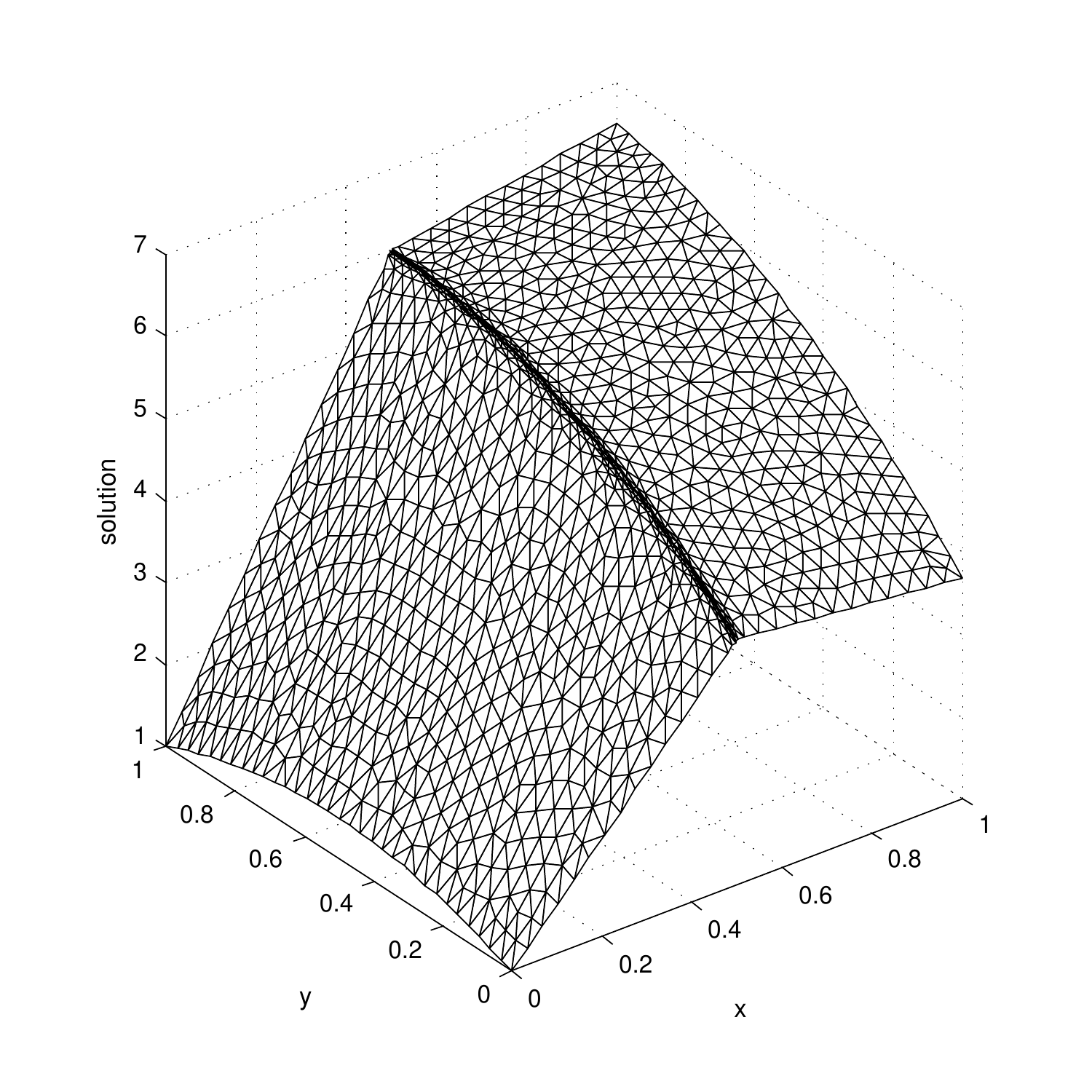}
\centerline{(a): $M_{adap}$, numerical solution, $u_{min}=0$}
\end{minipage}
\hspace{10mm}
\begin{minipage}[t]{3in}
\includegraphics[width=2.5in]{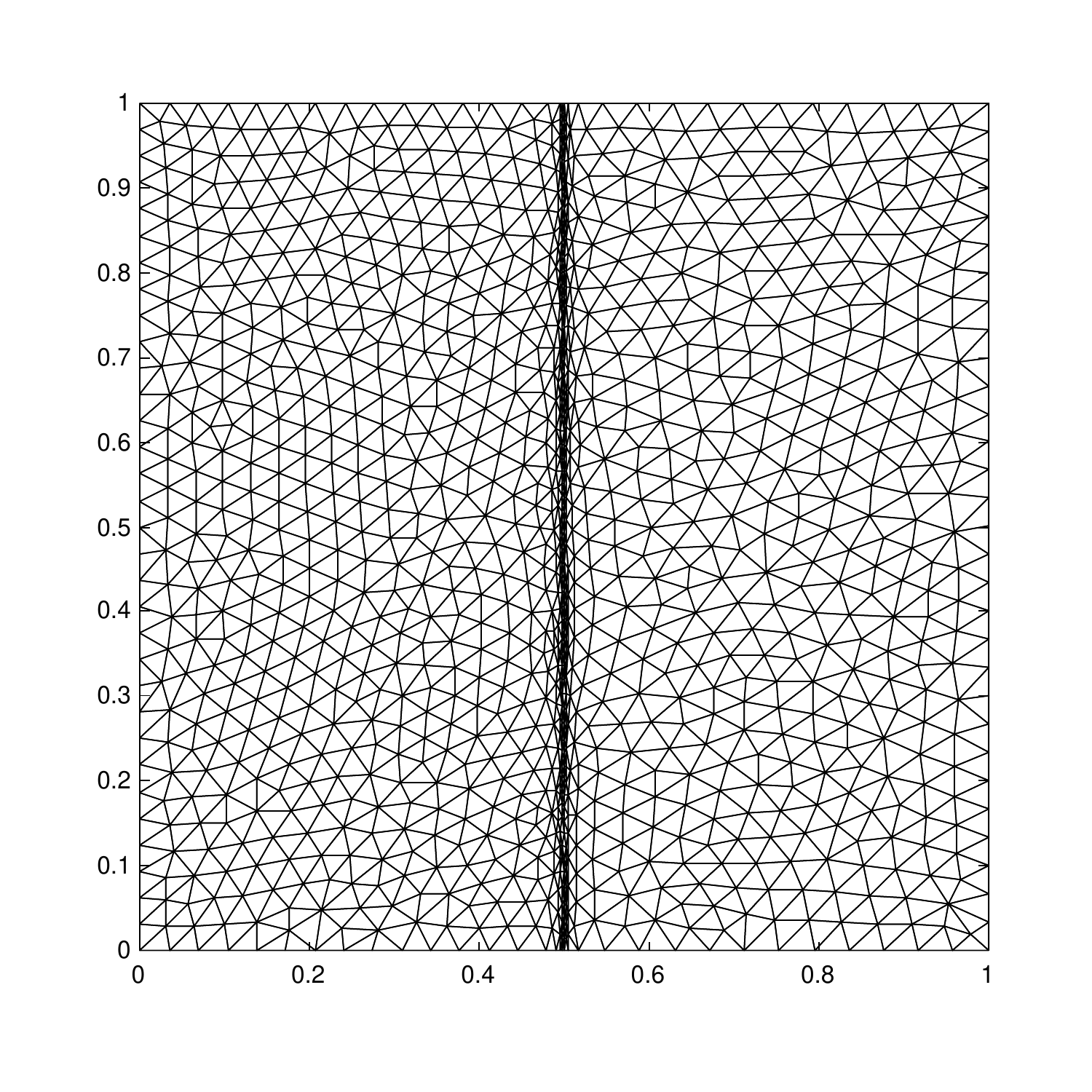}
\centerline{(b): $M_{adap}$, mesh, $N=2362$}
\end{minipage}
}
\hbox{
\begin{minipage}[t]{3in}
\includegraphics[width=2.5in]{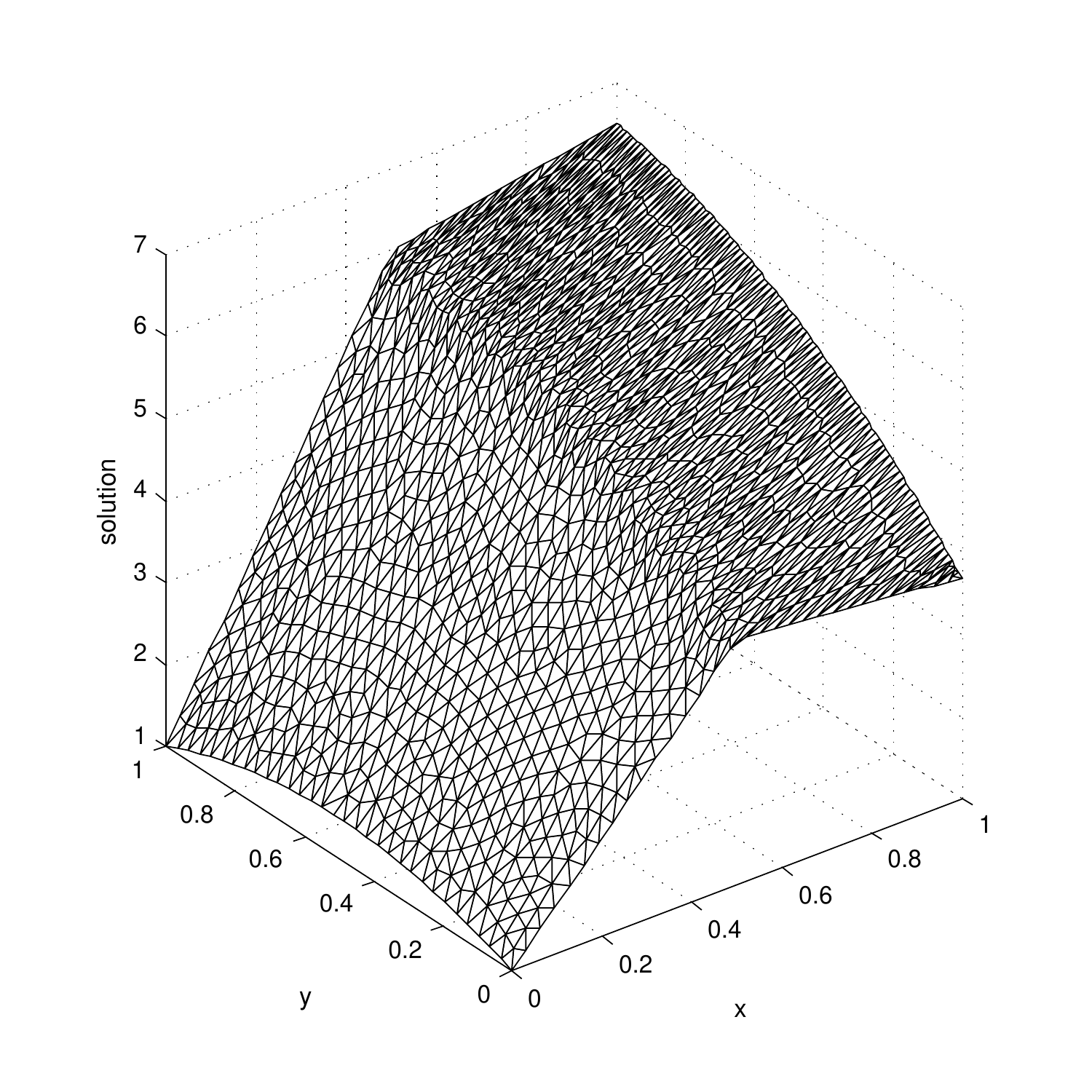}
\centerline{(c): $M_{DMP}$, numerical solution, $u_{min}=0$}
\end{minipage}
\hspace{10mm}
\begin{minipage}[t]{3in}
\includegraphics[width=2.5in]{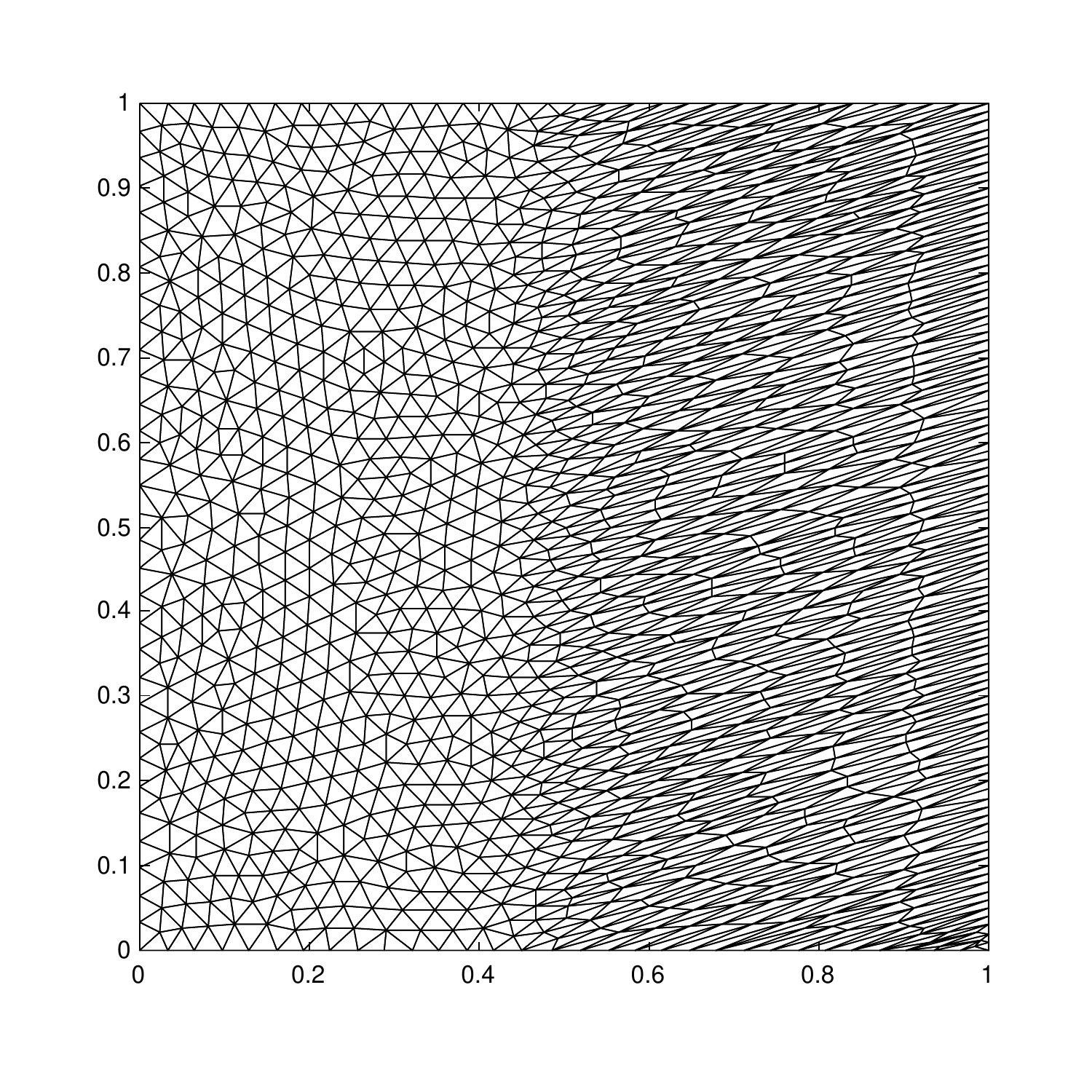}
\centerline{(d): $M_{DMP}$, mesh, $N=2415$}
\end{minipage}
}
\hbox{
\begin{minipage}[t]{3in}
\includegraphics[width=2.5in]{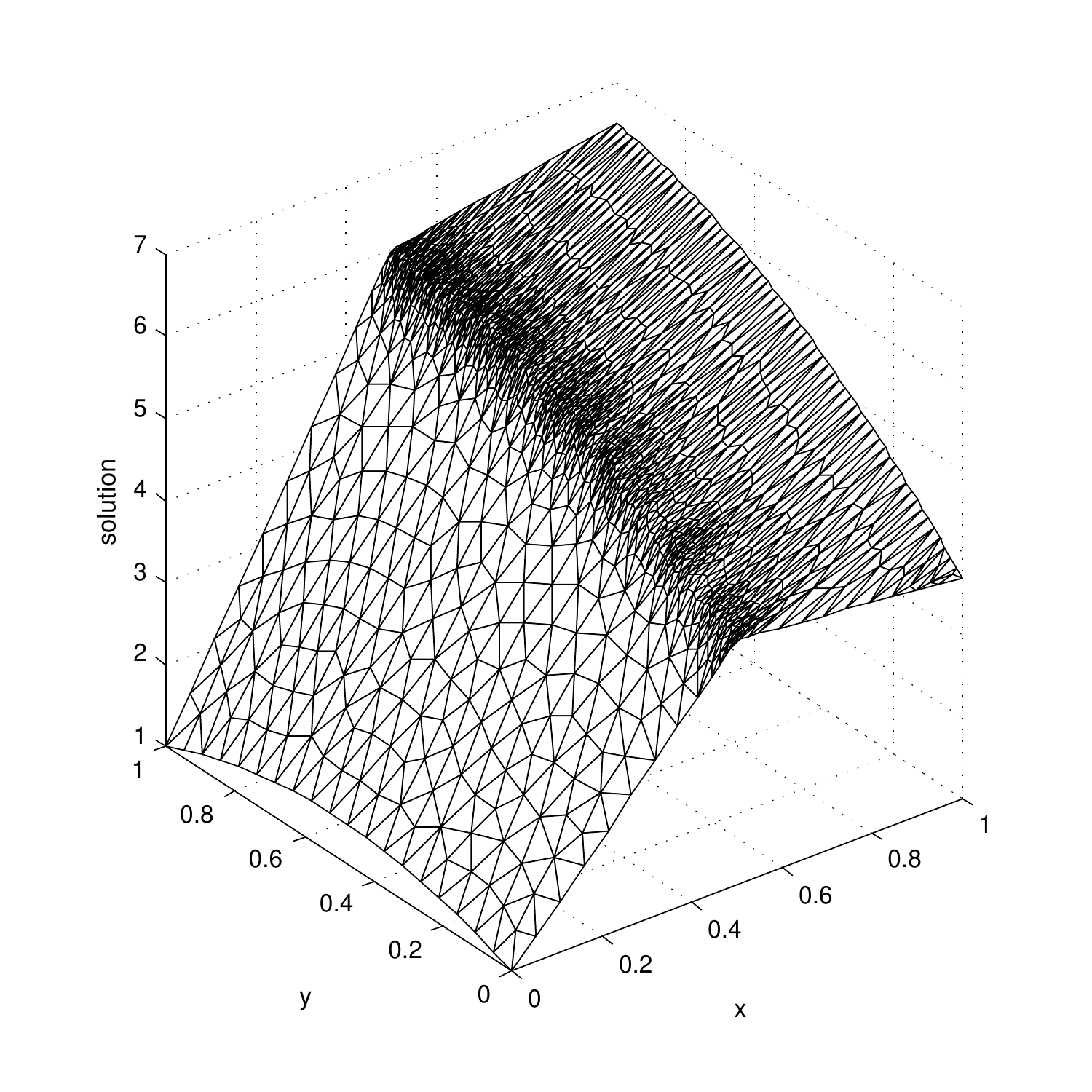}
\centerline{(e): $M_{DMP+adap}$, numerical solution, $u_{min}=0$}
\end{minipage}
\hspace{10mm}
\begin{minipage}[t]{3in}
\includegraphics[width=2.5in]{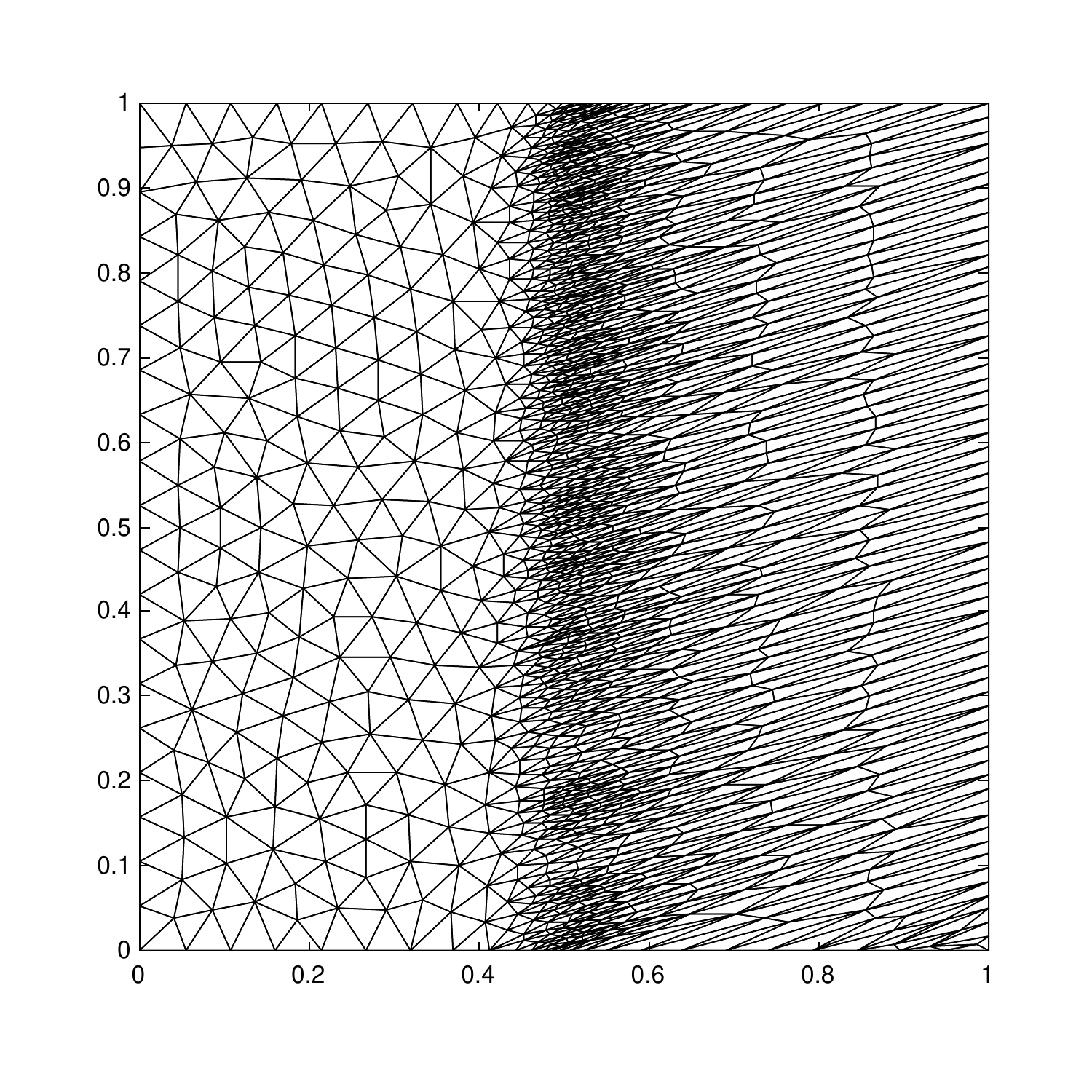}
\centerline{(f): $M_{DMP+adap}$, mesh, $N=2490$}
\end{minipage}
}
\caption{Example~\ref{ex3}. Numerical solutions and meshes obtained with three metric tensors.}
\label{ex3-solnmesh}
\end{figure}

\begin{figure}[thb]
\centering
\hbox{
\begin{minipage}[t]{3in}
\includegraphics[width=3in]{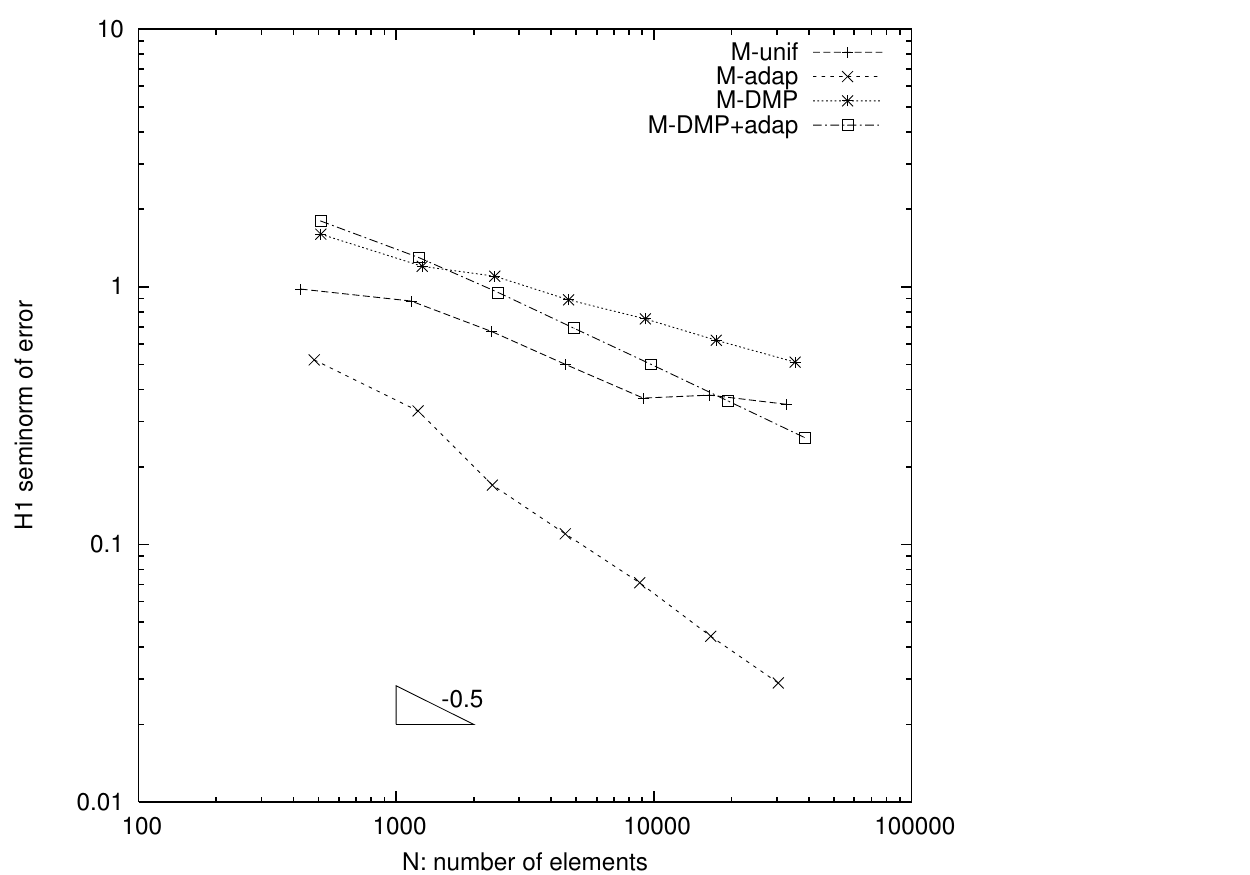}
\centerline{(a): $H^1$ semi-norm of error}
\end{minipage}
\hspace{5mm}
\begin{minipage}[t]{3in}
\includegraphics[width=3in]{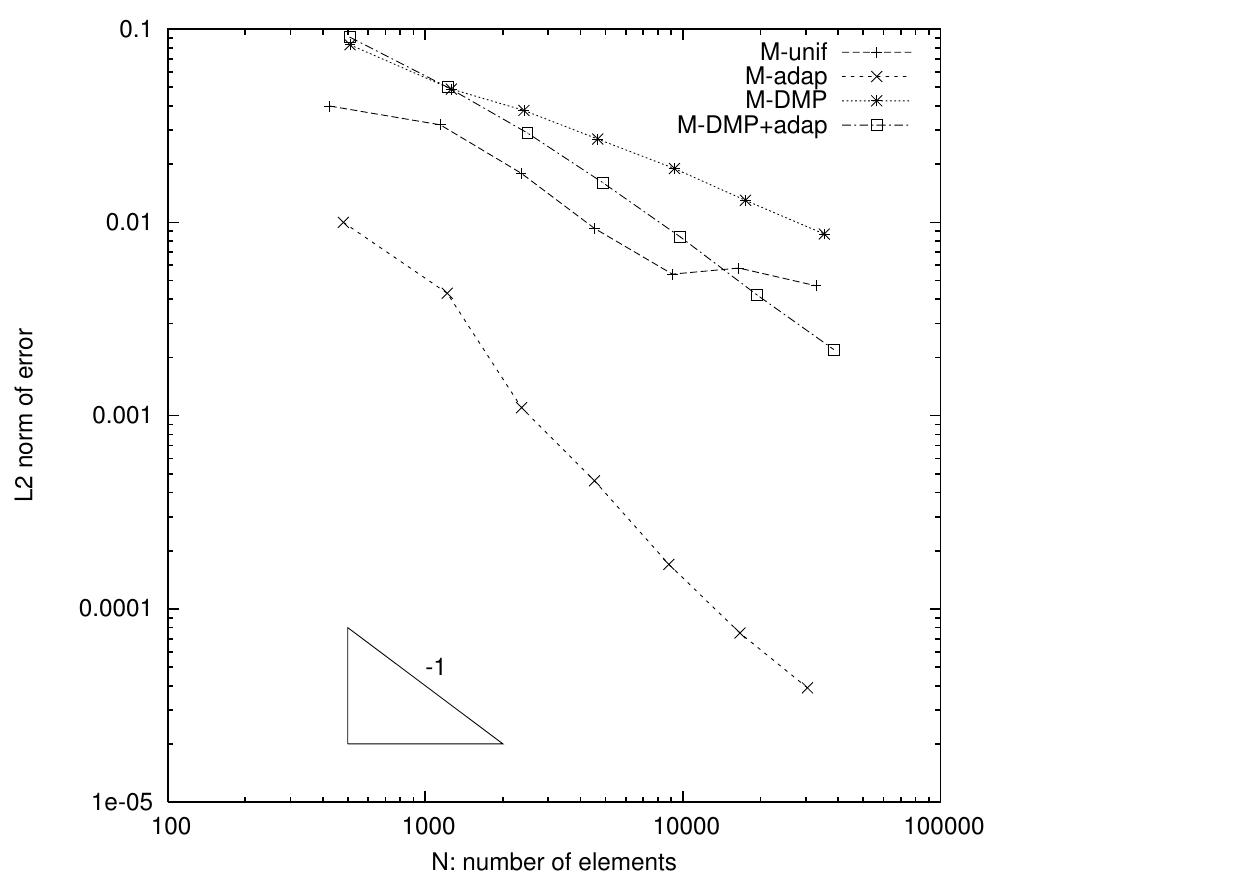}
\centerline{(b): $L^2$ norm of error}
\end{minipage}
}
\caption{Example~\ref{ex3}. The $H^1$ semi-norm and $L^2$ norm of solution error
are shown as functions of the number of elements for metric tensors $M_{unif}$,
$M_{adap}$, $M_{DMP}$, and $M_{DMP+adap}$.}
\label{ex3-error}
\end{figure}

\begin{figure}[thb]
\centering
\hbox{
\begin{minipage}[t]{3in}
\includegraphics[width=3in]{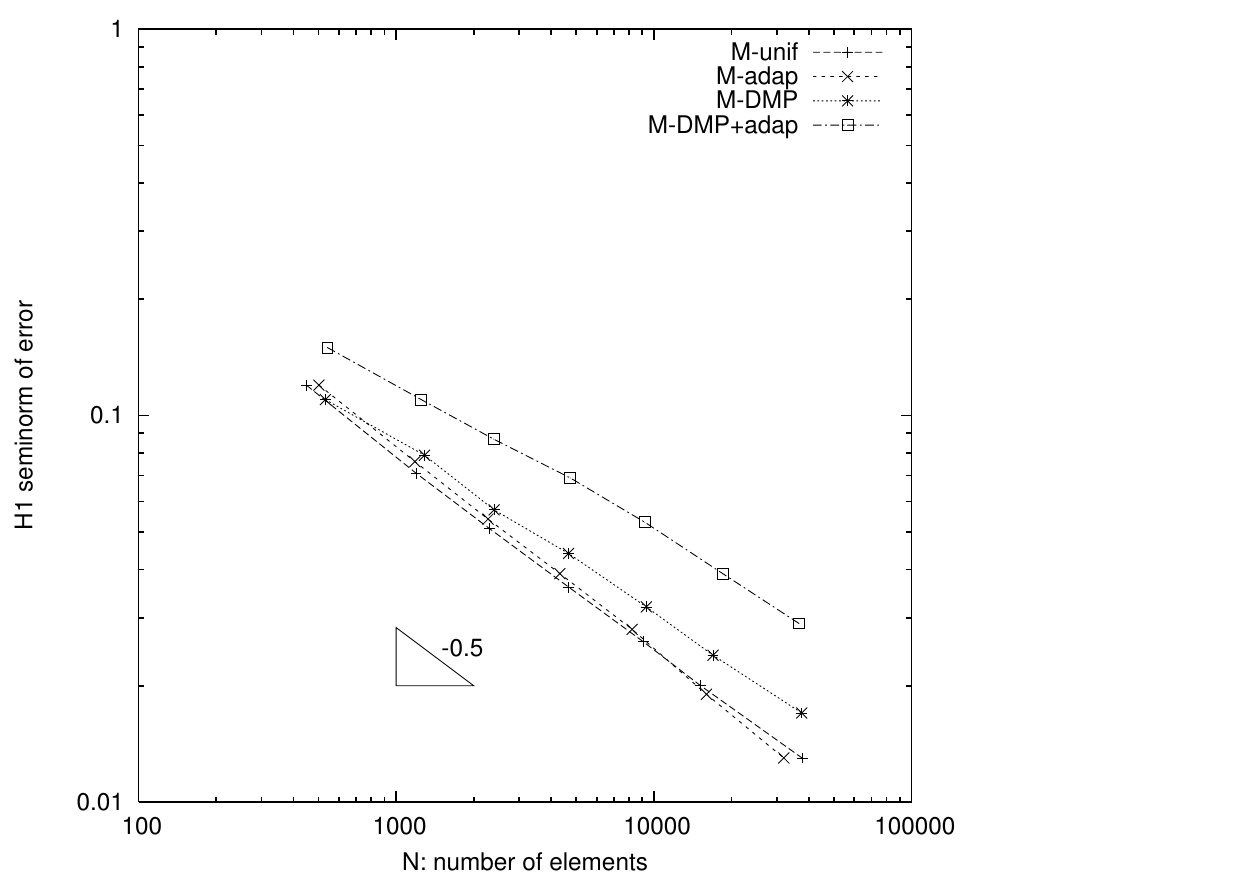}
\centerline{(a): $H^1$ semi-norm of error}
\end{minipage}
\hspace{5mm}
\begin{minipage}[t]{3in}
\includegraphics[width=3in]{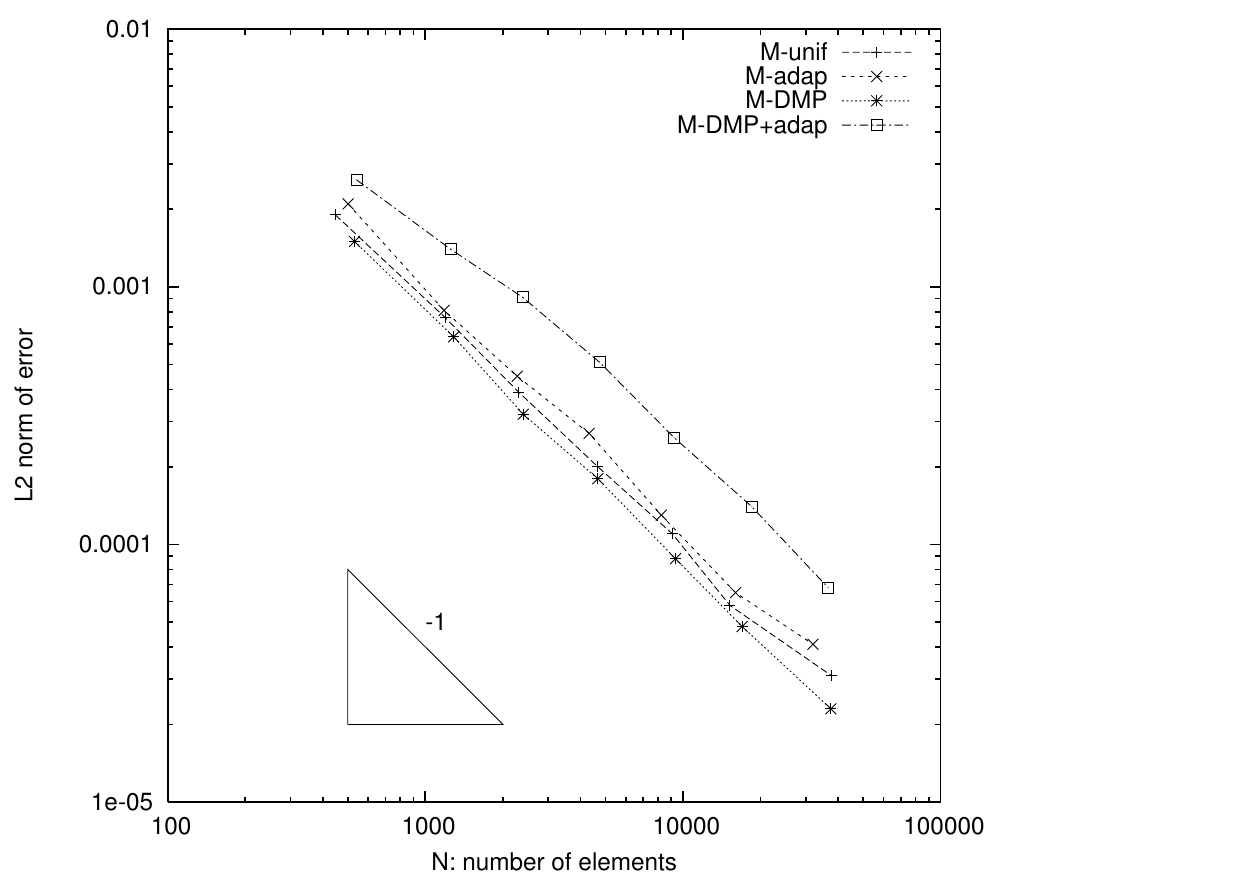}
\centerline{(b): $L^2$ norm of error}
\end{minipage}
}
\caption{Example~\ref{ex3}. The $H^1$ semi-norm and $L^2$ norm of solution error
are shown as functions of the number of elements for metric tensors $M_{unif}$,
$M_{adap}$, $M_{DMP}$, and $M_{DMP+adap}$. The interface ($x=0.5$) is predefined in the mesh.}
\label{ex3b-error}
\end{figure}

\begin{exam}
\label{ex3}
{\em
This example is given by (\ref{bvp-pde}) and (\ref{bvp-bc}) with
\[
\Omega =(0,1)\times (0,1), \quad
f(x,y) = \left \{ \begin{array}{rl}
4.0, & \text{ if } x<0.5 \\ -5.6, & \text{ if } x>0.5
\end{array} \right. ,
\quad u = u_{exact} \mbox{ on } \p \Omega ,
\]
\[
\mathbb{D}(x,y) = \left \{ \begin{array}{ll}
D_1, & \text{ if } x<0.5, \\ D_2, & \text{ if } x>0.5,
\end{array} \right. \quad
D_1 = \left ( \begin{array}{cc} 1 & 0 \\ 0 & 1 \end{array} \right ), \quad
D_2 = \left ( \begin{array}{cc} 10 & 3 \\ 3 & 1 \end{array} \right ) .
\]
The problem has the exact solution
\beq
u(x,y) = \left \{ \begin{array}{ll}
1 - 2y^2 + 4xy + 2y + 6x, & \text{ if } x \le 0.5 \\
-2 y^2 + 1.6 xy - 0.6 x + 3.2 y + 4.3, & \text{ if } x>0.5 .
\end{array} \right.
\label{ex3-exact-soln}
\eeq
Note that the value and primary diffusion direction of the diffusion matrix
change across the line $x=0.5$. This example has been studied in \cite{KSS09}. 

Solutions and meshes obtained with various metric tensors are shown in Fig. \ref{ex3-solnmesh}. For this example,
no overshoots and undershoots are observed for all numerical solutions.
The meshes obtained with $M_{DMP}$ and $M_{DMP+adap}$
show a better alignment with the primary diffusion direction than that obtained with $M_{adap}$.
Moreover, elements are concentrated along the line $x=0.5$ for the meshes obtained with
$M_{adap}$ and $M_{DMP+adap}$ whereas there is no concentration in the mesh
shown in Fig. \ref{ex3-solnmesh}(d) for $M_{DMP}$. The results are consistent with
what is expected from the construction of the metric tensors.

The exact solution is available for this example. 
The $H^1$ semi-norm and $L^2$ norm of the error are shown in Fig. \ref{ex3-error} as functions
of the number of mesh elements. Metric tensor $M_{adap}$ leads to far more accurate
results than the other three metric tensors, which produce comparable results for
the considered range of $N$. Moreover, $M_{adap}$ and $M_{DMP+adap}$ give the same convergence
rate, i.e., $|e^h|_{H^1(\Omega)} = O( N^{-0.5})$ and $\|e^h\|_{L^2(\Omega)} = O( N^{-1})$,
while $M_{unif}$ and $M_{DMP}$ result in a slower convergence rate, $|e^h|_{H^1(\Omega)}
= O( N^{-0.25})$ and $\|e^h\|_{L^2(\Omega)} = O( N^{-0.5})$.
This demonstrates the advantage of using adaptive meshes.
Interestingly, the results in \cite{KSS09} (Table 4) obtained
with a slope-limited scheme for triangular meshes also show a similar slow convergence.

It should be pointed out that the above results have been obtained when the interface ($x=0.5$) is not
predefined in the mesh. If the interface is predefined in the mesh, then the solution (\ref{ex3-exact-soln})
can be approximated accurately in the linear finite element space. As shown in Fig. \ref{ex3b-error},
all metric tensors produce comparable solutions and the same convergence rate
$|e^h|_{H^1(\Omega)} = O( N^{-0.5})$ and $\|e^h\|_{L^2(\Omega)} = O( N^{-1})$.
}\end{exam}

\section{Conclusions and comments}
\label{Sec-con}

In the previous sections we have developed a mesh condition (\ref{g-nonobtuse})
under which the linear finite element approximation of anisotropic diffusion problem
(\ref{bvp-pde}) and (\ref{bvp-bc}) validates the discrete counterpart of the maximum principle
satisfied by the continuous problem.
The condition is a generalization of the well known non-obtuse angle condition developed
for isotropic diffusion problems and requires that the dihedral angles of mesh elements
measured in a metric depending only on the diffusion matrix be non-obtuse.

We have also developed two variants of the anisotropic non-obtuse angle condition, (\ref{g-nonobtuse-2}) and
(\ref{g-nonobtuse-3}), which can be more convenient to use in actual mesh generation.
Indeed, metric tensor (\ref{M-DMP}) for use in
anisotropic mesh generation is derived based on (\ref{g-nonobtuse-3}) for accounting for
DMP satisfaction.  Moreover, an optimal metric tensor (\ref{M-DMP+adap}) accounting for
both DMP satisfaction and mesh adaptation is obtained
from (\ref{g-nonobtuse-3}) by minimizing an interpolation error bound.
Features of these metric tensors are illustrated in numerical examples.

It is worth pointing out that condition (\ref{g-nonobtuse}) has been derived based
on the local stiffness matrix on a mesh element. Like the non-obtuse angle condition for
isotropic diffusion problems, (\ref{g-nonobtuse}) may be relaxed by considering
the global stiffness matrix as a whole \cite{Let92}. Moreover,
we have restricted our attention to linear PDE (\ref{bvp-pde}) and
Dirichlet boundary condition (\ref{bvp-bc}). But the procedure developed in this work can be extended to
problems with nonlinear diffusion  $\mathbb{D} = \mathbb{D}(\V{x}, u, \nabla u)$ and
mixed boundary conditions (e.g., see \cite{KK05,KK06,KKK07,KL95}) without major modification.

Although the numerical examples have been presented in 2D, the anisotropic non-obtuse angle 
condition (\ref{g-nonobtuse}) and the corresponding metric tensor formulas (\ref{M-DMP}),
(\ref{M-DMP+adap}), and (\ref{M-DMP-2}) are $d$-dimensional ($d=1,2,3$). In 3D, a Delaunay triangulation
may not guarantee the satisfaction of DMP \cite{Let92}. Nevertheless, some polyhedrons can be
decomposed into tetrahedra satisfying the non-obtuse angle condition (\ref{non-obtuse})
and therefore the numerical solution satisfies DMP; e.g., see \cite{KL95}.
It is expected that this will also work for the anisotropic non-obtuse angle condition (\ref{g-nonobtuse})
for a given metric tensor $M$. On the other hand, the existence of the decomposition of an arbitrary polyhedron into non-obtuse tetrahedra is an open problem \cite{KL95}. It is also unclear if a 3D triangulation can be generated
to (approximately) satisfy the $M$-uniform mesh conditions (\ref{cond-equi} and \ref{cond-align}).
Those are interesting topics to investigate in the future. 

\vspace{20pt}

{\bf Acknowledgment.} The work was supported in part by the National Science Foundation (USA)
under grant DMS-0712935.


\end{document}